\documentclass[12pt,twoside]{report}
\usepackage{suthesis-2e}

\usepackage{amsfonts}
\usepackage{amsmath}
\usepackage{amssymb}
\usepackage{amsthm}
\usepackage{graphicx}

\usepackage{stdmath_thesis}

\graphicspath{{graphics_included/}}

    \title{Stability and Control of Functional Differential Equations}
        \author{Matthew Monnig Peet}
\dept{Aeronautics and Astronautics} % default is Computer Science, uncomment for other departments
    \principaladviser{Sanjay Lall}
    \firstreader{Geir Dullerud}
    \secondreader{Stephen Rock}
    \thirdreader{Matthew West}

\begin{document}
% first the preface sections.
% this includes the file preface.tex which should include the
% following commands
\beforepreface
% \prefacesection{preface}
% body of the preface
\prefacesection{Abstract}
This thesis addresses the question of stability of systems defined
by differential equations which contain nonlinearity and delay. In
particular, we analyze the stability of a well-known delayed
nonlinear implementation of a certain Internet congestion control
protocol. We also describe a generalized methodology for proving
stability of time-delay systems through the use of semidefinite
programming.

In Chapters~\ref{chp:ICC} and~\ref{chp:IQCresult}, we consider an
Internet congestion control protocol based on the decentralized
gradient projection algorithm. For a certain class of utility
function, this algorithm was shown to be globally convergent for
some sufficiently small value of a gain parameter. Later work gave
an explicit bound on this gain for a linearized version of the
system. This thesis proves that this bound also implies stability of
the original system. The proof is constructed within a generalized
passivity framework. The dynamics of the system are separated into a
linear, delayed component and a system defined by a nonlinear
differential equation with discontinuity in the dynamics.
Frequency-domain analysis is performed on the linear component and
time-domain analysis is performed on the nonlinear discontinuous
system.

In Chapter~\ref{chp:LinearCase}, we describe a general methodology
for proving stability of linear time-delay systems by computing
solutions to an operator-theoretic version of the Lyapunov
inequality via semidefinite programming. The result is stated in
terms of a nested sequence of sufficient conditions which are of
increasing accuracy. This approach is generalized to the case of
parametric uncertainty by considering parameter-dependent Lyapunov
functionals. Numerical examples are given to demonstrate convergence
of the algorithm. In Chapter~\ref{chp:nonlinearcase}, this approach
is generalized to nonlinear time-delay systems through the use of
non-quadratic Lyapunov functionals.

\prefacesection{Preface}

This thesis covers my work done in the Networked Systems and
Controls Lab at Stanford University under the direction of Professor
Sanjay Lall during the period from August of 2002 through January of
2006. With the exception of Chapter~\ref{chp:nonlinearcase}, almost
all of the work presented here has been previously published or
submitted to peer-reviewed academic journals or conferences. The
chapters of this thesis can be divided into those containing
original research and those containing an overview of existing
results. Specifically, Chapters~\ref{chp:IQCresult},
\ref{chp:LinearCase} and~\ref{chp:nonlinearcase} contain original
results, while most of the rest of the chapters can be considered as
survey material.

%The general topic of my doctoral studies has been algorithms for
%stability and control of nonlinear decentralized systems with delay.
%The application which has received the most attention during my
%studies has been analysis of protocols for internet congestion
%control.

Some attempt has been made in this thesis to keep the length of the
chapters minimal in order to improve readability. The presentation
style is mathematically oriented, with specifics of numerical
implementation suppressed. Proofs, when especially long, have been
moved to the appendices. To maintain focus on the theoretical
contribution, in-depth treatment of various special cases has been
omitted when such treatment can be inferred from previous
exposition.

%\newpage
%
%However everything should be taken with a grain of salt.
%
%    \prefacesection{Acknowledgements}
%        I would like to thank my mother and the little green men from
%        Mars.

% any other preface sections

% the last preface section (e.g., acknowledgement.tex)
% should look like
% \prefacesection{Acknowledgement}
% body
    \prefacesection{Acknowledgements}
There are many people who have contributed to the work contained in
this thesis. In particular, I would like to thank my advisor, Sanjay
Lall. Throughout my time at Stanford, Sanjay has always striven to
provide guidance in a field which often seems chaotic. In addition,
he has attempted to instill in me the discipline necessary to
produce results which are capable of withstanding outside scrutiny.
I would also like to thank my outstanding thesis committee
consisting of Geir Dullerud, G\"{u}nter Niemeyer, Steve Rock and
Matthew West. Their insightful comments have helped shape much of
the contents of this thesis. I would like to thank my Lab group,
consisting of Been-Der Chen, Randy Cogill, Ritesh Madan, Mike
Rotkowitz, and Ta-Chung Wang. These folks are the best of the best
and it has been a pleasure working with them. I would also like to
thank my wife, Yulia, who works in Computational Fluid Dynamics and
has provided support at home as well as in my research. Finally, I
would like to thank all my friends here at Stanford for helping to
make the life of a Ph.D. student bearable.
    %\afterpreface

\afterpreface

% now for the body of the thesis, modify the number of these lines as needed

% this includes chapter1.tex which should start with a \chapter{...}
% command
\chapter{Introduction}

\section{Research Goals}

The analysis of nonlinear systems is a subject of much recent
interest. Although linearization is a well-established approach to
analysis, results obtained in this manner are always, at best,
local. When considering changes in critical systems such as the
Internet, the guarantee of convergence associated with a global
stability result carries significant weight. Furthermore, some
dynamical systems, such as are found in biology or high performance
aircraft, are dominated by nonlinear behavior. In such cases, the
practical value of linear analysis is limited. Apart from
nonlinearity, the study of decentralized systems and systems with
delay is also the topic of much active research. Communications
systems, an issue of significant current interest, are often
decentralized and, by virtue of their geographic reach, inevitably
contain delay. Earth based telescopes are currently being built
which consist of thousands of mirror segments; each individually
actuated but required to move in a coordinated manner. Such systems
are difficult to model using a state space of reasonable dimension.
\\

The recent development of efficient algorithms for optimization on
the cone of positive semidefinite matrices has created a fundamental
shift in how research is conducted. Specifically, semidefinite
programming has been widely adopted in control theory and has led to
a greater understanding of linear finite dimensional centralized
systems. Unfortunately, this same understanding has not yet been
extended to nonlinear, infinite dimensional or decentralized
systems. In fact, many critical questions in analysis and control of
these systems have been shown to be NP-hard. The goal of our
research at Stanford has been to find new ways to address these
types of problems. One approach is to expand the definition of
solution. Instead of giving a single necessary and sufficient
condition, expressible as a semidefinite program, one can construct
a nested sequence of sufficient conditions, of increasing accuracy,
which converge to a necessary and sufficient condition and are
expressible as semidefinite programs. Examples of this technique as
well as others can be found in the chapters of this thesis.\\

\paragraph{Internet Congestion Control}

The application which has motivated much of the research in this
thesis has been stability analysis of proposed protocols for
Internet congestion control. The analysis of congestion control
protocols has received much attention recently. This work has been
motivated by concern about the ability of current protocols to
ensure stability and performance as the number of users and amount
of bandwidth continues to increase. Although the protocols that have
been used in the past have performed well as the Internet has
increased in size, as capacities and delays increase instability
will become a problem. In our work, we have considered global
stability with delay of congestion control protocols which attempt
to solve a distributed network optimization problem. These systems
are described by differential equations with delay and contain
non-static nonlinearity. Because of the nonlinearity and delays,
proving convergence is difficult. By combining frequency and time
domain techniques using a generalized passivity framework, in
Chapters~\ref{chp:ICC} and~\ref{chp:IQCresult}, we have shown for
certain protocols that global stability with delay holds under the
same conditions as local stability. Condensed versions of these
results also appear in publication in~\cite{peet_2004b}
and~\cite{peet_2006a}. These tight bounds, verified by experimental
evidence, allow one to accurately predict when congestion control
will fail.\\

\paragraph{An Approach to Analysis and Synthesis}

Our research at Stanford has produced a number of new results and
tools concerning the analysis of systems with delay, nonlinearity
and decentralized structure. Specifically, in
Chapter~\ref{chp:LinearCase}, we have proposed two new types of
refutation which can be used to construct the positive quadratic
Lyapunov-Krasovskii functionals necessary for stability of linear
time-delay systems. Furthermore, we have shown how these refutations
can be parameterized using the space of positive semidefinite
matrices. This has resulted in a nested sequence of sufficient
conditions, of increasing accuracy and expressible as semidefinite
programs, which prove stability of linear time-delay systems. Thus
for any desired level of accuracy, these results give a condition,
expressible as a semidefinite program, which will test stability to
that level of accuracy. These results appear in publication
in~\cite{peet_2006b} and~\cite{peet_2006c}. In addition, because of
the structure of the refutations, in Chapter~\ref{chp:nonlinearcase}
we have also been able to generalize these results to time-delay
systems with nonlinearities. Some of these results appear in
publication in~\cite{peet_2004a}.\\

\section{Prior Work}

The contents of this thesis utilize results from a number of
different areas of research. A complete survey of the results in any
one of these fields would be beyond the scope of this document. In
this section we briefly list some landmark contributions which have
directly influenced the direction of our research.

\paragraph{Stability Theory}
Early results include the work of the Russian mathematician A. M.
Lyapunov~\cite{lyapunov_1892}, who, in 1892, standardized the
definition of stability and generalized the potential energy work of
Lagrange~\cite{lagrange_1788} to systems of ordinary differential
equations of the form $\dot{x}(t)=f(x(t))$. In the century that has
followed, the use of Lyapunov functions to prove stability has
become commonplace and is known alternatively as the ``Direct method
of Lyapunov'' or ``Lyapunov's second method''.

In the 1960's, the input-output approach to stability analysis
emerged as an alternative to Lyapunov's method. This approach was
motivated by the development of complicated electronic systems for
which a detailed analysis of the stability of the internal states
was impractical. The work was pioneered by
Sandberg~\cite{sandberg_1964a,sandberg_1964a} and
Zames~\cite{zames_1966a,zames_1966b} and can be found in such works
as~\cite{willems_1971} and~\cite{desoer_1975}. The input-output
framework was of practical importance for a number of primary
reasons. First, it could be used with ``black-box'' frequency
sweeping techniques to quickly determine properties of complicated
linear systems even when those systems contained internal delay. In
addition, these frequency-domain properties could then be used to
predict the stability of the interconnection of systems through the
use of concepts such as passivity and small-gain. Development of the
input-output framework has had some relatively recent advances with
the introduction of the theory of ``Integral Quadratic
Constraints''(IQCs) by Rantzer and Megretskii~\cite{rantzer_1997}.
This work on IQCs generalizes the passivity framework by formulating
stability conditions which allow for the use of a broad class of
multipliers. This work will be discussed in more depth in
Chapter~\ref{chp:stabilitytheory}.

\paragraph{Semidefinite Programming}

The first stability problem to be posed as a linear matrix
inequality(LMI) can be attributed to Lyapunov himself, who stated
that the system
\[\dot{x}(t)=Ax(t)
\]
is stable if and only if there exists a $P\ge0$ such that
\[A^T P + P A < 0.
\]
Of course, in the absence of algorithms for this LMI, Lyapunov was
forced to express this conditional analytically as the unique
solution, $P$,  for arbitrary $Q > 0$, of the Lyapunov equation $A^T
P + PA =-Q$. The solution of LMIs arising in control theory through
the use of analytical solutions was continued in the 1940's and
1950's by the work of Lur'e, Postnikov, and others in the Soviet
Union.

The first general method for the solution of LMIs was developed in
the 1960's by Kalman, Yakubovich, and others. The famous
positive-real or KYP lemma showed that a certain class of LMI
problems could be solved graphically by considering a frequency
domain inequality. Today, of course, the KYP lemma is used more
often in the opposite sense, i.e. to convert a frequency-domain
inequality to a semidefinite program. The LMI associated with the
KYP lemma was later shown to admit an analytic solution through use
of the algebraic Riccati equation(ARE).

The most significant advance in the solution of LMI problems was the
recognition in the 1970's that many LMIs could be expressed as
convex optimization problems which could be solved using recently
developed efficient numerical algorithms for linear programming.
Thus the use of analytic solutions to LMI problems was replaced by
\emph{semidefinite programming} which is defined as optimization
over the convex cone of positive semidefinite matrices.

The algorithms for numerical optimization which were first used to
solve LMI problems were originally developed to solve linear
programming problems of the following form where the inequality is
defined by the positive orthant.
\begin{align*}
&\max c^T y:\\
&Ay \le b
\end{align*}

The solution of linear programming problems was first addressed by
the development of the Simplex algorithm by G. Dantzig in the
mid-40's. This algorithm performed well when applied to most
problems but was shown to fail in certain special cases. Concern
about this ``worst-case'' complexity prompted the search for what we
now refer to as polynomial-time algorithms, which have provable
bounds on worst-case performance. The first such polynomial-time
algorithm was introduced by Khachiyan in 1979 and is referred to as
the ellipsoid algorithm. This algorithm, however, proved inferior to
the Simplex algorithm in most practical cases. The first major
improvement over the Simplex algorithm was developed in 1984, when
N. Karmakar proposed a new type of algorithm which used what we now
refer to as an ``interior-point'' method. This algorithm was also of
polynomial-time complexity but performed far better than the
ellipsoid algorithm in practice. The final hurdle towards numerical
solution of LMI problems was overcome in 1988, when Nesterov and
Nemirovskii developed interior point algorithms which applied
directly to semidefinite programming problems. Today, semidefinite
programming programs can be solved simply and efficiently using
interior-point solvers such as SeDuMi~\cite{sturm_1999} combined
with user-friendly interfaces such as Yalmip. A summary of standard
LMI problems in control theory can be found in the seminal work by
Boyd et al.~\cite{boyd_1994}.

Recently, the success of semidefinite programming in addressing
problems in control has led to research into whether these
algorithms can be applied to convex problems in polynomial
optimization. Polynomial optimization problems arise in many areas
of nonlinear systems analysis and control and this topic will be
covered in significant depth in Chapter~\ref{chp:SOStheory}.

\paragraph{Internet Congestion Control}
The origins of research interest in Internet congestion control from
a mathematical perspective can be traced back to the the paper by
Kelly et al.~\cite{kelly_1998}, wherein the the congestion control
problem was first cast as a decentralized optimization problem. This
paper showed that certain congestion control protocols would
converge to the global optimum in the absence of delay through the
use of a Lyapunov argument. Subsequently, in Low and
Lapsley~\cite{low_1999}, it was shown that the dynamics of the
delayed Internet model with a certain class of control algorithms
could be interpreted as a decentralized implementation of the
asynchronous gradient projection algorithm to solve the dual to the
network optimization problem thus showing global convergence to
optimality for sufficiently small step size. However, no global
bound for the step size was given in this paper, making practical
interpretation difficult. This work was followed by the paper by
Paganini et al.~\cite{paganini_2001}, wherein it was shown that with
a certain set of pricing functions, a uniform bound of
$\alpha<\pi/2$ on a certain gain parameter $\alpha$ at the source
allows a proof of local stability for arbitrary topology and
heterogeneous time delays.  This work was similar in many respects
to work by Vinnicombe~\cite{vinnicombe_2000} for a different network
implementation. Our own work has shown this uniform bound of
$\alpha<\pi/2$ was shown to imply global stability for heterogeneous
time delays in a more limited topology in~\cite{peet_2004b}
and~\cite{peet_2006a}. All these results are discussed in more depth
in Chapters~\ref{chp:ICC} and~\ref{chp:IQCresult}.

\paragraph{Time-Delay Systems}
In 1963, the potential energy methods of Lagrange and Lyapunov were
generalized to systems of functional differential equations by N. N.
Krasovskii~\cite{krasovskii_1963}. Since this time, many
computationally tractable sufficient conditions have been given for
stability of both linear and nonlinear time-delay systems, all with
varying degrees of conservatism. An overview of some of these
results can be obtained from survey materials such as are found
in~\cite{gu_2003,hale_1993,kolmanovskii_1999,niculescu_2001}. These
results can be grouped into analysis either in the frequency-domain
or in the time-domain. Frequency-domain techniques can be applied to
linear systems only and typically attempt to determine whether all
roots of the characteristic equation of the system lie in the left
half-plane. This approach is complicated by the transcendental
nature of the characteristic equation, which imply the existence of
a possibly infinite number of roots. Time-domain techniques
generally use Lyapunov-based analysis. In the linear case, the
Lyapunov approach benefits from the existence of an
operator-theoretic version of Lyapunov's inequality, the existence
of a positive solution to which is necessary and sufficient for
stability. Computing solutions to this inequality, however, has
historically been problematic. A main contribution of this thesis is
to show how the operator-theoretic Lyapunov equation can be solved
directly using semidefinite programming. This work is described in
Chapter~\ref{chp:LinearCase} and the references~\cite{peet_2006b}
and~\cite{peet_2006c}. We note that one attempt to solve the
Lyapunov equation by considering piecewise-linear functions was made
in a series of papers by Gu et al. and is summarized
in~\cite{gu_2003}.

\section{Notation}
The following are used throughout this thesis.

\paragraph{Vector Spaces}
In this thesis, we use the following standard notation. $\R^{n
\times m}$ denotes the space of real $n \times m$ matrices. $\S^n$
denotes the space of symmetric $n \times n$ matrices. Let
$\R^+:=\{\, x\in\R \ \vert \ x\geq 0\,\}$. Let $\mathcal{C}(I)$
denote the set of continuous functions $u:I \rightarrow \R^n$ where
$I \subset \R$. We say $f \in \mathcal{C}(I)$ is \emph{bounded} if
there exists some $b \in \R^+$ such that $\norm{f(\theta)}_2 \le b$
for all $\theta \in I$. We use $\mathcal{C}_{\tau}$ to denote the
Banach space of continuous functions $u\in \mathcal{C}([-\tau,0])$
with norm $\norm{u} =\sup_{t\in [-\tau,0]} \, \norm{u(t)}_2$.

\begin{defn} For a given $\tau > 0$, $x \in \mathcal{C}([a,b))$, and
$t \in [a+\tau,b]$, where $b>a+\tau$, define $x_t \in
\mathcal{C}_\tau$ by $x_t(\theta)=x(t+\theta)$ for $\theta \in
[-\tau,0]$.
\end{defn}

$L_2(-\infty,\infty)$ is the Hilbert space of Lebesgue measurable
real vector-valued functions $x : \R \rightarrow \R^n$ with
inner-product $\ip{u}{v}_2=\int_{-\infty}^{\infty} u(t)^T v(t) dt$.
$L_2$ denotes $L_2[0,\infty)=\{\, x \in L_2(-\infty,\infty) \ \vert
\ x(t)=0 \text{ for all } t<0 \,\}$ and is a Hilbert subspace of
$L_2(-\infty,\infty)$. Similarly, $L_2[a,b]$ denotes the restriction
of $L_2(-\infty,\infty)$ to the interval $[a,b]$. Throughout, the
dimensions of $x(t)$ for $x \in L_2$ should be clear from context
and are not explicitly stated. We will occasionally also associate
with $\mathcal{C}_\tau$ an inner product space equipped with the
inner product associated with $L_2$. Thus for $x,y \in
\mathcal{C}_{\tau}$, $\ip{x}{y}$ denotes the $L_2$ inner product.
$\hat{L}_2$ denotes the Hilbert space of complex vector-valued
functions on the imaginary axis, $x : j \R \rightarrow \C^n$ with
inner-product $ \ip{\hat{u}}{\hat{v}}_2 =
\frac{1}{2\pi}\int_{-\infty}^{\infty} \hat{u}(\iw)^* \hat{v}(\iw)
d\omega$. $\hat{L}_{\infty}$ denotes the Banach space of
matrix-valued functions on the imaginary axis, $\hat{G} : j \R
\rightarrow \C^{m \times n}$ with norm
$\|\hat{G}\|_{\infty}=\text{ess} \sup_{\omega \in \R}
\bar{\sigma}\left( \hat{G}(\iw) \right)$ where
$\bar{\sigma}(\hat{G}(\iw))$ denotes the maximum singular value of
$\hat{G}(\iw)$. For a given $\tau>0$, let $M_2$ denote the product
space $\R^n \times L_2[-\tau,0]$ endowed with the inner product
\[\ip{x}{y}:=x_1^T y_1+ \ip{x_2}{y_2}_2,
\]
where we associate with $x\in M_2$ a pair $(x_1,x_2)$ where $x_1 \in
\R^n$ and $x_2 \in L_2[-\tau,0]$.

\paragraph{Maps and Operators}

An operator $A$ on an inner-product space is defined to be positive
on a subset $X$ if $\ip{x}{Ax}\ge0$ for all $x \in X$. A function
$x: \R \rightarrow \R$ is said to be \eemph{absolutely continuous}
if for any integer $N$ and any sequence $t_1,\dots, t_N$, we have
$\sum_{k=1}^{N-1} |x(t_k)-x(t_{k+1})|\rightarrow 0 $ whenever
$\sum_{k=1}^{N-1}|t_{k}-t_{k+1}|\rightarrow 0$. $P_T$ is the
truncation operator such that if $y=P_Tz$, then $y(t)=z(t)$ for all
$t\le T$ and $y(t)=0$ otherwise. $L_{2e}$ denotes the space of
functions such that for any $T>0$ and $y\in L_{2e}$, we have $P_T
y\in L_2$. We also make use of the space $W_2:=\{y\in L_2: \dot{y}
\in L_2\}$ with inner product
$\ip{x}{y}_{W_2}=\ip{x}{y}_{L_2}+\ip{\dot{x}}{\dot{y}}_{L_2}$ and
extended space $W_{2e}:=\{y\in L_{2e}:\dot{y} \in L_{2e}\}$. A
causal operator $H:L_{2e}\rightarrow L_{2e}$ is bounded if $H(0)=0$
and if it has finite gain, defined as
\[\norm{H}=\sup_{u\in L_2 \ne 0}\frac{\norm{Hu}}{\norm{u}}
\]
$\hat{u}$ denotes the either the Fourier or Laplace transform of
$u$, depending on $u$. We will also make use of the following
specialized set of transfer functions which define bounded linear
operators on $L_2$. $\mathcal{A}$ is defined to be those transfer
functions which are the Laplace transform of functions of the form
\[
g(t)=
\begin{cases}
h(t)+\sum_{i=1}^N g_i \delta(t-t_i) & \text{if } t\ge0\\
 0 & \text{otherwise}
\end{cases}
\]
where $h\in L_1$, $g_i\in \R$ and $t_i\ge0$.

\paragraph{Real Algebraic Geometry}

Denote by $\R[x]$ the ring of scalar polynomials in variables $x$.
Denote by $\R^{n \times m}[x]$ the set of $n$ by $m$ matrices with
scalar elements in $\R[x]$. Let $\S^n[x]$ denote the set of
symmetric $n$ by $n$ matrices with elements in $\R[x]$ and define
$\S_d^n[x]$ to be the elements of $\S^n[x]$ of degree $d$ or less.
Let $\mathcal{P}^{Y} \subset \R[x]$ denote the convex cone of scalar
polynomials which are non-negative on $Y$. Let
$\mathcal{P}^{+}\subset \R[x]$ denote the convex cone of globally
non-negative scalar polynomials. Let $\mathcal{S}_{n}^{Y}$ denote
the convex cone of elements $M \in \S^{n}[x]$ such that $M(x) \ge 0$
for all $x \in Y$. Let $\mathcal{S}_{n}^+$ denote the convex cone of
elements $M \in \S^{n}[x]$ such that $M(x) \ge 0$ for all $x$. We
let $Z_d[x]$ denote the $n+d \choose d$-dimensional vector of
monomials in variables $x\in \R^n$ of degree $d$ or less. Define
$\bar{Z}^n_d[x]:=I_n \otimes Z_d[x]$, where $I_n$ is the identity
matrix in $\S^{n}$ and $\otimes$ is the Kronecker product. Finally,
we note that for the sake of notational convenience, we will often
use the expression $M(x) \in X$ to indicate that the function $M \in
X$ for some set of functions $X$. This will hopefully not cause
substantial amounts of confusion.

\chapter{Functional Differential Equations}\label{chp:FDEtheory}

\section{Introduction}
The best way to introduce the concept of a functional differential
equation is through the use of an example.

\textbf{Example 1:} Perhaps the most easily understood example of a
system defined by a functional differential equation is that of the
dynamics of taking a shower. Specifically, define $\delta T(t)$ to
be the difference between water temperature and body temperature at
time $t$. We assume that the typical bather controls the water
temperature by turning the hot-water knob at rate $\omega(t)$,
proportional to the difference between water temperature and body
temperature so that $\dot{\omega}(t)=-\alpha \delta T(t)$. Ideally,
the position of the hot water knob is directly proportional to the
temperature of the water coming from the head, $\delta T(t)=\beta
\omega(t)$ for $\beta>0$. In this case, we have the following
ordinary differential equation.

\[\delta \dot{T}(t)= \alpha \omega(t) = -\alpha \beta \delta T(t)
\]

\noindent This is a linear, time-invariant system and consequently
we know that $\lim_{t \rightarrow \infty}\delta T(t)=0$ for any any
initial condition $\delta T(0)$ and any positive value of $\alpha$.
However, as most people know, there is occasionally a delay between
action on the hot-water knob and change of temperature at the head.
This delay occurs because hot water mixed at the knob must pass
through a length of pipe before arrival at the shower head. We
assume that the delay, $\tau$, is constant and can be calculated as
$\tau=L/v$; where $L$ is the linear distance from the tap to the
head and $v$ is the rate of the water flowing in the pipe. The
simplest description of the delayed system is given as follows.

\[\delta \dot{T}(t)=  -\alpha \beta \delta T(t-\tau)
\]

It can be shown~\cite{hale_1993} that a system of this form is
stable in the region $\alpha \beta \in (0,\frac{\pi}{2 \tau})$ and
unstable outside of this region. Thus we conclude that in the
presence of delay, any bather will get scalded if sufficiently
impatient. The reason that a large proportional feedback gain fails
in the shower example is that the temperature of the water at the
head does not provide an adequate representation of the state of the
system. In order to exactly predict exactly how the water
temperature will evolve over time, one needs to know the water
temperature, $T(\theta,t)$, at every point $\theta$ in the pipe from
the knob to the head for some time $t$. This information precisely
defines the state of the system at time $t$. Because the state at
time $t$, $T(\cdot,t)$ is a function, the dynamics of a shower are
defined by a function of a function or a \emph{functional}.
%\[
%\dot{T}(0,t)=f(T(\cdot,t))
%\]

In fact, the dynamics given above represent a simplification of the
more difficult problem of fluid flow described by partial
differential equations. We are able to represent the system using
such a simple model only because the controller and observer are
highly structured so that control and observation take place only at
discrete points in the flow. Such simplified models are quite common
and are collectively known as \emph{time-delay systems}.

\section{Definitions and the Concept of State}
To begin, we define the following.

\begin{defn}
Suppose we are given a $\tau\ge0$ and map $f: \mathcal{C}_\tau
\times \R^+ \rightarrow \R^n$. We say that a function $x\in
\mathcal{C}([-\tau,b))$ is a \emph{solution} on $[-\tau,b)$ to the
functional differential equation defined by $f$ with initial
condition $x_0 \in \mathcal{C}_\tau$ if $x$ is differentiable,
$x(\theta)=x_0(\theta)$ for $\theta \in [-\tau,0]$ and the following
holds for $t \ge 0$.
\begin{align}
\dot{x}(t)&=f(x_t,t)%\\
%x_t(\theta)&=x(t+\theta) \qquad \theta \in [-\tau,0]
\end{align}
\end{defn}

\begin{figure}[htb]
  \centerline{\includegraphics[width=0.5\textwidth]{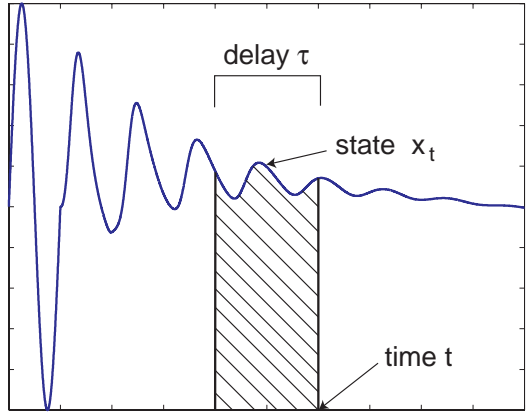}}
  \caption{An illustration of the state $x_t \in \mathcal{C}_\tau$}
  \label{fig:state_1}
\end{figure}

In this thesis, we make use of two different concepts of state
space. The first, and perhaps the most common is the Banach space
$\mathcal{C}_\tau$, equipped with the supremum norm. In this
scenario, the state of the system is simply the trajectory of the
system over the past $\tau$ seconds of time. It is in this space
that we will define solution maps and theorems of existence and
uniqueness. An example of such a state is illustrated in
Figure~\ref{fig:state_1}.

Another concept of state space which we find more satisfying for
applications such as stability is the space $M_2$, which was
thoroughly investigated by such work as that by Delfour and
Mitter~\cite{delfour_72a}.  For solution $x$, the state of the
system, $\phi(t)$ at time $t$ is defined as $\phi(t) = (x(t), x_t)$.
$\phi$ lies in the subspace $H_2 \subset M_2$ defined as follows.

\[H_2:=\{y \in M_2 \, : \, y_2 \in \mathcal{C}_\tau, \; y_1=y_2(0) \}
\]

This space better reflects the dependency of the system on both the
finite dimensional element $x(t)$ as well as the hereditary infinite
dimensional element $x_t$. In addition, use of this state space
allows us to embed the state space within the Hilbert space $M_2$, a
feature useful in describing the operators used to define Lyapunov
stability criteria.
%Although the following definitions and results can also be expressed
%in the space $M_2$, we present them here in $\mathcal{C}_\tau$, as
%this space will be more familiar to most readers.

\begin{defn} For a given functional, $f: \mathcal{C}_\tau \times
\R^+
\rightarrow \R^n$, suppose that for any $\phi \in \mathcal{C}_\tau$,
there exists a unique solution on $[-\tau,\infty)$ to the functional
differential equation defined by $f$. Define the \textbf{solution
map} $G_f: \mathcal{C}_\tau \rightarrow \mathcal{C}([-\tau,\infty))$
by $y=G_f(\phi)$ where $y$ is the solution to the functional
differential equation defined by $f$ with initial condition $\phi$.
\end{defn}

For the time-invariant case, where $f(x_t,t)=\tilde{f}(x_t)$, we
define the following, which maps the evolution of the state.

\begin{defn} For a given solution map, $G_f$, define the \textbf{flow
map} $\Gamma_f : \mathcal{C}_\tau \times \R^+ \rightarrow
\mathcal{C}_\tau$ by $\xi=\Gamma_f(\phi,\Delta t)$ if $\xi=y_{\Delta
t}$ where $y=G_f(\phi)$.
\end{defn}

The solution map, $G_f$, and flow map, $\Gamma_f$, provide
convenient notation for expressing the evolution of the system over
time. The following result gives conditions under which the solution
map is well-defined in the time-varying case.

\subsection{A Note on Existence of Solutions}
In this subsection, we give a theorem which provides conditions for
the existence and uniqueness of solutions defined on
$\mathcal{C}([-\tau,\infty))$.

\begin{thm}\label{thm:existuniq}
Suppose a functional $f:\mathcal{C}_\tau \times \R^+ \rightarrow
\R^n$ satisfies the following.
\begin{itemize}\item There exists a $K_1>0$ such that
\[\norm{f(x,t)-f(y,t)}_2 \le K_1 \norm{x-y}_{\mathcal{C}_\tau} \quad
\text{for all }x,y \in \mathcal{C}_\tau, t \ge 0
\]
\item There exists a $K_2>0$ such that
\[\norm{f(0,t)}\le K_2 \quad \text{for all }t \ge 0
\]
\item $f(x,t)$ is jointly continuous in $x$ and $t$. i.e. for every
$(x,t)$ and $\epsilon>0$, there exists a $\eta>0$ such that
\[\norm{x-y}_{\mathcal{C}_\tau}+\norm{t-s}\le \eta \; \Rightarrow \;
\norm{f(x,t)-f(y,s)}\le \epsilon.
\]
\end{itemize}
Then for any $\phi \in \mathcal{C}_\tau$, there exists a unique $x
\in \mathcal{C}[-\tau,\infty)$ such that $x$ is differentiable for
$t\ge0$, $x(t)=\phi(t)$ for $t\in [-\tau,0]$ and
\[\dot{x}(t)=f(x_t,t)\qquad \text{for all }t \ge 0
\]
\end{thm}

See Appendix~\ref{app:FDEtheory} for Proof.
%\begin{proof}
%We have already shown by Lemmas~\ref{lem:map}
%and~\ref{lem:contraction} that $\Lambda$, as defined by $f$, is a
%contraction on $V_\alpha$ for any $\alpha>K$. Therefore, there
%exists some $\alpha>K$, $x \in V_\alpha$ such that $x=\Lambda x$.
%Therefore, by Lemma~\ref{lem:existence}, we have that $x$ is a
%solution. By Lemma~\ref{lem:uniqueness}, we have that this solution
%is unique.
%\end{proof}

By use of Theorem,~\ref{thm:existuniq}, we can show that if a
functional, $f$, satisfies a global Lipschitz continuity condition,
then for any initial condition there exists a unique solution to the
functional differential equation defined by $f$ which is defined on
$[-\tau,\infty)$. In this case the solution map and, in the
time-invariant case, the flow map are well-defined.

\section{Concepts of Stability}

In this section we define stability of a functional differential
equation.

\subsection{Internal Stability}

Two concepts of stability will be used in this thesis. The first,
internal stability, defines stability as a property of solutions of
a functional differential equation for arbitrary initial conditions.
The second, input-output stability, is used to define stability of
an operator and describes the relationship between inputs and
outputs. In this subsection we give a definition of internal
stability.

\begin{defn}
  Assume that $f$ satisfies $f(0,t)=0$. The solution map $G_f$, defined by $f$, is
  \textbf{stable} on $X \subset \mathcal{C}_{\tau}$ if
\begin{enumerate}
\item[(i)] $G_f x$ is bounded for any $x \in X$
\item[(ii)] $G_f$ is continuous at $0$ with respect to the supremum norm on
  $\mathcal{C}([-\tau,\infty))$ and $\mathcal{C}_\tau$.
\end{enumerate}
\end{defn}

This is the usual notion of Lyapunov stability, which states that
for all $\eps >0$ there exists $\delta >0$ such that
$\norm{y}<\delta$ implies $\norm{G_f y} < \eps$.

\begin{defn}
The solution map $G_f$ defined by $f$ is \eemph{asymptotically
stable} on $X \subset \mathcal{C}_{\tau}$ if it is stable on $X$ and
$y=G_f x_0$ implies $\lim_{t \rightarrow \infty} y(t) = 0$ for any
$x_0 \in X$.
\end{defn}

\begin{defn}The solution map, $G_f$, is \eemph{globally stable} if
it is stable on $\mathcal{C}_\tau$
\end{defn}

\begin{defn}The solution map, $G_f$, is \eemph{globally asymptotically stable} if
it is asymptotically stable on $\mathcal{C}_\tau$.
\end{defn}

\textbf{Note:} For a system defined by a linear functional, global
stability is equivalent to stability on any open neighborhood of the
origin. See, for example~\cite{gu_2003,kolmanovskii_1999}.

\subsection{Input-Output stability}

We can associate with a functional, $f$, an input-output system.

\begin{defn}
For functional $f$ and function $g$, let $y=\Psi_{f,g} u$ if
$y(t)=g(x(t))$ for $t \in \R^+$ where $x=G_{\hat{f}}(0)$ and
\[\hat{f}(x_t,t)=f(x_t,t)+u(t).
\]
\end{defn}

We can now define input-output stability directly in terms of
properties of the operator $\Psi$.

\begin{defn} For normed spaces $X,Y$, the operator $\Psi$ is $Y$ stable
on $X$ if it defines a single-valued map from $X$ to $Y$ and there
exists some $\beta$ such that $\norm{\Psi u}_Y \le \beta \norm{u}_X$
for all $u \in X$.
\end{defn}

\textbf{Note:} We refer to $X$ stability on $X$ as simply $X$
stability.

\section{Time-Delay Systems} So far, we have only considered the
general case of a dynamic system defined by a functional. In this
section, we identify a specific class of functional differential
equation which will be of particular importance throughout this
thesis.

\begin{defn}
We say that a functional, $f$, defines a \eemph{Time-Delay System}
if $f$ can be represented using a function $p:\R^{n(K+1)+1}
\rightarrow \R^n$ in the following way where $\tau_i> \tau_{i-1}$
for $i=1, \ldots, K$ and $\tau_0=0$.
\begin{align*}
f(x_t)=\int_{-\tau_K}^0 p(x_t(-\tau_0),x_t(-\tau_1),\cdots,x_t(-\tau_K), x_t(\theta), \theta) d \theta\\
\end{align*}
\end{defn}

Stability of Time-Delay systems are generally classified using the
following definitions.
\begin{defn} A Time-Delay System is
\eemph{Delay-Independent Stable} if it is stable for arbitrary
$\tau_i \in \R$ for $i=1 \ldots K$.
\end{defn}

\begin{defn} A Time-Delay System is
\eemph{Delay-Dependent Stable} on $X$ if it is stable for $\tau_i
\in X_i$ for $i=1 \ldots K$ where the $X_i$ are compact subsets of
$\R$.
\end{defn}

\subsection{The Case of Linear Time-Delay Systems}

In this subsection, we consider the special case of linear
time-delay systems. Specifically, we consider functionals which can
be expressed in the following form, where $\tau_i > \tau_{i-1}$ for
$i=1,\ldots,K$ and $\tau_0=0$.

\begin{equation}
f(x_t)=\sum_{i=0}^K A_i x(t-\tau_i)+\int_{-
\tau_K}^{0}A(\theta)x(t+\theta)d \theta \label{eqn:general_lfde1}
\end{equation}

Here $A_i\in \R^{n \times n}$ and $A:\R \rightarrow \R^{n \times n}$
is bounded on $[-\tau_K,0]$. The following lemma, combined with
Theorem~\ref{thm:existuniq} shows that elements of this class of
system admit a unique solution for every initial condition $x_0 \in
\mathcal{C}_{\tau_K}$.

\begin{lem}\label{lem:linear_existence}
Let
\[f(x,t):=\sum_{i=1}^K A_i x(t-\tau_i)+\int_{-\tau_K}^0 A(\theta)x(t+\theta)d\theta\]
Where $A(\theta)$ is bounded on $[-\tau_K,0]$. Then there exists
some $k>0$ such that
\[\norm{f(x,t)-f(y,t)}_2 \le k \norm{x-y}_{\mathcal{C}_\tau}
\]
\end{lem}

See appendix for Proof.

By using Lemma~\ref{lem:linear_existence}, we can associate with
linear systems of this form a well-defined solution map $G_f$ and
flow map $\Gamma_f$.

\section{Conclusion}
In this chapter, we have introduced various concepts associated with
functional differential equations. These concepts include the state
of the system, the solution map and stability in both the internal
and input-output framework. We will make use of these definitions
throughout the remainder of this thesis.

\chapter{Stability of Functional Differential Equations}\label{chp:stabilitytheory}

In this chapter, we introduce two methods of proving stability of
functional differential equations. The first, which can be used to
prove internal stability, uses a generalization of Lyapunov theory
to functional differential equations. The second method, which is
used to prove input-output stability, is a generalization of the
notion of passivity.

\section{The Direct Method of Lyapunov}

Consider a solution map $G_f$ defined by a functional $f$. We
consider the state of the system at time $t$ to be defined by an
element $x_t\in \mathcal{C}_\tau$. Standard Lyapunov theory,
however, is defined using functions of the form $V(x(t))$. Such
functions capture only part of the energy of the state, $x_t$, i.e.
that part stored in $x_t(0)$. Therefore, any stability condition
derived from such functions will be inherently conservative. An
attempt to address this conservatism was made by
Krasovskii~\cite{krasovskii_1963} through the introduction of
Lyapunov functionals which depend on elements of $\mathcal{C}_\tau$.

%\begin{defn}Let $V:\mathcal{C}_{\tau} \rightarrow \R$ be a continuous
%function such that $V(0)=0$. We say that $V$ is a \eemph{candidate
%Lyapunov functional} if there exist non-decreasing continuous
%functions $u,v : \R^+ \rightarrow \R^+ $ such that $u(s),v(s)>0$ for
%$s \neq 0$, $u(0)=v(0)=0$, and $u(s)\rightarrow \infty$ as $s
%\rightarrow \infty$ and such that the following holds for all $\phi
%\in \mathcal{C}_\tau$.
%\begin{equation}
%u(|\phi(0)|) \le V(\phi) \le v(\norm{\phi})
%\end{equation}
%\end{defn}

\begin{defn} Let $V:\mathcal{C}_{\tau} \rightarrow \R$ be a continuous
function such that $V(0)=0$. For a given flow map, $\Gamma_f$,
defined by solution map, $G_f$, define the upper Lie derivative
of~$V$ as follows.
\[
\dot{V}(\phi):=\limsup_{h \rightarrow
0^+}\frac{1}{h}[V(\Gamma_f(\phi,h))-V(\phi)]
\]
\end{defn}
The following theorem follows from Gu~\cite{gu_2003}.

\begin{thm}\label{thm:Hale}
Let $\Omega \subset \mathcal{C}_\tau$ contain an open neighborhood
of the origin. Let $f : \Omega \rightarrow \R^n$ be continuous and
take bounded sets of $\Omega$ into bounded sets of $\R^n$. Let
$V:\mathcal{C}_{\tau} \rightarrow \R$ be a continuous map such that
$V(0)=0$. Let $u,v,w : \R^+ \rightarrow \R^+ $ be continuous
non-decreasing functions such that $u(s),v(s)>0$ for $s \neq 0$ and
$u(0)=v(0)=0$. Suppose that the following holds for any $\phi \in
\Omega$.
\begin{align*}
u(\norm{\phi(0)}) &\le V(\phi) \le v(\norm{\phi})\\
\dot{V}(\phi) &\le -w(|\phi(0)|)
\end{align*}
Then the solution map, $G_f$, defined by $f$ is stable on some open
neighborhood of the origin. If $w(s)>0$ for $s>0$, then the system
is asymptotically stable on some open neighborhood of the origin. If
$\Omega=\mathcal{C}_\tau$ and $\lim_{s \rightarrow \infty}= \infty$,
then $G_f$ is globally asymptotically stable.
\end{thm}

\subsection{Complete Quadratic Lyapunov Functionals}

There have been a number of results concerning necessary and
sufficient conditions for stability of linear time-delay systems in
terms of the existence of quadratic functionals. These results are
significant in that they allow us to restrict our search for a
Lyapunov-Krasovskii functional to a specific class without
introducing any conservatism. Consider a linear functional of the
following form.

\begin{equation}
f(x_t)=\sum_{i=0}^K A_i x(t-\tau_i)+\int_{-
\tau_K}^{0}A(\theta)x(t+\theta)d \theta \label{eqn:general_lfde2}
\end{equation}
We now make the additional assumption that $A$ is continuous on
$[-\tau_K,0]$. The following comes from Gu et al.~\cite{gu_2003}.

\begin{defn} We say that a functional $V:\mathcal{C}_\tau \rightarrow
\R$ is of the complete quadratic type if there exists a matrix $P
\in \S^n$ and matrix-valued functions $Q:\R \rightarrow \R^{n\times
n}$, $S:\R \rightarrow \S^n$ and $R: \R^2 \rightarrow \R^{n \times
n}$ where $R(\theta,\eta)=R(\eta,\theta)^T$ such that the following
holds.

\begin{align*}&V(\phi)=\phi(0)^T P \phi(0) +2 \phi(0)^T
\int_{-\tau_K}^{0} Q(\theta) \phi(\theta) d \theta \notag\\
&+ \int_{-\tau_K}^{0}\phi(\theta)^T S(\theta) \phi(\theta) d \theta
+\int_{-\tau_K}^{0} \int_{-\tau_K}^{0}\phi(\theta)^T R(\theta,
\eta)\phi(\eta) d \theta d \eta
\end{align*}
\end{defn}

\begin{thm} Suppose the system defined by
 a linear functional of the form~\eqref{eqn:general_lfde2} is asymptotically stable. Then
there exists a complete quadratic functional $V$ and $\eta>0$ such
that the following holds for all $\phi \in \mathcal{C}_\tau$.
\[V(\phi)\ge \eta \norm{\phi(0)}^2 \qquad  \text{and} \qquad \dot{V}(\phi)\le -\eta \norm{\phi(0)}^2
\]
Furthermore, the matrix-valued functions which define $V$ can be
taken to be continuous everywhere except possibly at points
$\theta,\eta=-\tau_i$ for $i=1 ,\ldots, K-1$.
\end{thm}

%\section{Linear Functional Differential Equations}
%
%
%
%%%%%%%%%%%%%%%%%%%%%%%%%%%%%%%%%%%%%%%%%%%%%%%%%%%%%%%%%%%%%%%%
%\subsection{Stability Theorem}
%%% Lyapunov-Krasovskii stability theorem
%Lyapunov Theory can be extended to time-delay systems through the
%use of Lyapunov-Krasovskii functionals. One such extension is
%defined by the following general theorem~\cite{hale_1993} which can
%be applied to the nonlinear as well as the linear case.
%\begin{thm}
%Consider a solution map $G$ defined by
%Equation~\eqref{eqn:general_lfde1}. Suppose $u,v,w:\R^+ \rightarrow
%\R^+$ are continuous nondecreasing functions, and that
%$u(s)>0,v(s)>0$ for $s>0$ and $u(0)=v(0)=0$. If there exists a
%continuous functional $V:\mathcal{C} \rightarrow \R$ such that both
%of the following hold for all $\phi \in \mathcal{C}_\tau$
%\begin{align}
%&u(|\phi(0)|)\le V(\phi) \le v(\norm{\phi}) \label{con:condition1}\\
%&\dot{V}(\phi)=\limsup_{\Delta t \rightarrow 0^+} \frac{1}{\Delta
%t}\left(V(\Gamma(\phi,\Delta t))-V(\phi)\right) \le -w(|\phi(0)|)
%\notag
%\end{align}
%Where $\Gamma$ is the flow map defined by $G$. Then the solution map
%$G$ is stable. If $w(s)>0$ for $s>0$, then the solution map is
%asymptotically stable.
%\end{thm}
%

\section{The Method of Integral Quadratic Contraints}
The theory of integral quadratic constraints(IQCs) is a method of
proving stability of the interconnection of stable operators.

\subsection{Theory of Integral-Quadratic
Constraints}\label{sec:IQC_theory}

Consider the following interconnection of operators.

\begin{figure}[htb]
  \centerline{\includegraphics[width=0.4\textwidth]{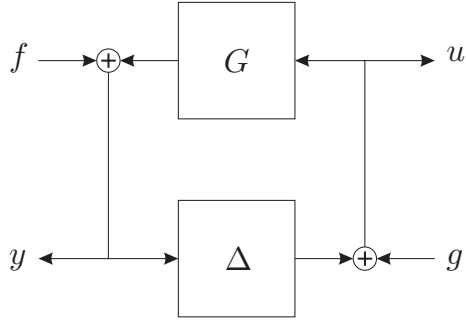}}
  \caption{Interconnection of systems}
  \label{fig:general_interconnection3}
\end{figure}
\begin{defn}
Let $G$ be a linear operator with transfer function
$\hat{G}\in\mathcal{A}$ and let the operator
$\Delta:L_{2}\rightarrow L_{2}$ be causal and bounded. Define inputs
$f \in W_2$ and $g \in L_2$. The \eemph{interconnection} of $G$ and
$\Delta$, denoted $\Phi_{F,G}$, is defined by $(y,u)=\Phi(f,g)$
where $y$ and $u$ are defined as follows.
\begin{align*}
y&=G u +f \notag \\
u&=\Delta y +g \notag
\end{align*}
\end{defn}

\begin{defn}[J\"{o}nsson~\cite{jonsson_1996}, p71]
The interconnection of $G$ and $\Delta$, $\Phi_{G,\Delta}$, is
\eemph{well-posed} if for every pair $(f,g)$ with $f \in W_2$ and $g
\in L_2$, there exists a solution $u\in L_{2e}$, $y \in W_{2e}$ and
the map $(f,g)\rightarrow (y,u)$ is causal.
\end{defn}

If the interconnection of $\Delta$ and $G$ is well posed, then the
interconnection defines an operator $\Phi_{G,\Delta} : W_2 \times
L_2 \rightarrow W_{2e} \times L_{2e}$. In this thesis, we use a
result by Rantzer and Megretski~\cite{rantzer_1997} which can be
interpreted as generalization of the classical notion of passivity.
Recall the the following classical passivity theorem from e.g.
Desoer and Vidyasagar(p. 182,~\cite{desoer_1975})
\begin{thm} The interconnection of $\Delta$ and $G$ is $L_2$ stable
on $L_2$ if there exists some $\epsilon>0$ such that for any $x \in
L_2$,
\begin{align*}
\left\langle \Delta x,x \right\rangle & \ge 0\\
\ip{x}{G x}& \le -\epsilon \norm{x}
\end{align*}
\end{thm}

Now given bounded linear transformations $\Pi_1, \Pi_2$, define the
following functional

\[\ipp{x}{y}:=\left\langle \Pi_1 \bmat{x\\y},\Pi_2 \bmat{x\\y}\right\rangle\]

Ignoring technical details for the moment, the result by Rantzer and
Megretski states that the interconnection of $\Delta$ and $G$ is
stable if there exists an $\epsilon>0$ such that for any $x \in
L_2$,
\begin{align*}
\ipp{x}{\Delta x} &\ge 0\\
\ipp{G x}{x} &\le -\epsilon \norm{x}
\end{align*}

This idea of generalized passivity is motivated from a geometric
standpoint by the topological separation argument, introduced by
Safonov~\cite{safonov_1982}. Consider the following definition of an
operator graph.

\begin{defn}
For an operator, $\rho : X \rightarrow X$, the \eemph{graph} of
$\rho$ is the set $\Phi(\rho):=\{(x,y): y=\rho(x), x \in X \} $. The
\eemph{inverse graph} of $\rho$ is the set $\Phi_i(\rho)=\{(x,y):
x=\rho(y), y \in X \}$.
\end{defn}

\begin{figure}[htb]
\centerline{\includegraphics[width=0.45\textwidth]{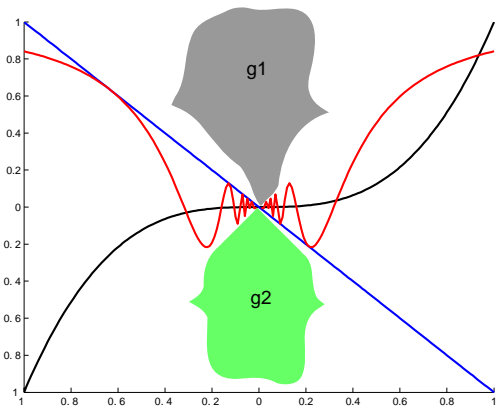}}
\caption{Separating functions prove separation of graphs $g_1$ and
$g_2$} \label{fig:g1andg2}
\end{figure}

Two graphs are separate if they intersect only at the origin.
Consider graphs $g_1$ and $g_2$.  If, for some functional $\sigma :
X \rightarrow \R$ and $\epsilon > 0$, we have $\sigma(x)\ge 0$ for
all $ x \in g_1 $ and $\sigma(y) \le -\epsilon \norm{y}$ for all $y
\in g_2$, then the graphs $g_1$ and $g_2$ can only intersect at the
origin and are said to be separated by the functional $\sigma$.  Now
consider the interconnection of operators $G$ and $\Delta$ when we
let the input $g=0$. We have the equation $f=(I-G \Delta)y$. The
question of $L_2$ stability becomes the question of existence of a
well-defined and bounded inverse of $(I-G \Delta)$ on $L_2$. Now
suppose there exists a nonzero $y\in L_2$ such that $f=(I-G
\Delta)y=0$, then $(I-G \Delta)$ has a nontrivial kernel and thus
can not have a bounded inverse. To ensure that $(I-G \Delta)$ has a
trivial kernel, we will consider the graph of $G$ and the inverse
graph of $\Delta$. The following well-known result shows that
separation of the graph and inverse graph of interconnected
operators is necessary for stability.

\begin{thm} Let $f=(I-G \Delta)y$. The following are equivalent
\begin{itemize}
\item $(I-G \Delta)y =0$ \; $\Rightarrow$ \; $y=0$ \item $\Phi(G)
\cap \Phi_i(\Delta) = 0$
\end{itemize}
\end{thm}
\begin{proof} $(\Leftarrow)$ If there exists a $y \ne 0$ such that $(I-G
\Delta)y=0$ then let $x=\Delta y$. Thus
\[\bmat{x \\ G x}=\bmat{\Delta y \\ y} \; \in \; \Phi(G) \cap
\Phi_i(\Delta)
\]
$(\Rightarrow)$ If there exists $y,x \ne 0$ such that
\[\bmat{x \\ G x}=\bmat{\Delta y \\ y}\]
then $(I-G \Delta)y = y - G \Delta y= G x -G \Delta y=G \Delta y- G
\Delta y = 0$
\end{proof}

Although the separation of graphs is clearly necessary for most
definitions of stability, it is unclear under what additional
conditions separation may also be sufficient for $L_2$ stability.
The result by Rantzer and Megretiski gives a particular class of
functionals for which graph separation is sufficient for stability
on $L_2$. Many classical theorems concerning the stability of the
interconnection of operators can be viewed as proving the separation
of graphs and inverse graphs by using these type of separating
functionals. For example, the small gain theorem can be expressed
using the separating functional $\sigma((x,y))=\norm{x}-k \norm{y}$
for some $k>0$. Similarly, classical passivity can be expressed
using the functional $\sigma((x,y))=\ip{x}{y}$. The complete class
of functionals shown by~\cite{rantzer_1997} to be sufficient for
$L_2$ stability is given by Definition~\ref{dfn:quad_cont}.

\begin{defn}[Rantzer~\cite{rantzer_1997}]
  The mapping $\sigma : L_2 \rightarrow \R$ is
\eemph{quadratically
    continuous} if for every $\delta > 0$, there exists a
  $\eta_{\delta}$ such that the following holds for all $x_1, x_2 \in L_2$.
  \[
  |\sigma(x_1)-\sigma(x_2)| \le \eta_{\delta}\|x_1-x_2\|^2 +\delta
\|x_2\|^2
  \] \label{dfn:quad_cont}
\end{defn}

This class includes the small gain and passivity functions.
Furthermore, for any bounded linear transformations $\Pi_1,\Pi_2$,
the function $\sigma(w)=\ip{\Pi_1 w}{\Pi_2 w} $ is quadratically
continuous. In this thesis we use the following generalization of
the work by Rantzer and Megretski as presented in the thesis work by
J\"{o}nsson~\cite{jonsson_1996}.

\begin{defn}Let $\Pi_B:j \R \rightarrow \C^{n\times n}$ be a bounded and measurable
function that takes Hermitian values and $\lambda \in \R$. We say
that $\Delta$ \eemph{satisfies the IQC} defined by $\Pi_B,\lambda$,
if there exists a positive constant $\gamma$ such that for all $y
\in W_2$ and $v=\Delta y \in L_2$,
\[\frac{1}{2 \pi}\int_{-\infty}^{\infty}\bmat{\hat{y}(\iw)\\\hat{v}(\iw)}^*
\Pi_B(\iw)\bmat{\hat{y}(\iw)\\\hat{v}(\iw)}d\omega + 2
\ip{v}{\lambda \dot{y}} \ge -\gamma|y(0)|^2\]
\end{defn}
\vskip 2mm
\begin{thm}Assume that
\begin{enumerate}
\item G is a linear causal bounded operator with
$s\hat{G}(s),\hat{G}(s)\in \mathcal{A}$ \label{IQCcondition3}
\item For all $\kappa \in [0,1]$, the interconnection of $\kappa
\Delta$ and $G$ is well-posed \item For all $\kappa \in [0,1]$,
$\kappa \Delta$ satisfies the IQC defined by $\Pi_B,\lambda$ \item
There exists $\eta>0$ such that for all $\omega \in \R$
\[\bmat{\hat{G}(\iw)\\I}^*\left(\Pi_B(\iw)+\bmat{0 & \lambda \iw^*\\\lambda \iw &
0}\right)\bmat{\hat{G}(\iw)\\I}\le - \eta I
\]\label{cdn:condition4}
\end{enumerate}
Then the interconnection of $G$ and $\Delta$ is $W_2 \times
L_2$ stable on $W_2 \times L_2$.\\
\label{thm:stability_thm}\end{thm}

\section{Conclusion}
In this chapter, we have introduced methods of constructing proofs
for internal and input-output stability. In the rest of this thesis,
we shall show how these methods can be applied to prove stability of
different kinds of time-delay systems.

\chapter{Internet Congestion Control}\label{chp:ICC}

\section{Introduction}
The analysis of Internet congestion control protocols has received
much attention recently. Explicit mathematical modeling of the
Internet has allowed analysis of existing protocols from a number of
different theoretical perspectives and has generated some
suggestions for improvement. This work has been motivated by concern
about the ability of current protocols to ensure stability and
performance of the Internet as the number of users and amount of
bandwidth continues to increase. Although the protocols that have
been used in the past have performed remarkably well as the Internet
has increased in size, analysis~\cite{low_2003} indicates that as
capacities and delays increase, instability will become a problem.

The purpose of this chapter is to outline the development of a
mathematical theory of Internet congestion control. Most of the
recent research activity in Internet congestion control protocols
can be traced back to the the paper by Kelly et
al.~\cite{kelly_1998}, wherein the the congestion control problem
was first cast as a decentralized optimization problem.
Subsequently, in Low and Lapsley~\cite{low_1999}, it was shown that
the dynamics of the Internet with a certain class of control
algorithms could be interpreted as a decentralized implementation of
the gradient projection algorithm to solve the dual to the network
optimization problem, thus showing global convergence to optimality
for sufficiently small step size. Finally, in Paganini et
al.~\cite{paganini_2001} it was shown that for a certain set of
utility functions, a global bound of $\alpha<\pi/2$ on a certain
parameter $\alpha$ at the source allows a proof of local stability
for arbitrary topology and heterogeneous time delays. These
fundamental results will form the background for our own work in
Chapter~\ref{chp:IQCresult}, wherein we show that the parameter
bound of Paganini et al. guarantees global stability of the
nonlinear implementation.

\section{The Internet Optimization Problem} \label{sec:problem}

\begin{figure}[htb]
  \centerline{\includegraphics[width=0.5\textwidth]{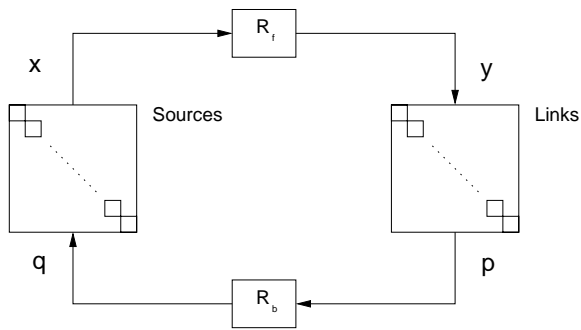}}
  \caption{The Internet is modeled as a collection of sources and links}
  \label{fig:internet_model}
\end{figure}

We view the Internet as an abstract collection of sources and links.
The term \eemph{source} refers to a connection between a user and a
particular destination. The source transmits data in packets. The
rate at which a source $i$ transmits packets is dictated by the
source's round-trip time, $\tau_i$, as well as window size, $w_i$.
The \emph{round-trip time} is defined as the time between
transmission of a packet and receipt of an acknowledgement for that
packet from the destination. The \emph{window size} is defined as
the number of packets which are allowed to be simultaneously
unacknowledged. In this thesis, we assume that packet losses do not
affect the source transmission rates, since any lost packet will
presumably be detected by the user after one round-trip time and
resent. We assume that acknowledgements contribute to delay but do
not contribute to congestion at the links. We assume a fixed bit
size for all packets and that $\tau_i$ is known at least for the
purposes of determining data rate.  The packet transmission rate,
$x_i$, at source $i$ can be controlled by the window size according
to
\begin{equation}
 \label{eqn:source_rate}
 w_i = x_i \tau_i.
\end{equation}

The term \eemph{link} refers to a single congested resource such as
a router. Packets arriving at a link enter an entrance queue. A link
can process packets in the queue at some rate capacity $c_j$. If too
many data packets arrive in a given period of time, the size of the
queue may grow and some packets may experience a queueing delay
while in the queue. In this thesis, we assume that the dynamics from
this variable queueing delay are negligible and we only model the
delay due to the fixed propagation time. Links must also be able to
feed back information. This can be done either through the ECN bit
in the packet header, through packet dropping schemes or through
measurement of variations in queueing delay at the source.  The
value of the congestion indicator at link $j$ is denoted $p_j$. We
also assume that the congestion indicator received at each source is
the summation of the indicators of all links in the source's route.
This value is denoted $q_i$.

Sources and links are related by routing tables which specify the
route or set of links, $J_i$ through which the packets from source
$i$ to its destination must pass. The rate of packets received at a
link $j$ is then the sum of the rates of all sources using that link
and is denoted by $y_j$. The set of users for link $j$ is denoted
$I_j$. Ignoring delay for the moment, we have the following
equations.
\begin{align*}
y &= R x, \quad q = R^T p,
\end{align*}
where
\[
R_{ji}=\begin{cases}
1 & \text{if source $i$ uses link $j$}\\
0 & \text{otherwise}
\end{cases}
\]

%%%%%%%%%%%%%%%%%%%%%%%%%%%%%%%%%%%%%%%%%%%%%%%%%%%%%%%%%%%%%%%%%%%%%%%%%%%%%%

\section{Optimization Model} The following model for
optimizing flow rates in a network was proposed by Kelly et
al.~\cite{kelly_1998}.
\begin{align*}
\text{maximize} \qquad & \sum_{i}^N U_i(x_i) \\
\text{subject to} \qquad &  x \ge 0, \quad R x \leq c
\end{align*}

Assume that the $U_i$ are continuously differentiable strictly
concave non-decreasing functions. If all sources utilize at least
one link, then the problem has a unique optimum. Note that, as $N$
increases, the problem becomes progressively more difficult to solve
using a centralized algorithm. We now consider the dual problem with
dual variable $p\in\R^M$, where $M$ is the number of links, which is
given by
\begin{align*}
\text{minimize} \qquad & h(p)\\
\text{subject to} \qquad & p \geq 0
\end{align*}
where the dual function $h$ is given by
\begin{align*}
h(p)&=\max_{x\ge 0}  \sum_{i} \bl( U_i(x_i)\br) - p^T(Rx-c)\\
&= \sum_{i} \bbl( U_i \bl( x_{\text{opt},i}(p)\br)\bbr) - p^T(R x_{\text{opt}}(p)  - c )\\
x_{\text{opt},i}(p) &= \max\{0,U_i'^{-1}(\sum_{j} R_{j,i} p_j)\}\\
&= \max\{0,U_i'^{-1}(q_i(p))\} \\
q(p)&=R^T p
\end{align*}
The map $U_i'^{-1}: \R^+ \rightarrow \{ \R \cup \infty \}$ is well
defined since $U_i'\in \mathcal{C}$ and $U_i$ is strictly concave.
We would like to construct a dynamical system which converges to the
solution of the dual problem. One such system is given by the
gradient projection algorithm. In discrete-time, this is
\begin{equation*}
  p_j(t+1)=\max\{ 0,  p_j(t) - \gamma_j D_j h(p(t)) \},
\end{equation*}
where $D_j$ denotes the partial derivative with respect to the
$j$'th argument and $\gamma_j$ is the step size. Since the $U_i$ are
strictly concave, $h(p)$ is continuously differentiable with the
following derivatives~\cite{bertsekas_1989}.
\begin{align*}
D_j h(p)&=c_j-\sum_{i \in I_j}
x_{\text{opt},i}\\
&=c_j-y_{\text{opt},j}(p)\\
y_{\text{opt}}(p)&=R x_{\text{opt}}(p).
\end{align*}
If $\gamma$ is sufficiently small, the discrete-time gradient
projection algorithm will converge to the solution of the dual
problem~\cite{low_1999}. Because of convexity of the problem, strong
duality implies that $x_{\text{opt}}(p)$ will converge to the unique
optimum of the primal problem. A continuous-time implementation of
this algorithm in the network framework is as follows.
\begin{align*}
  \label{eqn:continous_gp}
\dot{p}_j(t) &=
\begin{cases}
  \gamma_j(y_j(t)-c_j) & p_j(t) > 0  \\
  \max \{ 0, \gamma_j(y_j(t)-c_j) \}  & p_j(t) \leq 0
\end{cases} \\
x_{i}(t)&=\max\{0,U_i'^{-1}(q_i(t))\}\\
\quad y(t) &= R x(t),\quad q(t) = R^T p(t)
\end{align*}

$\gamma_j$ now denotes a gain parameter, corresponding to step-size
in discrete time. This algorithm has the remarkable property that it
is decentralized, corresponding to the separable structure of the
constraints. $p_j$ is computed at each of $M$ links. Link $j$
requires only knowledge of $y_j$ to compute this value. $x_{i}$ is
computed at each of $N$ sources. Source $i$ requires only knowledge
of $q_i$ to compute this value.

%%%%%%%%%%%%%%%%%%%%%%%%%%%%%%%%%%%%%%%%%%%%%%%%%%%%%%%%%%%%%%%%%%%%%%%%%%%%%%

\section{Stability Properties}

To ensure that the continuous-time gradient projection algorithm
will converge when implemented with the current Internet framework,
we must also consider the delay in transmitting packets from the
source to the link and then receiving acknowledgements at the
source. The delay from source $i$ to link $j$ is denoted
$\tau^f_{ij}$ and the delay from link $j$ to source $i$ is denoted
$\tau^b_{ij}$. For any source $i$, the total round trip time is
fixed, i.e. $\tau_i= \tau^f_{ij}+\tau^b_{ij}$ for all $j \in J_i$.
We express these delays in the frequency domain by replacing the
entries of the routing matrix $R$ with forward and backward delay
transfer functions $\hat{R}^f$ and $\hat{R}^b$, giving
\begin{align*}
\hat{y}(s)&=\hat{R}^f(s) \hat{x}(s), \quad \hat{q}(s)=\hat{R}^b(s)^T \hat{p}(s)\\
\hat{R}^f_{ji}(s)&=\begin{cases}e^{-\tau^f_{ij}s} & \text{if source $i$ uses link $j$}\\
0 & \text{otherwise}\end{cases}\\
\hat{R}^b_{ji}(s)&=\begin{cases}e^{-\tau^b_{ij}s} & \text{if source $i$ uses link $j$}\\
0 & \text{otherwise}\end{cases}
\end{align*}
The work by Paganini et al.~\cite{paganini_2001} introduced a class
of utility functions under which this system was shown to have a
stable linearization about its positive equilibrium point for a
fixed gain parameter $\gamma_j = 1/c_j$. This class was given by the
set of $U_i$ such that
\[
\frac{d}{d q_i} U_i'^{-1}(q_i)=-\frac{\alpha_i}{M_i \tau_i}
U_i'^{-1}(q_i),
\]
where $M_i$ is a bound on the number of links in the path of source
$i$ and $\alpha_i < \pi/2$. In particular, the choice of
\[
U_i(x)=\frac{M_i \tau_i}{\alpha_i} x
\bbbl(1-\ln\frac{x}{x_{\text{max},i}} \bbbr),
\]
with restricted domain $x\le x_{\text{max},i}$ was suggested
in~\cite{paganini_2001} as a strictly concave utility function such
that the function $U_i'^{-1}(q)=x_{\text{max},i}
e^{-\frac{\alpha_i}{M_i \tau_i}q} \ge 0$ has the necessary
derivative.

\section{Recent Work}
\begin{figure}[htb]
  \centerline{\includegraphics[width=0.4\textwidth]{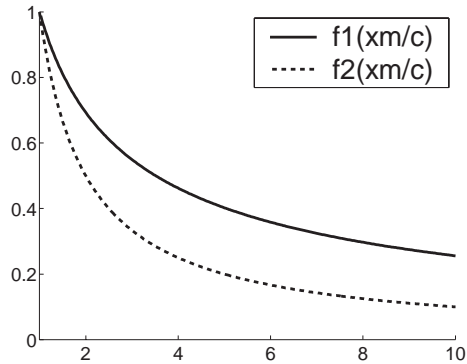}}
  \caption{Plot of $f_1(x)$ and $f_2(x)$ vs. $x$}
  \label{fig:f1f2plot}
\end{figure}

Some efforts have been made to extend this local stability result to
the global case. For a single source and a single link, the paper by
Wang and Paganini~\cite{wang_2002} has shown this implementation to
be globally stable for $\alpha \le f_1(x_{\text{max}}/c)$, where
\[f_1(x)=\frac{\ln x}{x-1}.
\]
In addition, the paper by
Papachristodoulou~\cite{papachristodoulou_2004a} has shown this
implementation to be stable when $\alpha \le f_2(x_{\text{max}}/c)$,
where
\[f_2(x)=\frac{1}{x}.
\]
$f_1$ and $f_2$ are illustrated in Figure~\ref{fig:f1f2plot}. Note
that when $x_{\text{max}}=c$, both these conditions become $\alpha
\le 1$ which is more restrictive than the local stability bound of
$\alpha<\pi/2$. The results presented in the next chapter attempt to
eliminate the gap between local and global stability results by
showing global stability for $\alpha < \pi/2 f_1(x_{\text{max}}/c)$.

%%%%%%%%%%%%%%%%%%%%%%%%%%%%%%%%%%%%%%%%%%%%%%%%%%%%%%%%%%%%%%%%%%%%%%%%%%%%%%

\chapter{Stability of Internet Congestion
Control}\label{chp:IQCresult}

\section{Introduction}

In this chapter, we extend some of the linear stability results of
Paganini et al.~\cite{paganini_2001} discussed in the previous
chapter to the dynamics of a nonlinear implementation. Many
algorithms have been proposed for Internet congestion control, some
of which have been shown to be globally stable in the presence of
delays, nonlinearities and discontinuities. These proofs can be
grouped into several categories according to methodology. In
particular, Lyapunov-Razumhikin theory has been used to show global
stability
in~\cite{ying_2004,deb_2003,hollot_2001,wang_2002,wang_2003},
Lyapunov-Krasovskii functionals have been used to show global
stability
in~\cite{alpcan_2003,mazenc_2003,papachristodoulou_2004a,papachristodoulou_2004b}
and an input-output approach was taken in~\cite{wang_2002,fan_2003}.
In all of these cases, stability has been proven with varying
degrees of conservatism with respect to restrictions on system
parameters or delays.

In general, stability analysis of nonlinear, discontinuous
differential equations with delay is quite difficult. Although
frequency domain techniques have been shown to be effective when
applied to linear systems with delay, these tools fail in the
presence of nonlinearity. In addition, although time-domain analysis
of nonlinear finite dimensional systems has had some success,
analysis of the infinite dimensional systems associated with delay
has been more problematic. In this thesis, we are able to obtain
improved results by decomposing the nonlinear, discontinuous,
delayed system into an interconnection of a linear system with delay
and a nonlinear, discontinuous system without delay. We analyze the
subsystems separately and prove a passivity result for each. One
benefit of such an approach is that it allows us to use
frequency-domain arguments in addressing the infinite dimensional
linear system. We can then use time-domain arguments in the analysis
of the single state nonlinear system. These individual results are
then be combined using the newly developed generalized passivity
framework of Rantzer and Megretski as outlined in
Chapter~\ref{chp:stabilitytheory}. This approach yields improved
results by allowing us to decompose the original difficult problem
into simpler subproblems, each of which may be solved with less
conservatism. To apply the passivity framework of Rantzer and
Megretski to the internet congestion control problem, we must first
transform what is essentially a question of internal stability into
a problem posed in the input-output framework. Furthermore, because
input-output stability and internal stability are not equivalent,
once a passivity result has been proven, we must use further
analysis to show that this implies asymptotic stability. All these
issues are addressed in the following sections.

\section{Reformulation of the Problem}\label{sec:prelim}
In this section we reformulate the question of internal stability of
the proposed congestion control algorithm as a question of
input-output stability of the interconnection of a linear system
with delay and a nonlinear system without delay. This approach was
motivated, in part, by the work of Wang\cite{wang_2002} and
J\"{o}nsson\cite{jonsson_2000}. If we consider the problem of a
single link and a single source, then from the development in
Chapter~\ref{chp:ICC}, we have that $y(t)=x(t-\tau^f)$ and
$q(t)=p(t-\tau^b)$ where $\tau^f+\tau^b=\tau$. Given an initial
condition $x_0\in \mathcal{C}_\tau$, the dynamics can now be
summarized as $p(t)=x_0(t)$ for $t \in [-\tau,0]$ and the following
for $t\ge 0$.
\begin{align}
\label{eqn:dynamics_p} \dot{p}(t) &=
\begin{cases}
\frac{x_{\text{max}}}{c}e^{-\frac{\alpha}{\tau} p(t-\tau)}-1
& p(t) > 0 \\
\text{max}\{0, \frac{x_{\text{max}}}{c}e^{-\frac{\alpha}{\tau}
p(t-\tau)}-1\} & p(t) \leq 0
\end{cases}
\\
\label{eqn:dynamics_x} x(t)&=x_{\text{max}} e^{-\frac{\alpha}{\tau}
p(t-\tau^b)}
\end{align}
Since the dynamics of Equation~\eqref{eqn:dynamics_p} are decoupled
from those of~\eqref{eqn:dynamics_x} and stability of $x$ follows
from that of $p$, we need only consider stability of
Equation~\eqref{eqn:dynamics_p}. Now consider the equilibrium point
of Equation~\eqref{eqn:dynamics_p}, $p_0=\frac{\tau}{\alpha} \ln
\frac{x_{\text{max}}}{c}$. As is customary, we change to variable
$z$, where $z(t)=p(t)-p_0$ so that the origin is an equilibrium
point.  Now we have $z(t)=x_0(t)-p_0$ for $t \in [-\tau,0]$ and the
following for $t \ge 0$.
\begin{equation}
\dot{z}(t)=
\begin{cases}
e^{-\frac{\alpha}{\tau} z(t-\tau)}-1 & z(t) > -p_0 \\
\text{max}\{0, e^{-\frac{\alpha}{\tau} z(t-\tau)}-1\} & z(t) \leq
-p_0
\end{cases}
\label{eqn:dynamics_z}
\end{equation}
For convenience and efficiency of presentation, we will refer to the
solution map defined by Equation~\eqref{eqn:dynamics_z} as
$A:\mathcal{C_\tau} \rightarrow \mathcal{C}$. Implicit in these
dynamics is the constraint $z(t)\ge -p_0$. If we assume that any
initial condition will satisfy this constraint, we can include the
constraint in the dynamics without altering the solution map. For
convenience, we define the following bounded continuous functions.
\begin{align*}
f_1(y) &= \min \bl\{
e^{\frac{\alpha}{\tau} y}-1, e^{\frac{\alpha}{\tau} p_0}-1\br\} \\
f_2(y) &= \max \bl\{0, f_1(y) \br\} \\
f_c(x,y) &=
\begin{cases}
f_1(y) & \text{if }x > -p_0\\
f_2(y) & \text{otherwise}
\end{cases}
\end{align*}

\begin{figure}[htb]
\centerline{\includegraphics[width=0.6\textwidth]{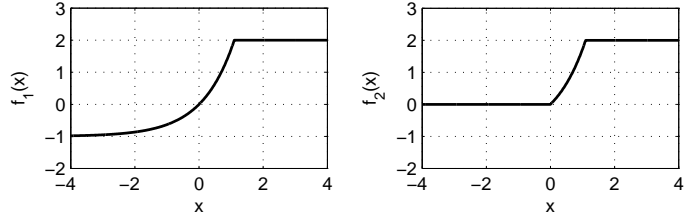}}
\caption{$f_1$ and $f_2$} \label{fig:f1andf2}
\end{figure}

%where the domains are constrained such that $y\geq -p_0$.
These functions are illustrated in Figure~\ref{fig:f1andf2}. We now
have the following equation for $t\ge 0$.

\begin{equation}
\label{eqn:dynamics_z2} \dot{z}(t)=f_c(z(t),-z(t-\tau))
\end{equation}

Before proceeding, we must mention well-posedness of the solution
map, $A$. We use a method of steps. Given any absolutely continuous
solution $z(t)$ on some interval $[T_1, T_1+\tau]$, we observe that
$\dot{z}(t)=f_c(z(t),z(t-\tau))=f_c(z(t),t)$ is a function only of
time and state $z(t)$ for the interval $[T_1+\tau, T_1+2\tau]$. From
boundedness of the $f_c$, continuity with respect to $z(t-\tau)$ and
noting upper semi-continuity of the associated differential
inclusion with respect to $z(t)$, we can now establish via
Fillipov~\cite{filippov_1988}[p77] the existence and uniqueness of a
continuous solution $z(t)$ over the interval $[T_1+\tau, T_1+2
\tau]$. Assuming a continuous initial condition, this implies the
existence and uniqueness of the solution map $A: \mathcal{C}_\tau
\rightarrow \mathcal{C}$.

%%%%%%%%%%%%%%%%%%%%%%%%%%%%%%%%%%%%%%%%%%%%%%%%%%%%%%%%%%%%%%%%%%%%%%%%%%%%%%

\subsection{Separation into subsystems}\label{sec:separation}

Equation~\eqref{eqn:dynamics_z2} is a delay-differential equation
defined by a nonlinear, discontinuous function. To aid in the
analysis, we will reformulate the problem as the interconnection of
two subsystems where the $W_2$ stability on $W_2 \times L_2$
stability of this interconnection implies asymptotic stability on
$X$ of the original formulation for some set $X$. Define the map $G$
by $w=G u$ if
\begin{equation*}
w(t)= \int_{t-\tau}^{t}u(\theta)d\theta.
\end{equation*}
Note that $G$ is a linear operator which can be represented by the
convolution with $g(t)=\text{step}(t)-\text{step}(t-\tau) \in L_1$.
This implies that $\hat{G} \in \mathcal{A}$. Moreover, $G$ can be
represented in the frequency domain by $\hat{G}(s)=\frac{1-e^{-\tau
s}}{s}$ which implies $G$ is a bounded operator on $L_2$ since
$\norm{\hat{G}(\iw )}_\infty = \tau$. In addition, $s \hat{G}(s) \in
\mathcal{A}$ since it can be represented by convolution with
$\delta(t)-\delta(t-\tau)$. Define the map $\Delta_z$ by $z=\Delta_z
y$ if $z(0)=0$ and
\begin{equation*}
\dot{z}(t) = f_c \bl( z(t),y(t)-z(t) \br).
\end{equation*}

We define the map $\Delta$ by $v=\Delta y$ if $v(t)=\dot{z}(t)$
where $z=\Delta_z y$. Addressing well-posedness, if $y\in W_2$, then
$y$ is absolutely continuous on any finite interval(See p.~25 in
J\"{o}nsson~\cite{jonsson_1996}). From boundedness of $f_c$,
continuity with respect to $y(t)$, and upper semi-continuity of the
associated differential inclusion with respect to $z(t)$, we can
again establish the existence and uniqueness of an absolutely
continuous solution $z$ and thus of the map $\Delta_z$.
Well-posedness of $\Delta$ follows immediately. Further properties
of $\Delta$ will be derived in later sections.

If we now form the interconnection of $G$ and $\kappa \Delta$ for
$\kappa \in [0,1]$, as defined above with a single input $f\in W_2$,
we can construct a map from input $f$ to outputs $y,u$. For
convenience and efficiency of presentation, we will denote the
interconnection map for $\kappa=1$ by $B:W_2 \rightarrow W_{2e}
\times L_{2e}$. Furthermore, for $\kappa=1$ we denote the map from
input $f$ to internal variable $z$ by $B_z$. For $t \le 0$, we let
$u(t)=y(t)=z(t)=f(t)=0$ and for $t \ge 0$, the interconnection
dynamics combine as follows.
\begin{align*}
u(t)&=\kappa \dot{z}(t)\\
y(t) &= \int_{t-\tau}^t u(t)dt+f(t)\\
&= \kappa \bl( z(t)-z(t-\tau)\br)+f(t)\\
\dot{z}(t)&=f_c(z(t),y(t)-z(t))\\
&=f_c(z(t),f(t)-\kappa z(t-\tau)-(1-\kappa)z(t))
\end{align*}

As before, from continuity with respect to $f(t)$ and $z(t-\tau)$,
upper semi-continuity with respect to $z(t)$ and boundedness of
$f_c$, we can conclude existence and uniqueness of $z \in L_{2e}$.
Since $z$ has bounded derivative and $z(t)=0$ for $t\le 0$, we now
have that the map $B_z: W_2 \rightarrow W_{2e}$ is well-posed.
Furthermore, this implies that $y \in W_{2e}$ and $u \in L_{2e}$
which yields well-posedness of the interconnection for any $\kappa
\in [0,1]$ and specifically of the map $B: W_2 \rightarrow W_{2e}
\times L_{2e}$.

%%%%%%%%%%%%%%%%%%%%%%%%%%%%%%%%%%%%%%%%%%%%%%%%%%%%%%%%%%%%%%%%%%%%%%%%%%%%%%
\section{Input-Output Stability}
In this section we use the IQC defined by $\Pi_{B}$,
$\lambda=\frac{2}{\pi}$, where
\begin{equation*}
\Pi_B=\bmat{0 & \beta \\ \beta & -\frac{4}{\pi}-2}
\end{equation*}
and $\beta=\alpha/(\alpha_\text{max}\tau)$ to establish $W_2 \times
L_2$ stability on $W_2$ of the interconnection for any $\tau \ge 0$,
$ 0 < \alpha < \pi/2 \alpha_\text{max}$. Here we define
\[ \alpha_\text{max}= \ln(x_{\text{max}}/c)/ \bl((x_{\text{max}}/c)-
1\br).
\]

%We also show that $W_2 \times L_2$ stability on $W_2$ of the
%interconnection implies asymptotic stability of the original
%formulation of the congestion control protocol.

\subsection{$\Delta$ satisfies the IQC}\label{sec:delta_IQC}

In this subsection we show that if $\alpha>0$, then $\Delta$ and
consequently $\kappa\Delta$ are bounded and satisfy the IQC defined
by $\Pi_B,\lambda=\frac{2}{\pi}$ for all $\kappa \in [0,1]$. The
methods used in this subsection were motivated by those in
J\"{o}nnson~\cite{jonsson_2000} and Wang~\cite{wang_2002}. For
$\gamma=4 \beta/ \pi
> 0$, we prove the following for all $y \in W_2$, $v=\Delta y$.
\begin{align*}\frac{1}{2
\pi}\int_{-\infty}^{\infty}\bmat{\hat{y}(\iw)\\\hat{v}(\iw)}^*
\bmat{0 & \beta \\ \beta & -\frac{4}{\pi}-2}
\bmat{\hat{y}(\iw)\\\hat{v}(\iw)} + \frac{4}{\pi}
\ip{v}{\dot{y}}\\
\ge -\gamma|y(0)|^2
\end{align*}
By Parseval's equality, this is equivalent to
\begin{align*}
\frac{2}{\pi}\ip{v}{\dot{y}-v}+\ip{v}{\beta y-v } \geq
-\frac{\gamma}{2} |y(0)|^2
\end{align*}

\begin{figure}[htb]
\centerline{\includegraphics[width=0.45\textwidth]{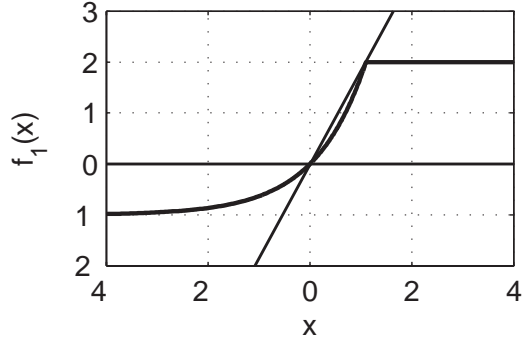}}
\caption{The nonlinearity $f_1$ satisfies a sector bound}
\label{fig:sector}
\end{figure}

A critical result used in the analysis of this section is the
existence of a sector bound on the nonlinearity $f_1$ and
consequently on $f_2$.

\begin{lem}$0 \le f_i(x)x \le \beta x^2$ for $i=1,2$ where
$\beta=\frac{e^{\frac{\alpha}{\tau} p_0}-1}{p_0}$.
\end{lem}
See Appendix~\ref{app:IQCresult} for Proof.

This key feature is illustrated in Figure~\ref{fig:sector}.

\begin{lem}\label{lem:IQC1}
\label{lem:z_in_L2} If $v=\Delta y$ with $y\in W_{2}$, then
\begin{enumerate}
\item $v \in L_2$ with norm bound $\beta \norm{y}$,
\item $\ip{v}{\beta y-v} \ge 0$
\end{enumerate}
\end{lem}
\begin{proof}
We start by noting the following.
\[
f_c(x,y) =
\begin{cases}
f_1(y) & \text{if } x > -p_0 \text{ or } y \geq 0 \\
0 & \text{otherwise}
\end{cases}
\]
As a consequence of the above sector bounds, we have
\[
f_c(x,y)^2 \leq  \beta y f_c(x,y).
\]
Let $z=\Delta_z y$, then this implies
\begin{align*}
\dot{z}(t)^2    & =f_c\bl( z(t), y(t)-z(t) \br)\dot{z}(t)\\
&\leq \beta \bl(y(t)-z(t)\br) \dot{z}(t)
\end{align*}
Now for any $T\ge 0$, we have
\begin{align}
\norm{P_T v}^2
&=\int_{0}^{T}v(t)^2 dt = \int_{0}^{T}\dot{z}(t)^2 dt \notag \\
&\le \beta \int_{ 0}^{T} \dot{z}(t)(y(t)-z(t)) dt \notag \\
&= \beta \int_{ 0}^{T} \dot{z}(t)y(t)dt -\frac{\beta}{2} (z(T)^2- z(0)^2) \notag \\
&\le \beta \ip{P_T \dot{z} }{ y} \label{eqn:analysis_3} \\
&\le \beta \norm{P_T \dot{z}} \norm{y} = \beta \norm{P_T v} \norm{y}
\notag
\end{align}
Therefore, $\norm{P_T v} \le \beta \norm{y}$ for all $T \ge 0$. Thus
$v \in L_2$ with norm bounded by $\beta \norm{y}$. Statement 2
follows from line~\ref{eqn:analysis_3} by letting $T \rightarrow
\infty$.
\end{proof}
%%%%%%%%%%%%%%%%%%%%%%%%%%%%%%%%%%%%%%%%%%%%%%%%%%%%%%%%%%%%%%%%%%%%
\begin{lem}
\label{lem:z->0} Let $z=\Delta_z y$ with $y \in W_{2}$, then
$\lim_{t \rightarrow \infty } z(t)=0$.
\end{lem}

See Appendix~\ref{app:IQCresult} for Proof.

\begin{lem}\label{lem:IQC2}
If $v=\Delta y$ with $y \in W_{2}$, then $\ip{v}{\dot{y}-v} \ge
-\beta|y(0)|^2$.
\end{lem}

See Appendix~\ref{app:IQCresult} for Proof.

\begin{lem}\label{lem:IQC} $\kappa \Delta$ satisfies the IQC
defined by $\Pi_B,\lambda=\frac{2}{\pi}$ for any $\kappa \in [0,1]$
\end{lem}

\begin{proof}By Lemmas~\ref{lem:IQC1} and~\ref{lem:IQC2}, we have
the following.
\begin{align*}
&\frac{1}{2
\pi}\int_{-\infty}^{\infty}\bmat{\hat{y}(\iw)\\\hat{v}(\iw)}^*
\bmat{0 & \beta \\ \beta & -\frac{4}{\pi}-2}
\bmat{\hat{y}(\iw)\\\hat{v}(\iw)} + \frac{4}{\pi} \ip{v}{ \dot{y}}
\\
&=\frac{4}{\pi}\ip{v}{\dot{y}-v}+2\ip{v}{\beta y-v }
\\
&\ge - \frac{4\beta}{\pi}|y(0)|^2
\end{align*}

We conclude as a consequence that $\kappa \Delta$ satisfies the IQC
defined by $\Pi_B,\lambda=\frac{2}{\pi}$ for any $\kappa \in [0,1]$,
since
\begin{align*}
&\frac{2}{\pi}\ip{\kappa v}{\dot{y}-\kappa v}+\ip{\kappa
v}{\beta y-\kappa v } \\
&\ge \kappa\bl( \frac{2}{\pi}\ip{ v}{\dot{y}-v}+\ip{v}{\beta y-v }\br) \\
&\ge - \kappa \frac{2 \beta}{\pi}|y(0)|^2 \ge-
\frac{2\beta}{\pi}|y(0)|^2
\end{align*}

\end{proof}
%To summarize this section, we have shown that $\Delta$ is bounded
%and that for any $y\in W_2$ and $v=\Delta y$, we have $\ip{v}{\beta
%y-v} \ge 0$ and $\ip{v}{\dot{y}-v} \ge -\beta |y(0)|^2$. Therefore,
%we conclude that $\Delta$ satisfies the IQC defined by $\Pi_B,
%\lambda=\frac{2}{\pi}$, since
%
%\begin{align*}
%\frac{1}{2
%\pi}\int_{-\infty}^{\infty}\bmat{\hat{y}(\iw)\\\hat{v}(\iw)}^*
%\bmat{0 & \beta \\ \beta & -\frac{4}{\pi}-2}
%\bmat{\hat{y}(\iw)\\\hat{v}(\iw)} + \frac{4}{\pi} \ip{v}{ \dot{y}}
%\\
%\ge - \frac{4\beta}{\pi}|y(0)|^2
%\end{align*}

\subsection{Properties of $G$}\label{sec:G_IQC}

Recall that we define the map $G$ as follows. $w=G u$ if
\[w(t)=
\int_{t-\tau}^{t}u(\theta)d\theta.
\]
\begin{lem}\label{lem:condition4}Suppose $0 < \alpha <\pi
\alpha_\text{max}/2$. Then $\hat{G}$ satisfies
condition~\ref{cdn:condition4} of Theorem~\ref{thm:stability_thm}.
\end{lem}
\begin{proof}
Recall that $G$ can be represented as a transfer function in the
frequency domain by $\hat{G}=\frac{1-e^{-\iw \tau}}{\iw}$. Now,
examine the term
\begin{align*}
&\bmat{\hat{G}(\iw)\\I}^*\left(\Pi_B+\bmat{0 & \lambda \iw^* \\ \lambda \iw & 0 }\right)\bmat{\hat{G}(\iw)\\I}\\
&=\bmat{\frac{1-e^{-\iw \tau}}{\iw}\\1}^*\bmat{0 & \beta+\frac{2}{\pi}\iw^* \\
\beta+\frac{2}{\pi}\iw & -\frac{4}{\pi}-2}\bmat{\frac{1-e^{-\iw \tau}}{\iw} \\1}\\
&=2\cdot \text{Real} \left(\beta\frac{1-e^{-\iw
\tau}}{\iw}-\frac{2}{\pi} e^{-\iw
\tau}-1\right) \\
&= 2\left(\beta \tau \frac{\sin (\omega \tau)}{\omega
\tau}-\frac{2}{\pi} \cos(\omega \tau)-1\right)=2 p(\omega \tau)
\end{align*}
\begin{figure}[htb]
\centerline{\includegraphics[width=0.35\textwidth]{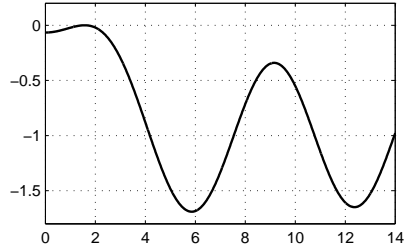}}
\caption{Plot of $p_0(\omega)=\frac{\pi}{2}
\frac{\sin(\omega)}{\omega}-\frac{2}{\pi} \cos(\omega ) -1$ vs.
$\omega$.} \label{fig:plot_of_fn}
\end{figure}
Define $p_0(\omega)=p(\omega)$ for $\beta \tau=\pi/2$. The plot of
$p_0(\omega)=\frac{\pi}{2} \frac{\sin
(\omega)}{\omega}-\frac{2}{\pi} \cos(\omega)-1 $ is given in
Figure~\ref{fig:plot_of_fn}. As one can see, the function is
non-positive near the origin. In fact, for any $\beta \tau <\pi/2 $,
there exists an $\eta>0$ such that $p(\omega)<-\eta$ for all
$\omega$. To see this, consider the following domains.

\paragraph{$\omega > 2 \pi$:} We have that $|
\frac{\sin (\omega)}{\omega} | \le 1/ 2 \pi$ for $|\omega| >2 \pi$
and $-1-\frac{2}{\pi}\cos (\omega))\le -1+2/ \pi $ for all $\omega$.
Therefore, $p(\omega) \le \frac{1}{4} -(1+\frac{2}{\pi}\cos
(\omega))\le \frac{2-3\pi/4}{\pi}<-.1$.

\paragraph{$\pi-.1 \le \omega \le 2 \pi$:} For the interval $[\pi-.1,2\pi]$,
we have $\frac{\sin (\omega)}{\omega} \le .04$. Therefore $p(\omega)
\le -1+2/\pi+.04\frac{\pi}{2} < -.3$.

\paragraph{$0 \le \omega < \pi-.1$:} One can see from plot of $p_0$
in Figure~\ref{fig:plot_of_fn} that $p_0(\omega)\le 0$ on the
interval $[0,\pi-.1]$. Since $\beta \tau <\pi/2$, we can let
$\epsilon = \pi/2 - \beta \tau > 0$. Then
$p(\omega)=p_0(\omega)-\epsilon \frac{\sin (\omega)}{\omega}\le
-\epsilon \frac{\sin (\omega)}{\omega} < -.2 \epsilon $ on
$[0,\pi-.1]$.

Therefore, for any $\beta \tau < \frac{\pi}{2}$, let $\eta=\min\{.1
, .2(\pi/2-\beta \tau)\}$. Then $p(\omega) < -\eta$ for all $\omega
\in \R$. Since $\beta \tau = \tau(e^{\frac{\alpha}{\tau}p_0}-1)/p_0
=\alpha/\alpha_\text{max}$, if $\alpha < \pi\alpha_\text{max}/2$, we
have that $\beta \tau < \pi/2$, and hence if $0 < \alpha <\pi
\alpha_\text{max}/2$, condition~\ref{cdn:condition4} of
Theorem~\ref{thm:stability_thm} is satisfied.
\end{proof}
\subsection{Stability of the Interconnection}\label{sec:stable_intercon}

We conclude our discussion of input-output stability with the
following Theorem concerning the stability of the the
interconnection of $\Delta$ and $G$.

\begin{thm}
Suppose $\alpha \in (0 , \pi\alpha_\text{max}/2)$. Then the map $B$,
defining the interconnection of $\Delta$ and $G$ is $W_2 \times L_2$
stable on $W_2$.
\end{thm}

\begin{proof}
We have shown that $G$ is a linear causal bounded operator with
$\hat{G}(s),s\hat{G}(s)\in \mathcal{A}$. We have also shown that the
interconnection of $G$ and $\kappa \Delta$ is well-posed for all
$\kappa \in [0,1]$ and Lemma~\ref{lem:IQC} shows that $\kappa
\Delta$ satisfies the IQC defined by $\Pi_B, \lambda=\frac{2}{\pi}$
for all $\kappa \in [0,1]$.  Finally, by Lemma~\ref{lem:condition4}
we have that for all $\alpha \in (0 , \pi\alpha_\text{max}/2)$,
condition~\ref{cdn:condition4} of Theorem~\ref{thm:stability_thm} is
satisfied. We can therefore use Theorem~\ref{thm:stability_thm} to
prove $W_2 \times L_2$ stability on $W_2$ of the interconnection for
any $\alpha \in (0 , \pi\alpha_\text{max}/2)$.
\end{proof}

\section{Asymptotic Stability}\label{sec:asymptotic}
In this section, we show that $W_2 \times L_2$ stability on $W_2$ of
the interconnection implies asymptotic stability of the original
formulation of the congestion control protocol. Recall that the
original solution map $A$ is defined by the following.
\begin{align}
\dot{z}(t)&=f_c(z(t),-z(t-\tau)) \qquad &t \ge 0  \\ z(t) &= x_0(t)
\qquad \qquad &t \in [-\tau,0]  \label{eqn:statedyn}
\end{align}

The interconnection maps $B$ and $B_z$, however, are defined by the
following differential equation.
\begin{align}
\dot{z}(t)&=f_c(z(t),f(t)-z(t-\tau)) \qquad &t \ge 0\label{eqn:interdyn}\\
z(t)&=0 \qquad \qquad &t\le 0
\end{align}

In the previous section, we have proven $W_2 \times L_2$ stability
on $W_2$ of the map $B$ where $B$ represents a reformulation of the
problem in the input-output framework. We would like to show,
however, that for some $X \subset \mathcal{C}_\tau$ this result also
implies asymptotic stability on $X$ of the solution map $A$, where
$A$ represents the original formulation of the problem. This is done
in the following Theorem.

\begin{thm}
Suppose $\alpha \in (0 , \pi\alpha_\text{max}/2)$. Then the
delay-differential equation~\eqref{eqn:dynamics_p} describing the
algorithm proposed by Paganini et al.~\cite{paganini_2001} is
asymptotically stable on $X=\{x: x \in W_2 \cap \mathcal{C}_\tau,
x(t) \ge -p_0, x(-\tau)>-p_0 \}$.
\end{thm}

\begin{proof}
We have already shown $W_2 \times L_2$ stability on $W_2$ of the map
$B$. Let $x_0 \in X$ be an arbitrary initial condition.
Theorem~\ref{thm:AppendixA} in the Appendix states that for any
initial condition $x_0 \in X$, there exists some $f \in W_2$ and
$T>0$ such that $A(x_0,t)=B_z(f,t+T)$ for all $t \ge 0$. Let
$(y,u)=B(f)$, then $y \in W_2$. Furthermore, recall $B_z f=\Delta_z
y$ where $B_z$ is the map to internal variable $z$. By
Lemma~\ref{lem:z->0}, if $y \in W_2$, then $\lim_{t \rightarrow
\infty}\Delta_z(y,t)=0$. Therefore, $\lim_{t \rightarrow
\infty}A(x_0,t)=\lim_{t \rightarrow \infty}B_z(f,t+T)=\lim_{t
\rightarrow \infty}\Delta_z(y,t+T)=0$.
\end{proof}

\section{Implementation}\label{sec:implementation}

To implement the proposed algorithm in the Internet framework, the
window size can be adjusted to deliver the required packet rate as
given by equation~\eqref{eqn:source_rate}.  In implementation, the
delay is unlikely to be known. In this case, a bound on the expected
delay size can be used. Overestimation of the delay will result in
an increased stability margin.

Modification of the link can take many forms. Price information from
the link must be fed back to the source. Since queues themselves
integrate excess rate, price of a congested resource can be computed
directly using the queueing delay. However, this approach results in
non-empty equilibrium queues and the possibility of unmodeled
dynamics due to variable queueing delays. If a link instead uses a
virtual capacity to avoid non-empty equilibrium queues, then
explicit integration of incoming packets would be required and
another mechanism must be used to feed back price information. An
example of direct feedback of price information using packet marking
is given by ECN. In one of the proposed implementations, packets are
randomly marked at each link with probability $1-\phi^{-p_j(t)}$ for
some fixed $\phi>1$. Thus, assuming no duplications, if $\nu$ is the
percentage of marked packets received at the source, then the
aggregate price can be measured as $q_i(t)= -\frac{\log
(1-\nu)}{\log (\phi)}$. This variant is known as random exponential
marking.

\section{Conclusion} \label{sec:conclusion}
To summarize, for the case of a single source with a single link, we
have demonstrated both input-output and asymptotic stability. We
have used a generalized passivity framework to decompose a difficult
nonlinear, discontinuous, infinite-dimensional problem into separate
linear, infinite-dimensional and nonlinear, finite-dimensional
subproblems, each of which is amenable to existing analysis
techniques. A key feature of the analysis of the nonlinear subsystem
was the existence of a sector bound on the nonlinearity.

In addition to the case of a single source with a single link, the
result presented in this chapter applies directly to the case of
multiple sources with identical fixed delay. The proof may also be
easily adapted to alternate implementations of the proposed
linearized protocols so long as they admit a similar sector bound on
the nonlinearity. Although we would like to have proven a result in
the case of multiple heterogeneous delays and arbitrary topology, we
have as yet not been able to construct an appropriate sector bound.
The root of the problem seems to lie in the discontinuity due to the
projection used in the gradient projection algorithm.

\chapter{Sum-of-Squares and Convex
Optimization}\label{chp:SOStheory}
%Just as semidefinite programming introduced a revolution in the
%analysis and control of linear finite-dimensional systems, so has
%sum-of-square(SOS) optimization started an avalanche of results in
%the analysis of nonlinear problems. The basic premise is to answer
%difficult, often NP-hard, questions through the use of a sequence of
%computationally tractable sufficient conditions, of increasing
%accuracy and complexity, and which approach necessity. At the core
%of all these results is a semidefinite programming approach to
%question of polynomial non-negativity first introduced in the thesis
%by Pablo Parrilo~\cite{parrilo_2000}. This topic is discussed in the
%first section of this chapter for scalar polynomials and extended to
%matrices of polynomials in the subsequent section.

\section{Convexity and Tractability}
Consider the following optimization problem for $f_i \in \R[x]$ for
$i=0 \ldots n_f$.

\begin{align*}
&\max f_0(x): \\
&f_i(x) \ge 0, \qquad i=1, \ldots, n_f
\end{align*}

The problem as formulated is clearly, in general, not convex. In
fact, many $NP$-hard problems can be formulated in this manner.
However, if we define $Y:=\{x:f_i(x)\ge 0,\; i=1\ldots n_f\}$, then
the problem can be recast as an equivalent convex optimization
problem as follows.

\begin{align*}
&\min \gamma: \\
&\gamma - f_0(x) \in \mathcal{P}^Y
\end{align*}

Although this reformulation is convex, we must conclude that it too
is intractable, unless we believe that $P=NP$. The point that we
make with this section is that convexity, by itself, is not a
sufficient condition for tractability of a problem. What we need in
addition is an efficient affine test for membership in the convex
cone. Construction of such membership tests for various convex sets
is the goal of Chapters~\ref{chp:SOStheory},~\ref{chp:LinearCase},
and~\ref{chp:nonlinearcase}.

\section{Sum-of-Squares Decomposition}\label{sec:SOSDecomp}
Any polynomial, $f\in \R[x]$, of degree $d$ and $n$ variables can be
expressed as the linear combination of $n+d \choose d$ monomials.

\begin{equation}
f(x)=\sum_{i=1}^{n+d \choose d} c_i
  x_1^{\gamma_{i,1}}\cdots x_n^{\gamma_{i,n}} \quad \gamma_{i,j} \in \Z^+
\label{eqn:polynomial}
\end{equation}

The question of whether $f \in \mathcal{P}^+$, that is, $f(x)\ge 0$
for all $x\in \R^n$, is NP hard for polynomials of degree $4$ or
more. Thus there is unlikely to exist a computationally tractable
set membership test for $\mathcal{P}^+$ unless $P=NP$. However,
there may exist some convex cone $\Sigma \subset \mathcal{P}^+$ for
which the set membership test is computationally tractable. One such
cone is defined as follows.

\begin{defn}
A polynomial $s \in \R[x]$ satisfies $s\in \Sigma_s$ if it can be
represented in the following form for some finite set of polynomials
$g_i\in \R[x]$, $i=1 \ldots m$.
\[s(x)=\sum_{i=1}^m g_i(x)^2
\]
An element $s \in \Sigma_s$ is referred to as a
\eemph{sum-of-squares polynomial}. The set $\Sigma_s^d$ is defined
as the elements of $\Sigma_s$ of degree $d$ or less.
\end{defn}

\begin{lem} $\Sigma_s \subset \mathcal{P}^+$
\end{lem}

In~\cite{parrilo_2000}, it was shown that the set membership test $f
\in \Sigma_s$ is computationally tractable.

\begin{lem}\label{lem:SOS} A degree $2d$ polynomial $f \in \R[x]$ satisfies $f \in
\Sigma_s$ if and only if there exists some matrix $Q \in \S^{n_z}$
where $n_z= {n+d \choose d}$, such that $Q \ge 0$ and
\[f(x)=Z_d[x]^T Q Z_d[x].
\]
\end{lem}

The equality constraint in Lemma~\ref{lem:SOS} is affine in the
monomial coefficients of $f$. Therefore, set membership in
$\Sigma_s$ can be tested using semidefinite programming. The
question of how well $\Sigma_s$ represents $\mathcal{P}^+$ has been
a topic on ongoing research for some time. It was conclusively shown
through use of the Motzkin polynomial that in general $\Sigma_s \ne
\mathcal{P}^+$. However, if we define $\Sigma_r$ to be the convex
cone of sums of squares of quotients of polynomials, then Artin
showed that $\Sigma_r=\mathcal{P}^+$. Not surprisingly, however,
there is no tractable set membership test for $\Sigma_r$. As an
interesting corollary of this result, it was shown that $f \in
\mathcal{P}^+$ if and only if there exists some $g \in \Sigma_s$
such that $fg \in \Sigma_s$. This does not constitute a tractable
set membership test, unfortunately, since no bound is given on the
degree of $g$. In addition, the introduction of the auxiliary
variable $g$ destroys the convexity of the constraint since the
resulting feasibility problem is bilinear. Building on this work,
however, Reznick showed that if $f \in \mathcal{P}^+$, then there
exists some $d \in Z^+$ such that $(x^T x)^d (f(x)+1) \in \Sigma_s$.
This result provides a test for membership of $f$ in the interior of
$\mathcal{P}^+$, but not a tractable one, as there are, in general,
still no bounds on $d$. What the result of Reznick does provide,
however, is a sequence of tractable tests for membership in the
interior of $\mathcal{P}^+$, indexed by the integer $d$.

%It has been shown~\cite{reznick_1996} that if $f(x)>0$ for all $x
%\in \R^n$, then there exists an $m\ge0$ and a \eemph{sum of squares}
%(SOS) polynomial $s$ such that $(x^T x)^m f(x)=s(x)$. Additionally,
%in the case of 1 or 2 variables, quartic polynomials in three
%variables and all quadratic polynomials, existence of a SOS
%representation is also necessary for global non-negativity.
%
%In~\cite{parrilo_2000}, it has been shown that the existence of a
%bounded degree SOS representation of a polynomial is equivalent to a
%semidefinite program with equality constraints expressed in terms of
%the monomial coefficients in Equation~\eqref{eqn:polynomial}.More
%specifically, if $\deg(f)$ is even, let $z$ be the vector of
%monomials of degree less than or equal to $\deg(f)/2$. Then $f$ is
%SOS if and only if there exists a $Q \ge 0$ such that $f(x)=z^T Q
%z$. Equating coefficients gives a set of affine constraints on $Q$.
%In this manner, SOS constraints can be used in the place of
%non-negativity constraints, and can be tested using semidefinite
%programming.

\section{Positivstellensatz Results}\label{sec:PSresults}

Membership in the cone $\mathcal{P}^K$ is harder to approximate than
$\mathcal{P}^+$. Consider a region $\hat{K_f}$ of the following form
for $f_i \in \R[x]$, $i=1, \ldots, n_K$.
\[K_f:=\{x \in \R^n : f_i(x) \ge 0, \; i=1, \ldots, n_K\}.
\]
Then $f_0 \in \mathcal{P}^{K_f}$ if and only if the following set is
empty.

\[K_f^*:=\{x \in \R^n: -f_0(x)\ge 0 ,\; f_0(x)\ne 0,\; f_i(x) \ge 0 ,\; i=1, \ldots, n_K\}
\]

A condition for emptiness of the set $K_f^*$ comes from the
following simplified Positivstellensatz result from
Stengle~\cite{stengle_1974}.

\begin{thm}\label{thm:Stengle}
Let $K_f$ be given as above and let $I$ denote the set of subsets of
$\{0,\ldots,n_K\}$. Then $-f_0 \in \mathcal{P}^{K_f}$ if and only if
there exist $s \in \Sigma_s$, $s_J \in \Sigma_s$ for $J \in I$ and
$k \in \Z^+$ such that
\[s(x)+\sum_{J\in I} s_J(x)\prod_{i \in J}f_i(x)+ f_0(x)^{2k} =0
\]
\end{thm}

%\begin{thm}\label{thm:Stengle} Given $f_i \in \R[x]$, $i=0 \ldots n_K$ and $g \in \R[x]$, define
%\[K^*:=\{x \in \R^n : g(x) \ne 0, \; f_i(x) \ge 0, \; i=0 \ldots n_K\}.
%\]
%Let $I$ denote the set of subsets of $\{0,\ldots,n_K\}$. Then $K^* =
%\emptyset$ if and only if there exist $s_J \in \Sigma_s$, $J \in I$
%and $k \in \Z^+$ such that
%\[\sum_{J\in I} s_J(x)\prod_{i \in J}f_i(x)+ g(x)^{2k} \in \Sigma_s
%\]
%\end{thm}

The Positivstellensatz condition given by Stengle is not tractable,
due to the introduction of bilinear terms $s_J f_0$ and since no
degree bound exists on either $k$ or the degrees of the $s_J$. In
the work by Schm\"{u}dgen, the convexity problem was addressed by
considering the case when $K_f$ is compact. This result can be
stated as follows.

\begin{thm}\label{thm:schmudgen}
Let $\hat{I}$ be the set of subsets of $\{1,\ldots,n_K\}$ and
suppose $K_f$, as defined above, is compact. Suppose $f$ lies in the
interior of $\mathcal{P}^{K_f}$. Then there exist $s_J \in \Sigma_s$
for $J \in \hat{I}$, such that the following holds.
\[f(x)-\sum_{J\in \hat{I}} s_J(x)\prod_{i \in J}f_i(x)\in \Sigma_s
\]
\end{thm}

Theorem~\ref{thm:schmudgen} differs from Stengle's result by
essentially allowing us to set $s_J=0$ when $0 \not \in J$. This
means that the resulting constraint is convex. However,
Schm\"{u}dgen's result still does not constitute a tractable
condition since there exist no bounds on the degrees of the $s_J$.
The theorem does, however, suggest a method for constructing a
sequence of sufficient conditions for membership in the set
$\mathcal{P}^{K_f}$ by considering elements of $\Sigma_s$ of bounded
degree.

\begin{defn}
For a given degree bound, $d$, and compact semi-algebraic set,
$K_f$, of the form given above and $\hat{I}$ as defined above, we
define the convex cone $\Lambda_d^{K_f}$ of functions $f \in \R[x]$
such that there exist functions $s_J \in \Sigma^d_s$ for $J \in
\hat{I}$, of degree $d$ or less, such that the following holds.
\[f(x)-\sum_{J\in \hat{I}} s_J(x)\prod_{i \in J}f_i(x)\in \Sigma_s
\]
\end{defn}

\begin{lem} For any $d\in \Z^+$ and any $K_f$, $\Lambda_d^{K_f} \subset \mathcal{P}^{K_f}$
\end{lem}

Testing membership in $\Lambda_d^{K_f}$ is a tractable problem for
any fixed $d$ and $K_f$ and thus provides a tractable sufficiency
test for membership in $\mathcal{P}^{K_f}$. However, the complexity
of the test scales poorly in $n_k$. This problem was addressed to
some extent by the work of Putinar~\cite{putinar_1993}, which showed
that under certain additional restrictions on the set $K_f$, that we
can assume that $s_J=0$ for all $J$ except $J=1, \ldots, n_K$. This
theorem is stated as follows.

\begin{defn} We say that $f_i \in \R[x]$ for $i=1,\ldots, n_K$ define a \eemph{P-compact}
set $K_f$, if there exist $h\in \R[x]$ and $s_i\in \Sigma_s$ for
$i=0,\ldots,n_K$ such that the level set $\{x \in \R^n : h(x)\ge
0\}$ is compact and such that the following holds.
\[h(x)-\sum_{i=1}^{n_K}s_i(x)f_i(x) \in \Sigma_s
\]
\end{defn}

The condition that a region be P-compact may be difficult to verify.
However, some important special cases include:
\begin{itemize}
\item Any region $K_f$ such that all the $f_i$ are linear.
\item Any region $K_f$ defined by $f_i$ such that there
exists some $i$ for which the level set $\{x:f_i(x)\ge 0\}$ is
compact.
\end{itemize}

Thus any compact region can be made P-compact by inclusion of a
redundant constraint of the form $f_{n_K+1}:=\beta-x^T x$ for
sufficiently large $\beta$. However, for such an inclusion, we must
have some explicit bound for the size of the original region.

\begin{thm}\label{thm:putinar}
Suppose $K_f$, as defined above, is P-compact. Suppose $f$ lies in
the interior of $\mathcal{P}^{K_f}$. Then there exist $s_i \in
\Sigma_s$ for $i=1 ,\ldots, n_K$, such that the following holds.
\[f(x)-\sum_{i=1}^{n_K} s_i(x) f_i(x)\in \Sigma_s
\]
\end{thm}

Theorem~\ref{thm:putinar} is similar to the Positivstellensatz
result of Schm\"{u}dgen in that it does not represent a tractable
test for membership in $\mathcal{P}^{K_f}$ since no bounds exist for
the degree of the $s_i$. However, Putinar's Positivestellensatz
provides a basis for a less computationally intensive sequence of
tractable sufficient conditions.

\begin{defn}
For a given degree bound, $d$, and semi-algebraic set, $K_f$, of the
form given above, we define the convex cone $\Upsilon_d^{K_f}$ of
functions $f \in \R[x]$ such that there exist $s_i \in \Sigma^d_s$
for $i =1, \ldots, n_K$, such that the following holds.
\[f(x)-\sum_{i=1}^{n_K} s_i(x) f_i(x)\in \Sigma_s
\]
\end{defn}

\begin{lem} For any $d\in \Z^+$ and any $K_f$, $\Upsilon_d^{K_f} \subset \mathcal{P}^{K_f}$
\end{lem}

Membership in the set $\Upsilon_d^{K_f}$ is a tractable problem for
any $d\in \Z^+$ and $K_f$ and can be tested using semidefinite
programming. Moreover, for any fixed $d$, the test associated with
$\Upsilon_d^{K_f}$ is less computationally difficult than that
associated with $\Lambda_d^{K_f}$, especially for large values of
$n_K$. However, such statements should be taken with a grain of
salt, since, for a fixed $d$, $\Lambda_d^{K_f}$ may well provide a
better approximation of $\mathcal{P}^{K_f}$ than $\Upsilon_d^{K_f}$.

%Let $k=1$, $s_J=0$ for $0 \not \in J$, $s_J=1$ for $0 \in J$ and $J
%\neq \{0\}$, and $s=0$. Then we have the following condition.
%
%\begin{lem}\label{lem:Stenglerelaxed}
%Let $\hat{I}$ be the set of subsets of $\{1,\ldots,n_K\}$. Then $f
%\in \mathcal{P}^{K_f}$ if there exists a $s_0 \in \Sigma_s$ such
%that the following holds.
%\[f(x)-\sum_{J\in \hat{I}} \prod_{i \in J}f_i(x)\in \Sigma_s
%\]
%\end{lem}
%
%For any given $f_i$, Lemma~\ref{lem:Stenglerelaxed} provides a
%tractable sufficient condition for membership in
%$\mathcal{P}^{K_f}$. Unfortunately, the number of such variables
%scales poorly with $n_K$.
%
% Now suppose we consider the question of whether a polynomial
%$f$ satisfies $f(x)\ge 0$ for all $x \in X$, where $X :=\{x : p_i(x)
%\ge 0, \; i=1,\ldots, k\}$ with scalar polynomials
%$\{p_i\}_{i=1}^k$. The following lemma gives a sufficient condition
%for non-negativity of $f$ on $X$.
%
%\begin{lem}Suppose that for a polynomial $f$ there exist sum of squares
%polynomials $s_i$, for $i=0\ldots k$ such that
%\[f(x)=s_0(x)+\sum_{i=1}^k p_i(x)s_i(x).
%\]
%Then $f(x)\ge 0$ for all $x \in X$.
%\end{lem}
%
%Positivstellensatz results from Putinar\cite{putinar_1993} and
%Schm\"{u}dgen\cite{schmudgen_1991} show that the above condition is
%necessary when $f$ is strictly positive and the $p_i$ satisfy
%specific conditions.

\section{Applications}
In this section, we demonstrate some applications of the techniques
provided in Sections~\ref{sec:SOSDecomp} and~\ref{sec:PSresults} to
stability of nonlinear ordinary differential equations. Our concepts
of stability are the same as those defined for functional
differential equations, albeit with a state space defined in $\R^n$
as opposed to $\mathcal{C}_\tau$.

\subsection{Stability of Ordinary Differential
Equations}\label{sec:nonlinear_ODE} Consider a time-invariant
ordinary differential equation of the following form.
\[\dot{x}(t)=f(x(t))
\]
Here $x(t)\in \R^n$ and $f$ is a polynomial of degree $m$ such that
$f(0)=0$.
\begin{thm}\label{thm:Lyapunov_global}
The system defined by the polynomial function $f$ is globally stable
if there exists an $\alpha>0$ and a Lyapunov function, $V(x)$, with
continuous derivatives such that $V(0)=0$ and the following holds.
\begin{align*}
V(x)&\ge \alpha x^T x \qquad &\text{for all }x \in \R^n\\
\nabla V(x)^T f(x) &\le 0 \qquad &\text{for all }x \in \R^n
\end{align*}
If, in addition, $V(x)^T f(x) \le -\alpha x^T x$ for all $x \in
\R^n$, then the system is globally asymptotically stable.
\end{thm}

By using the methods of the Section~\ref{sec:SOSDecomp}, we can pose
the question of existence of Lyapunov function proving global
stability of the system as a convex feasibility problem in the cone
$\mathcal{P}^+$. Then, by approximating $\mathcal{P}^+$ by
$\Sigma_s$, we can construct a sequence of sufficient conditions, of
increasing accuracy, indexed by integers $d,n\ge0$, which can be
expressed as semidefinite programs.

\begin{thm}\label{thm:nonlinear_SOS}
The origin is globally asymptotically stable if there exists a $V\in
\R[x]$, of degree $d$, and an $\alpha >0$, such that the following
conditions hold.
\begin{align*}
V(x)-\alpha x^T x &\in \Sigma_s\\
-(\nabla V(x)^T f(x) + \alpha x^T x) &\in \Sigma_s
\end{align*}
\end{thm}

%For any particular values of $d,n$, the conditions in
%Theorem~\ref{thm:nonlinear_SOS} constitute linear constraints on the
%coefficients of $v$ and $s$. The stability condition can therefore
%be expressed as a semidefinite program of polynomial-time
%complexity. Furthermore, as $d,n \rightarrow \infty$, the condition
%increases in accuracy.

\subsection{Local Stability and Parametric
Uncertainty}\label{sec:nonlinear_ODE_PD}
Consider a time-invariant ordinary differential equation of the
following form.
\[\dot{x}(t)=f(p,x(t))
\]
Here $x(t)\in \R^n$, $p\in \R^{d_p}$ and $f$ is a polynomial of
degree $m$ such that $f(p,0)=0$ for all $p \in K \subset \R^{d_p}$.
We first give a local, parameter-dependent version of
Theorem~\ref{thm:Lyapunov_global}.

\begin{thm}\label{thm:Lyapunov_local}
Let $\Omega \subset \R^n$ contain an open neighborhood of the origin
and $f\in \R[p,x]$ satisfy $f(p,0)=0$ for all $p\in K$. Suppose
there exists an $\alpha>0$ and a Lyapunov function, $V(p,x)$, with
continuous derivatives such that $V(p,0)=0$ for $p\in K$ and the
following holds.
\begin{align*}
V(p,x)&\ge \alpha x^T x \qquad &\text{for all }x \in \Omega, \; p \in K\\
\nabla_x V(p,x)^T f(p,x) &\le 0 \qquad &\text{for all }x \in \Omega,
\; p \in K
\end{align*}
Then, for any $p\in K$, the system defined by $f$ is stable on some
open neighborhood of the origin. If, in addition, $V(p,x)^T f(p,x)
\le -\alpha x^T x$ for all $x \in \Omega$ and $p \in K$, then for
any $p \in K$, the system is asymptotically stable on some open
neighborhood of the origin.
\end{thm}

Suppose the sets $\Omega$ and $K$ are given by $K_x$ and $K_p$,
respectively, where $K_x:=\{x\in \R^n : g_i(x) \ge 0 , i=1, \ldots,
n_g\}$ for $g_i \in \R[x]$ and $K_p:=\{p \in \R^{d_p}: h_i(x) \ge 0
, i=1, \ldots, n_h\}$ for $h_i\in \R[x]$. Now let the set
$K_{\times}$ be the intersection of the sets $K_x$ and $K_p$ in the
product space $\R^{n+d_p}$, defined as follows.
\[K_{\times}:=\{(p,x)\in\R^{n+d_p}: p \in K_p,\;x\in K_x\}
\]
%Then the stability conditions associated with the local
%parameter-dependent Lyapunov theorem can be expressed as a convex
%program in the convex cone $\mathcal{P}^{K_{\times}}$.

By using the methods of the Sections~\ref{sec:SOSDecomp}
and~\ref{sec:PSresults}, we can pose the stability conditions
associated with Theorem~\ref{thm:Lyapunov_local} as a convex
feasibility problem in the convex cone $\mathcal{P}^{K_{\times}}$.
Then, for a given $d \in Z^+$, by approximating
$\mathcal{P}^{K_{\times}}$ by $\Lambda_d^{K_{\times}}$ or
$\Upsilon_d^{K_{\times}}$, we can construct a sequence of sufficient
conditions, of increasing accuracy, indexed by integers $d,r\ge0$,
which can be expressed as semidefinite programs.

%\begin{defn}
%For polynomials $\{g_i\}_{i=1}^{n_g}$ and $\{p_i\}_{i=1}^{n_p}$,
%define semialgebraic sets $K_g:=\{x: g_i(x) \ge 0 , i=1 \ldots
%n_g\}$ and $K_p:=\{x: p_i(x) \ge 0 , i=1 \ldots n_p\}$.
%\end{defn}

\begin{thm}\label{thm:nonlinearSOS_local} Let $f$ be a polynomial
such that $f(p,0)=0$ for all $p \in K_p$ and $K_x$ contain an open
neighborhood of the origin. Suppose there exists a $V \in \R[p,x]$,
of degree $r$, and an $\alpha
>0$ such that the following conditions hold.
\begin{align*}
V(p,x)-\alpha x^T x & \in \Upsilon_d^{K_{\times}} \\
-(\nabla_x V(p,x)^T f(p,x) + \alpha x^T x) & \in
\Upsilon_d^{K_{\times}}
\end{align*}
%Furthermore, we can take $X$ to be as follows.
%\begin{align*}
%X&=\{x:v(x,y)\le \gamma \text{ for all }y
%\in G \}\\
%\gamma&:=\min_{x \not \in P, y \in G} v(x,y)
%\end{align*}
Then, for any $p\in K_p$, the origin is asymptotically stable in
some open region about the origin.
\end{thm}

%\begin{thm}\label{thm:nonlinearSOS_local} Suppose that $f$ is a polynomial
%such that $f(p,0)=0$ for all $p \in K_p$. The origin is
%asymptotically stable in some region, $X$, for all $y \in G$ if
%there exists a polynomial, $v(x,y)$, of degree $d$, an $\alpha
%>0$, and sum-of-squares polynomials $\{s_i\}_{i=0}^{n_p+n_g}$,
%$\{\hat{s}_i\}_{i=0}^{n_p+n_g}$, of degree $d$, such that the
%following conditions hold.
%\begin{align*}
%v(x,y)-\alpha x^T x &= s_0(x,y)+\sum_{i=1}^{n_p} p_i(x)s_i(x,y)+\sum_{i=n_p+1}^{n_p+n_g} g_i(y)s_i(x,y) \\
%\nabla_x v(x,y)^T f(x,y) + \alpha x^T x&=
%\hat{s}_0(x,y)+\sum_{i=1}^{n_p}
%p_i(x)\hat{s}_i(x,y)+\sum_{i=n_p+1}^{n_p+n_g} g_i(y)\hat{s}_i(x,y)
%\end{align*}
%Furthermore, we can take $X$ to be as follows.
%\begin{align*}
%X&=\{x:v(x,y)\le \gamma \text{ for all }y
%\in G \}\\
%\gamma&:=\min_{x \not \in P, y \in G} v(x,y)
%\end{align*}
%\end{thm}
%
%\begin{proof}
%If the conditions of the theorem hold, then $v(x,y)>\alpha x^T x$
%for all $y \in G$, $x \in P$. In addition, $\nabla_x v(x,y)^T f(x,y)
%\le -\alpha x^T x$ for all $y \in G$, $x \in P$. Furthermore,
%$X\subset V_\gamma(y)$ for any $y \in G$, where $V_\gamma(y)$ is the
%$\gamma$-level set of $v(x,y)$ at value $y$. Finally, $V_\gamma(y)
%\subset P$ for any $y \in G$.
%\end{proof}

\section{The SOS Representation of Matrix
Functions}\label{subsec:matrixfns}

The previous sections discussed tractable tests for membership in
the convex cones $\mathcal{P}^+$ and $\mathcal{P}^K$ of
non-negativity scalar polynomials. In this section, we present an
extension of these results to the convex cones $\mathcal{S}_n^+$ and
$\mathcal{S}_n^K$ of non-negative matrices of polynomials. A simple
test for $M \in \mathcal{S}_n^+$ or $M \in \mathcal{S}_n^K$ can be
performed through the introduction of $n$ auxiliary variables.

\begin{lem}\label{lem:SOS_mat_test_SOS}
For $M \in \S^n[x]$, let $f_M \in \R[x,y]$ be defined as
$f(x,y):=y^T M(x)y$.
\begin{enumerate}
\item $M \in \mathcal{S}_n^+$ if and only if $f_M \in
\mathcal{P}^+$.
\item $M \in \mathcal{S}_n^K$ if and only if $f_M \in
\mathcal{P}^{K'}$, where $K':=\{(x,y):x \in K\}$.
\end{enumerate}
\end{lem}

Thus, we can determine whether $M \in \mathcal{S}_n^+$ by testing
whether $f_M \in \mathcal{P}^+$ using the techniques discussed in
the previous section. The problem with this approach is that the
size of the resulting semidefinite program scales very poorly in the
dimension of $M$. This is due to the introduction of the auxiliary
variables $y$. Instead, we will consider the following direct
approximation to the cone $\mathcal{S}_n^+$.

\begin{defn}
$M \in \S^n[x]$ satisfies $M\in \bar{\Sigma}_s$ if it can be
represented in the following form for some finite set of polynomial
matrices $G_i\in \R^{n \times n}[x]$, $i=1, \ldots, m$.
\[M(x)=\sum_{i=1}^m G_i(x)^T G_i(x)
\]
An element $M \in \bar{\Sigma}_s$ is referred to as a
\eemph{sum-of-squares matrix function}. Furthermore, define
$\bar{\Sigma}_s^d$ to be elements of $\bar{\Sigma}_s$ of degree $d$
or less.
\end{defn}

We now provide a tractable test for membership in $\bar{\Sigma}_s$.

\begin{lem}\label{lem:SOS_mat_test}Suppose $M \in \S_{2d}^n[x]$ for $x \in \R^{n_2}$.
$M \in \bar{\Sigma}_s$ if and only if there exists some matrix $Q
\in \S^{n_z \cdot n}$, where $n_z= {n_2+d \choose d}$ such that $Q
\ge 0$ and the following holds.
\[M(x)=\left(\bar{Z}^n_d[x]\right)^T Q \bar{Z}^n_d[x]
\]
\end{lem}
See Appendix~\ref{app:SOStheory} for Proof.

The complexity of the membership test associated with
Lemma~\ref{lem:SOS_mat_test} is considerably lower than that
associated with Lemma~\ref{lem:SOS_mat_test_SOS}. Specifically, the
variables associated with the first test are of order $n{n_2+d
\choose d}$ as opposed to order ${n+n_2+d+1 \choose d+1}$.
Furthermore, the use of the test associated with
Lemma~\ref{lem:SOS_mat_test} does not increase the conservatively,
as evidenced by the following lemma.

\begin{lem} $M(x) \in \bar{\Sigma}_s$ if and only
if $y^T M(x) y \in \Sigma_s$.
\end{lem}
See Appendix~\ref{app:SOStheory} for Proof.

%Thus for a given vector of monomials, $Z(x)$, of length $q$, if
%there exists a positive semidefinite matrix $P \in \S^{nq}$ such
%that
%\[M(x)=\bmat{Z(x)&&\\&\ddots&\\&&Z(x)}^T P \bmat{Z(x)&&\\&\ddots&\\&&Z(x)}
%\]
%Then $M(x) \ge 0 \text{ for all }x \in \R^p$. We note for a $M$ of
%degree $d$, that the semidefinite program associated with this
%approach uses variables in $\S^{n(d/2)^p}$, a considerable
%improvement over the standard SOS approach which would use variables
%in $\S^{(d/2+1)^{n+p}}$.

As for $\Sigma_s$ and $\mathcal{P}^+$, the question of how well
$\bar{\Sigma}_s$ approximates $\mathcal{S}_n^+$ is an open question.
However, in the case of a single variable, we have the following
result due to Choi et al.~\cite{choi_1980}.

\begin{lem}\label{lem:SOS_mat_1d} For $x \in \R$, $M \in \S^n[x]$, the following are equivalent
\begin{itemize}\item $M \in \bar{\Sigma}_s$ \item $M \in \mathcal{S}_n^+$.
\end{itemize}
\end{lem}

Lemma~\ref{lem:SOS_mat_1d} will find considerable application in
Chapter~\ref{chp:LinearCase}.

\section{Positivstellensatz Results for Matrix Functions}

In this section, we consider membership in the set
$\mathcal{S}^{K_f}_n$ for some semi-algebraic set $K_f$, defined by
polynomials $f_i \in \R[x]$ as follows.
\[K_f:=\{x : f_i(x)\ge 0 , \; i=1\ldots n_f\}
\]

Scherer and Hol~\cite{scherer_2005} have proven the following
generalization of Putinar's Positivstellensatz.

\begin{thm}~\label{thm:scherer} Suppose that $K_f$ is P-compact and $M \in
\mathcal{S}^{K_f}_n$. Then there exist $\eta>0$ and $S_i \in
\bar{\Sigma}_s$ for $i=1,\ldots, n_f$ such that the following holds.
\[M(x)-\sum_{i=1}^{n_f}f_i(x)S_i(x) -\epsilon I \in \bar{\Sigma}_s
\]
\end{thm}

Theorem~\ref{thm:scherer} does not represent a tractable test for
$\mathcal{S}^{K_f}_n$, since no degree bounds are given on the
$S_i$. However, the result does suggest a sequence of sufficient
conditions, indexed by $d \in \Z^+$, and defined as follows.

\begin{defn}
We say that $M \in \bar{\Upsilon}^{K_f}_{n,d}$ if there exist
$\eta>0$ and $S_i \in \bar{\Sigma}^d_s$ for $i=1,\ldots, n_f$ such
that the following holds.
\[M(x)-\sum_{i=1}^{n_f}f_i(x)S_i(x) \in \bar{\Sigma}_s
\]
\end{defn}

\begin{lem}
For any $d\ge0$ and $K_f$, $\bar{\Upsilon}^{K_f}_{n,d} \subset
\mathcal{S}^{K_f}_n$.
\end{lem}

Therefore, for any $d \ge 0$, we can construct a tractable test for
membership in the convex cone $\bar{\Upsilon}^{K_f}_{n,d}$,
expressible as a semidefinite program. Furthermore, as $d$
increases, $\bar{\Upsilon}^{K_f}_{n,d}$ will constitute a better
approximation for $\mathcal{S}^{K_f}_n$.

\subsection{Linear ODEs with Parametric
Uncertainty}\label{sec:Linear_ODE}
The most common application for the sum-of-squares matrix
representation is the case of a linear finite-dimensional system
which contains uncertainty. Consider a system of ordinary
differential equations of the following form where $A:\R^{n_y}
\rightarrow \R^{n \times n}$.
\begin{equation}\dot{x}(t)=A(y)x(t) \qquad y \in G
\label{eqn:LTI_uncertain}
\end{equation}
The following is a standard result in analysis of LTI systems.

\begin{thm}~\label{thm:linear_lyapunov} The system defined by Equation~\ref{eqn:LTI_uncertain}
is globally asymptotically stable for all $y \in G$ if and only if
there exists a matrix function $P:G \rightarrow \S^n$ such that the
following conditions hold.
\begin{align*}
P(y)-I&\in \mathcal{S}^G_n\\
-(A(y)^T P(y)+P(y)A(y)+I) &\in \mathcal{S}^G_n\\
\end{align*}
\end{thm}

Now suppose that $G$ is a semi-algebraic set defined by polynomials
$\{g_i\}_{i=1}^{n_g}$. Then for integer $d \ge 0$, we can
approximate the conditions~\ref{thm:linear_lyapunov} using the
following.

\begin{thm}\label{thm:linear_lyapunov_test} The system defined by Equation~\ref{eqn:LTI_uncertain}
is globally asymptotically stable for all $y \in G$ if there exists
a matrix of polynomials $P\in \S^n[y]$ such that the following
conditions hold.
\begin{align*}
P(y)-I&\in \Upsilon^G_{n,d}\\
-(A(y)^T P(y)+P(y)A(y)+I) &\in \Upsilon^G_{n,d}\\
\end{align*}
\end{thm}

The conditions of Theorem~\ref{thm:linear_lyapunov_test} can be
expressed as a semidefinite program.

\section{Conclusion}

This chapter has detailed a method for expressing may difficult
problems in control as convex optimization problems. These
optimization problems are, in general, not tractable due to the lack
of an efficient membership test for the convex cone $\mathcal{P}^+$.
However, we have shown that by approximating $\mathcal{P}^+$ by
$\Sigma_s$, one can construct tractable approximations. The
techniques and concepts of this chapter will play an integral role
as we seek to generalize the ODE applications presented in
subsections~\ref{sec:nonlinear_ODE},~\ref{sec:nonlinear_ODE_PD},
and~\ref{sec:Linear_ODE} to differential equations which contain
delay.

\chapter{Stability of Linear Time-Delay
Systems}\label{chp:LinearCase}
\section{Introduction}

The study of stability of systems of differential equations which
contain delays has been an active area of research for some time. A
complete summary of the results in this field is beyond the scope of
this section. However, an overview of these results can be obtained
from various survey papers and books on the subject, see for
example~\cite{gu_2003,hale_1993,kolmanovskii_1999,niculescu_2001}.
The previous results on this subject can be grouped into analysis
either in the frequency-domain or in the time-domain.
Frequency-domain techniques apply only to linear systems and
typically attempt to determine whether all roots of the
characteristic equation of the system lie in the left half-plane.
Time-domain techniques generally use Lyapunov-based analysis, an
approach which was extended to infinite dimensional systems by
Krasovskii in~\cite{krasovskii_1963}. Stability results in this area
are grouped into delay-dependent and delay-independent conditions.
If a delay-dependent condition holds, then stability is guaranteed
for a specific value or range of values of the delay. If a
delay-independent condition holds, then the system is stable for all
possible values of the delay. A particularly interesting result on
the linear delay-independent case appears in~\cite{bliman_2002b}.

In this chapter, we show how to compute solutions to an
operator-theoretic version of the Lyapunov inequality using
semidefinite programming. This inequality, defined by the derivative
of the complete quadratic functional, can be posed as a convex
feasibility problem over certain infinite-dimensional convex cones
defined by positive operators on certain subspaces of
$\mathcal{C}_\tau$. Our result is expressed as a nested sequence of
sufficient conditions which are of increasing accuracy and can be
tested using semidefinite programming. Our approach is based on
parameterizing the convex cones mentioned above using polynomial
functions of bounded degree.

\section{Positive Operators}

In this section, we present the convex cones of positive multiplier
and integral operators which define the complete quadratic
functional and some its derivative forms. We then use polynomials of
bounded degree to parametrization certain convex subsets of these
cones in terms of positive semidefinite matrices. These results will
allow us to express the Lyapunov stability conditions in terms of
semidefinite programming problems. To begin, consider the complete
quadratic functional.
\begin{align*}&V(x) = x(0)^T P x(0)+2 \int_{-\tau}^0
x(0)^T Q(\theta) x(\theta) d\theta\\ &+\int_{-\tau}^0 x(\theta)^T
S(\theta) x(\theta) d\theta+\int_{-\tau}^0\int_{-\tau}^0 x(\theta)^T
R(\theta,\omega)x(\omega)d \theta d \omega
\end{align*}
We can associate with $V$ and $\epsilon>0$ an operator $A:
\mathcal{C}_\tau \rightarrow \mathcal{C}_\tau$ such that
$V(x)=\ip{x}{Ax}$, where $A$ is defined as follows.
\begin{align*}(Ay)(\theta)&=\bmat{P-\epsilon I & Q(\theta)\\Q(\theta)^T
&S(\theta)}y(\theta)+\int_{-\tau}^0
\bmat{0&0\\0&R(\theta,\omega)}y(\omega)d \omega\\
&=(A_1y)(\theta)+(A_2y)(\theta), \notag
\end{align*}
The operator $A$ is a combination of multiplier operator $A_1$ and
integral operator $A_2$. The complete quadratic functional, $V$
satisfies the positivity condition of the stability theorem with
$u(\phi(0))=\epsilon \phi(0)^2$ if the integral operator, $A_2$, is
positive on $\mathcal{C}_\tau$ and there exists some $\epsilon>0$
such that the multiplication operator, $A_1$, is positive on the
subspace $X \subset \mathcal{C}_\tau$ where
\[X:=\{ x\in\mathcal{C}_\tau  \ \vert\ ,
  x=\bmat{x_1 \\ x_2} \text{and }x_1(\theta)=x_2(0) \text
      { for all }\theta \}.
\]
We now define the specific sets of operators which define the
complete quadratic functional and its derivative forms.

\begin{defn} For a matrix-valued function $M:\R \rightarrow \S^{n}$,
we denote by $A_M$ the multiplication operator such that
\[(A_M x)(\theta)=M(\theta)x(\theta)
\]
\end{defn}
\begin{defn} For a matrix-valued function
$R:\R \times \R \rightarrow \R^{n \times n}$ such that
$R(\theta,\omega)=R(\omega,\theta)^T$, we denote by $B_R$ the
integral operator such that
\[(B_R x)(\theta)=\int_{-\tau}^0 R(\theta,\omega)x(\omega)d \omega
\]
\end{defn}
%In the case of a single or distributed delay, we need the sets
%$H_1^+$ and $H_2^+$, defined as follows.
\begin{defn} For a continuous matrix-valued function $M$, we say $M\in H_1^+$
if $\ip{x}{A_M x} \ge 0$ for all $x \in X$.
\end{defn}
\begin{defn} For a continuous matrix-valued function
$R:\R \times \R \rightarrow \R^{n \times n}$ such that
$R(\theta,\omega)=R(\omega,\theta)^T$, we say $R\in H_2^+$ if
$\ip{x}{B_R x} \ge 0$ for all $x \in \mathcal{C}_\tau$.
\end{defn}
%For the case of multiple discrete delays, we also need the sets
%$\tilde{H}_1^+$ and $\tilde{H}_2^+$, defined as follows.
\begin{defn} For a piecewise-continuous matrix-valued function $M$,
we say $M\in \tilde{H}_1^+$ if $\ip{x}{A_M x} \ge 0$ for all $x \in
X$.
\end{defn}
\begin{defn} For a piecewise-continuous matrix-valued function
$R:\R \times \R \rightarrow \R^{n \times n}$ such that
$R(\theta,\omega)=R(\omega,\theta)^T$, we say $R\in \tilde{H}_2^+$
if $\ip{x}{B_R x} \ge 0$ for all $x \in \mathcal{C}_\tau$.
\end{defn}
%To parameterize the derivatives of the complete quadratic
%functional, we need the additional sets $H_3^+$ and $\tilde{H}_3^+$,
%defined as follows.
\begin{defn} For a continuous matrix-valued function $M$, we say $M\in H_3^+$
if $\ip{x}{A_M x} \ge 0$ for all $x \in X_3$, where
\begin{align*}X_3:=\{& x\in\mathcal{C}_\tau  \ \vert\ ,
  x=\bmat{x_1^T & x_2^T & x_3^T}^T \text{ and }x_1(\theta)=x_3(0) \text{ and } \\
  &x_2(\theta)=x_3(-\tau) \text
      { for all }\theta \}.
\end{align*}
\end{defn}
\begin{defn} For a piecewise-continuous matrix-valued function $M$, we say $M\in \tilde{H}_3^+$
if $\ip{x}{A_M x} \ge 0$ for all $x \in \tilde{X}_3$, where
\begin{align*}&\tilde{X}_3:=\bbl\{ x\in\mathcal{C}_\tau  \ \vert\ ,
  x=\bmat{x_1^T & \ldots & x_{K+2}^T}^T \text{ and } \\
  &x_i(\theta)=x_{K+2}(-\tau_{i-1}) \text
      { for all } \theta,\; i=1 ,\ldots,
K+1 \bbr\}.
\end{align*}
\end{defn}
%The convex cones $H_i^+$ and $\tilde{H}_i^+$ are
%infinite-dimensional. We will use the results of the following
%subsections to parameterize certain subsets of these cones using the
%space of positive semidefinite matrices.

%%%%%%%%%%%%%%%%%%%%%%%%%%%%%%%%%%%%%%%%%%%%%%%%%%%%%%%%%%%%5
\subsection{Multiplier Operators and Spacing Functions}

In this subsection, we consider the convex cones $H_1^+$ and $H_3^+$
which define positive multiplication operators on $X_1$ and $X_3$
respectively. We use the following theorem.

\begin{thm} Suppose $M:\R \rightarrow \S^{2n}$ is a continuous matrix-valued function. Then the following are equivalent
\begin{enumerate}
\item There exists some $\epsilon>0$ such that the following holds for all
$x \in \mathcal{C}_\tau$.
\[\int_{-\tau}^0 \bmat{x(0)\\x(\theta)}^T M(\theta)
\bmat{x(0)\\x(\theta)} d \theta \ge \epsilon \norm{x}^2_2
\]
\item There exists some $\epsilon'>0$ and some continuous matrix-valued function
$T:\R \rightarrow \S^n$ such that the following holds.
\begin{align*}
&\int_{-\tau}^0 T(\theta) d \theta = 0\\
&M(\theta)+\bmat{T(\theta) & 0 \\ 0 & -\epsilon' I} \ge 0 \qquad
\text{for all } \theta \in [-\tau,0]
\end{align*}
\end{enumerate}\label{lem:H1_LMI}
\end{thm}
See Appendix~\ref{app:LinearCase} for Proof.

\begin{defn}
Given $n$ and $\tau$, we refer to any piecewise-continuous
matrix-valued function $T:\R \rightarrow \R^n$ such that
$\int_{-\tau}^0 T(\theta) d\theta=0$ as a \eemph{spacing function},
denoted $T \in \Omega$.
\end{defn}

This theorem states that any element of $H_1^+$ which satisfies a
certain additional positivity condition can be represented by the
combination of a positive semidefinite matrix-valued function and a
spacing function. This allows us to search for elements of $H_1^+$
by simultaneously searching over the set of positive semidefinite
matrix-valued functions and the set of spacing functions. For any $d
\in \Z^+$, the following gives a tractable condition for membership
in $H_1^+$.

\begin{defn}
Given $n$ and $\tau$, define $\Omega_p^d$ to be the set of functions
$T \in \S^n_d[\theta]$ such that the following holds.
\[\int_{-\tau}^0 T(\theta) d \theta = 0
\]
\end{defn}

\begin{defn}For a given $d\ge0$, $M \in \S^{2n}[\theta]$ satisfies $M \in
G^{d}_{1}$ if there exists a $T\in\Omega_p^d$ such that
\[M(\theta)+\bmat{T(\theta) & 0 \\ 0 & 0} \in
\mathcal{S}_{2n}^+.
\]
\end{defn}

Membership in the convex cone $G^{d}_{1}$ can be implemented as a
semidefinite programming constraint by noting that the constraint
that a polynomial integrate to $0$ is affine in the coefficients and
that for the single variable case, $\bar{\Sigma}_s=\mathcal{S}_n^+$.
We now consider $H_3^+$. We can use a simple extension of the
previous theorem.
\begin{lem} Let $M:\R \rightarrow \S^{3n}$ be a continuous
matrix-valued function. Then the following are equivalent.
\begin{enumerate}
\item There exists an $\epsilon>0$ such that the following holds for
all $x \in \mathcal{C}_\tau$.
\[\int_{-\tau}^0\bmat{x(0)\\x(-\tau)\\x(\theta)}^T M(\theta)\bmat{x(0)\\x(-\tau)\\x(\theta)}\ge \epsilon \norm{x}_2^2
\]
\item There exists some $\epsilon'>0$ and a continuous matrix-valued
function $T:\R \rightarrow \S^{2n}$ such that the following holds.
\begin{align*}
&\int_{-\tau}^0 T(\theta) d \theta = 0\\
&M(\theta)+\bmat{T(\theta) & 0 \\ 0 & -\epsilon'I} \ge 0 \qquad
\text{for all } \theta \in [-\tau,0]
\end{align*}
\end{enumerate}\label{lem:H3_LMI}
\end{lem}
See Appendix~\ref{app:LinearCase} for Proof.

\begin{defn}For a given $d\ge0$, $M \in \S^{3n}[\theta]$ satisfies $M \in
G^{d}_{3}$ if there exists a $T\in \S^{2n}[\theta]$ such that
$T\in\Omega_p^d$ and
\[M(\theta)+\bmat{T(\theta) & 0 \\ 0 & 0} \in
\mathcal{S}_{3n}^+.
\]
\end{defn}

Similar to Theorem~\ref{lem:H1_LMI}, Lemma~\ref{lem:H3_LMI} allows
us to use $G^{d}_{3} \subset H_3^+$ to replace $M \in H_3^+$ with a
semidefinite programming constraint.

\subsubsection{Piecewise-Continuous Spacing Functions}

We now consider the sets $\tilde{H}_1^+$ and $\tilde{H}_3^+$ of
piecewise continuous functions which define positive multiplication
operators on $X_1$ and $\tilde{X}_3$ respectively. The following
lemma is a generalization of Theorem~\ref{lem:H1_LMI}.
\begin{lem}
Suppose $S_i:\R \rightarrow \S^{2n}$, $i=1 ,\ldots, K$ are
continuous symmetric matrix-valued functions with domains
$[-\tau_i,-\tau_{i-1}]$ where $\tau_i > \tau_{i-1}$ for $i=1,\ldots,
K$ and $\tau_0=0$. Then the following are equivalent.
\begin{enumerate}
\item There exists an $\epsilon>0$ such that the following holds
for all $x\in \mathcal{C}_{\tau_K}$.
\[
\sum_{i=1}^K \int_{-\tau_i}^{-\tau_{i-1}} \bmat{x(0)\\x(\theta)}^T
S_i(\theta)
\bmat{x(0)\\x(\theta)}d \theta \ge \epsilon \norm{x}_2^2\\
\]
\item There exists an $\epsilon'>0$ and continuous symmetric matrix
valued functions, $T_i:\R \rightarrow \S^n$, such that
\begin{align*}
&S_i(\theta)+\bmat{T_i(\theta) & 0\\0&-\epsilon'I} \ge 0 \qquad
\text{for }\theta \in [-\tau_i,-\tau_{i-1}],\; i=1 ,\ldots,
K \; \text{and} \\
&\sum_{i=1}^K \int_{-\tau_i}^{-\tau_{i-1}}T_i(\theta)=0.
\end{align*}
\end{enumerate}
\label{lem:H1_pc_LMI}
\end{lem}
See Appendix~\ref{app:LinearCase} for Proof.

Lemma~\ref{lem:H1_pc_LMI} allows us to represent the constrain $S
\in \tilde{H}_1^+$ by considering a piecewise-continuous function
$S$ to be defined by $S(\theta):=S_i(\theta)$ for $\theta \in
[-\tau_{i},-\tau_{i-1}]$.

\begin{defn}Given $n$ and $\tau$, define $\tilde{\Omega}^d_p$ to be the set of
piecewise-continuous functions $T:\R \rightarrow \S^n$ such that
$T(\theta)=T_i(\theta)$ for $T_i \in \S^n_d[\theta]$ and
$\theta\in[-\tau_{i},-\tau_{i-1}]$ where $i=1,\ldots,K$ and such
that
\[\sum_{i=1}^K \int_{-\tau_i}^{-\tau_{i-1}}T_i(\theta)=0.
\]
\end{defn}
%We can use $\bar{\Sigma}_s$ and $\tilde{\Omega}^d_p$ to construct
%elements of $\tilde{H}_1^+$.

\begin{defn}
For a given $d\ge0$, we say that $M : \R \rightarrow \S^{2n}$
satisfies $M \in \tilde{G}^{d}_{1}$ if there exists $M_i\in
\S_d^{2n}[\theta]$  such that $M(\theta)=M_i(\theta)$ for $\theta
\in [-\tau_i,-\tau_{i-1}]$ and there exists a
$T\in\tilde{\Omega}_p^d$, defined by $T_i\in \S_d^{n}[\theta]$ on
the interval $[-\tau_i,-\tau_{i-1}]$, and such that
\[M_i(\theta)+\bmat{T_i(\theta) & 0\\0&0} \in
\mathcal{S}_{2n}^+ \qquad i=1,\ldots, K.
\]
\end{defn}
Similar to Theorem~\ref{lem:H1_LMI}, Lemma~\ref{lem:H1_pc_LMI}
allows us to use $\tilde{G}^{d}_{1} \subset \tilde{H}_1^+$ to
replace $M \in \tilde{H}_1^+$ with a semidefinite programming
constraint. Now for $\tilde{H}_3^+$, we have the following lemma.

\begin{lem}
Suppose $S_i:\R \rightarrow \S^{n(K+2)}$ are continuous matrix
valued functions with domains $[-\tau_i,-\tau_{i-1}]$ for
$i=1,\ldots, K$ where $\tau_i>\tau_{i-1}$ for $i=1,\ldots, K$ and
$\tau_0=0$. Then the following are equivalent.
\begin{enumerate}
\item There exists an $\epsilon>0$ such that the following holds
for all $x\in \mathcal{C}_{\tau_K}$.
\[
\sum_{i=1}^K \int_{-\tau_i}^{-\tau_{i-1}}
\bmat{x(-\tau_0)\\\vdots\\x(-\tau_K)\\x(\theta)}^T S_i(\theta)
\bmat{x(-\tau_0)\\\vdots\\x(-\tau_K)\\x(\theta)}d \theta \ge \epsilon \norm{x}_2^2\\
\]
\item There exists an $\epsilon'>0$ and continuous matrix
valued functions $T_i:\R \rightarrow \S^{n(K+1)}$ such that
\begin{align*}
&S_i(\theta)+\bmat{T_i(\theta) & 0\\0&-\epsilon'I} \ge 0 \qquad  \text{for } \theta \in [-\tau_i,-\tau_{i-1}], \; i=1 ,\ldots, K \text{ and }\\
&\sum_{i=1}^K \int_{-\tau_i}^{-\tau_{i-1}}T_i(\theta)=0.
\end{align*}
\end{enumerate}
\label{lem:H3_pc_LMI}
\end{lem}
See Appendix~\ref{app:LinearCase} for Proof.

%Similar to Lemma~\ref{lem:H1_pc_LMI}, Lemma~\ref{lem:H3_pc_LMI}
%allows us to construct elements of $\tilde{H}_3^+$ using elements of
%$\bar{\Sigma}_s$ and $\tilde{\Omega}^d_p$.

\begin{defn}
For a given $d\ge0$, we say that $M:\R \rightarrow \S^{n(K+2)}$
satisfies $M \in \tilde{G}^{d}_{3}$ if there exist $M_i\in
\S_d^{n(K+2)}[\theta]$ such that $M(\theta)=M_i(\theta)$ on
$[-\tau_i,-\tau_{i-1}]$ and $T_i\in \S_d^{n(K+1)}[\theta]$ such that
$T(\theta)=T_i(\theta)$ on $[-\tau_i,-\tau_{i-1}]$ implies
$T\in\tilde{\Omega}_p^d$ and such that
\[M_i(\theta)+\bmat{T_i(\theta) & 0\\0&0} \in
\mathcal{S}_{n(K+2)}^+ \qquad i=1,\ldots, K.
\]
\end{defn}

Similar to Theorem~\ref{lem:H1_LMI}, Lemma~\ref{lem:H3_pc_LMI}
allows us to use $\tilde{G}^{d}_{3} \subset \tilde{H}_3^+$ to
replace $M \in \tilde{H}_3^+$ with a semidefinite programming
constraint.

\subsubsection{Parameter-Dependent Spacing Functions}
In this subsection, we consider the case which arises when the
dynamics of the system depend on some uncertain time-invariant
vector of parameters.
\begin{lem}\label{lem:H1_LMI_pd}
For a given $d \in \Z^+$ and $S \in \S^{2n}_d[\theta,y]$, suppose
there exists some $T\in \S^n_d[\theta,y]$ such that the following
holds.
\begin{align*}&\int_{-\tau}^0
T(\theta,y) d
\theta=0 \\
&S(\theta,y)+\bmat{T(\theta,y) & 0\\
0&0} \in \bar{\Sigma}_s
\end{align*}
Then $S(\cdot,y) \in H_1^+$ for all $y\in \R^p$.
\end{lem}
Now suppose we wish to impose the condition that $S(\cdot,y) \in
H_1^+$ for all $y \in \{y: p_i(y) \ge 0, i=1 ,\ldots, N\}$ for
scalar polynomials $p_i$.
\begin{lem}
Suppose there exist functions $S_i$ such that $S_i(\cdot,y) \in
H_1^+$ for $i=0,\ldots, N$ for all $y$ and such that the following
holds.
\[S(\theta,y)=S_0(\theta,y)+\sum_{i=1}^N p_i(y) S_i(\theta,y),
\]
Then $S(\cdot,y) \in H_1^+$ for all $y \in \{y: p_i(y) \ge 0, i=1
,\ldots, N\}$\label{lem:H1_LMI_pd2}
\end{lem}
The purpose of these lemmas is to show how our representation of
$H_1^+$ can be used to prove stability in the presence of parametric
uncertainty. Clearly, these results are motivated by the
construction $\bar{\Upsilon}^K_d$ presented in
Chapter~\ref{chp:SOStheory}. The generalization of
Lemmas~\ref{lem:H1_LMI_pd} and~\ref{lem:H1_LMI_pd2} to $H_3^+$,
$\tilde{H}_1^+$ and $\tilde{H}_3^+$ is similar and therefore not
explicitly discussed.

\subsection{Positive Finite-Rank Integral Operators}

We now consider the set $H_2^+$ of continuous matrix-valued kernel
functions which define positive integral operators on
$\mathcal{C}_\tau$. These operators are compact and satisfy the
following inequality for all $x \in \mathcal{C}_\tau$.
\[\int_{-\tau}^0\int_{-\tau}^0  x(\theta)^T k(\theta,\omega) x(\omega)
d \theta d \omega \ge 0
\]

Because $\mathcal{C}_\tau$ is an infinite dimensional space, we
cannot fully parameterize the operators which are positive on this
space using a finite-dimensional set of positive semidefinite
matrices. Instead, we will consider the subset of finite-rank
operators.

\begin{thm} Let $A$ be a compact Hermitian operator which
is positive on $Y_d \subset \mathcal{C}_\tau$, where
\begin{align*}Y_d:=\{&p \in \mathcal{C}_\tau: p \text{ is an $n$-dimensional vector of}\\
& \text{univariate polynomials of degree $d$ or less.}\}
\end{align*}
Then there exists a positive semidefinite matrix $Q\in \S^{n \times
(d+1)}$ such that
\begin{align*}
\ip{x}{Ax}&= \int_{-\tau}^0 \int_{-\tau}^0 x(\theta)^T k(\theta,\omega) x(\omega)d \theta d \omega \qquad \forall\, x \in Y_m\\
k(\theta,\omega)&:=\bmat{Z_d[\theta]&&\\&\ddots&\\&&Z_d[\theta]}^T Q
\bmat{Z_d[\omega]&&\\&\ddots&\\&&Z_d[\omega]}\\
&=\bar{Z}^n_d(\theta)^T Q \bar{Z}_d^n(\omega),
\end{align*}
where recall $Z_d[\theta]$ is the $(d+1)$-dimensional vector of
powers of $\theta$ of degree $d$ or less.
\end{thm}

Before proving the theorem, we quote the following lemma which
follows directly from the Spectral Theorem~\cite{young_1988}.

\begin{lem} Let $\{e_i\}_{i=1}^p$ be a basis for some finite dimensional
subspace $Y$ of an inner product space $Z$. Let $A$ be some compact
Hermitian operator which is positive on $Y$. Then there exists some
$K\ge0$ such that the following holds for all $x \in Y$.
\[\ip{x}{Ax}=\sum_{i,j=1}^p
K_{ij}\ip{e_i}{x}\ip{e_j}{x}
\]\label{lem:H2_LMI}
\end{lem}

\begin{proof}
Define $e_i$ to be the transpose of the $i^{\text{th}}$ row of
$\bar{Z}^n_d$, then $\{e_i\}_{i=1}^{n (d+1)}$ forms a basis for the
finite dimensional subspace $Y_m \subset \mathcal{C}_\tau$. By
Lemma~\ref{lem:H2_LMI}, there exists some $Q\ge0$ such that the
following holds where $p=n(d+1)$.
\begin{align*}\ip{x}{Ax}&=\sum_{i,j=1}^p
Q_{ij}\ip{e_i}{x}\ip{e_j}{x}\\
&=\int_{-\tau}^0 \int_{-\tau}^0 \sum_{i,j=1}^p Q_{ij}e_i(\theta)^T
x(\theta)
e_j(\omega)^T x(\omega) d \theta d \omega\\
&=\int_{-\tau}^0 \int_{-\tau}^0 x(\theta)^T \sum_{i,j=1}^p
\left(e_i(\theta) Q_{ij}
e_j(\omega)^T\right) x(\omega) d \theta d \omega\\
&=\int_{-\tau}^0 \int_{-\tau}^0 x(\theta)^T \bar{Z}^n_d[\theta]^T Q
\bar{Z}^n_d[\omega] x(\omega) d \theta d \omega
\end{align*}
\end{proof}
\begin{defn}
For a given integer $d\ge0$, we denote the set of kernel functions
which defined compact Hermitian operators which are positive on
$Y_d$ by $G_2^{2d}$.
\end{defn}
This theorem shows how we can use semidefinite programming to
represent the subset $G_2^{2d}$ which consists of finite rank
elements of the convex cone $H_2^+$ with polynomial eigenvectors of
bounded degree. Naturally, as we increase the degree bound, $d$, the
rank of the operators in $G_2^{2d}$ will increase, as will the
computational complexity of the problem.

\subsubsection{Piecewise-Continuous Positive Finite-Rank Operators}

In this subsection, we consider the set $\tilde{H}_2^+$ of
matrix-valued kernel functions which are discontinuous only at
points $\theta,\omega = \{-\tau_i\}_{i=1}^{K-1}$ and which define
positive integral operators on $\mathcal{C}_\tau$. Although the
results of the previous section can be applied directly, the kernel
functions defined in such a manner will necessarily be continuous
since they are constructed using a finite number of continuous
functions. In order to allow for the construction of operators
defined by piecewise-continuous matrix valued functions, we use the
following lemma.

\begin{lem} Let $M$ be a matrix valued function
$M: \R^2 \rightarrow \S^{n}$ which is discontinuous only at points
$\theta,\omega=-\tau_{i}$ for $i=1, \ldots, K-1$ where the $\tau_i$
are increasing and $\tau_0=0$. Then $M \in \tilde{H}_2^+$ if and
only if there exists some continuous matrix valued function $R: \R^2
\rightarrow \R^{nK \times nK}$ such that $R \in H_2^+$ and the
following holds where $I_i=[-\tau_i,-\tau_{i-1}]$,
$\Delta_i=\tau_i-\tau_{i-1}$.
\begin{align*}
&M(\theta,\omega)=M_{ij}(\theta,\omega) \qquad \text{for all }\;
\theta \in
I_i, \quad \omega \in I_j\\
&M_{ij}(\theta,\omega)=
R_{ij}\left(\frac{\tau_K}{\Delta_i}\theta+\tau_{i-1}\frac{\tau_K}{\Delta_i},
\frac{\tau_K}{\Delta_j}\omega+\tau_{j-1}\frac{\tau_K}{\Delta_j}\right)\\
&R(\theta,\omega)=\bmat{R_{11}(\theta,\omega) & \ldots&  R_{1K}(\theta,\omega) \\
\vdots & & \vdots\\
R_{K1}(\theta,\omega) &\ldots & R_{KK}(\theta,\omega)}
\end{align*} \label{lem:H2_pc_LMI}
\end{lem}
See Appendix~\ref{app:LinearCase} for Proof.

\begin{defn}For an integer $d\ge0$, we denote by $\tilde{G}_2^d$ the
set of piecewise continuous matrix valued functions $M: \R^2
\rightarrow \S^{n}$ such that there exist some $R: \R^2 \rightarrow
\R^{nK \times nK}$ such that $R \in G_2^d$ and the following holds
where $I_i$, $\Delta_i$ are as defined above.
\begin{align*}
&M(\theta,\omega)=M_{ij}(\theta,\omega) \qquad \text{for all }\;
\theta \in
I_i, \quad \omega \in I_j\\
&M_{ij}(\theta,\omega)=
R_{ij}\left(\frac{\tau_K}{\Delta_i}\theta+\tau_{i-1}\frac{\tau_K}{\Delta_i},
\frac{\tau_K}{\Delta_j}\omega+\tau_{j-1}\frac{\tau_K}{\Delta_j}\right)\\
&R(\theta,\omega)=\bmat{R_{11}(\theta,\omega) & \ldots&  R_{1K}(\theta,\omega) \\
\vdots & & \vdots\\
R_{K1}(\theta,\omega) &\ldots & R_{KK}(\theta,\omega)}
\end{align*}
\end{defn}

Lemma~\ref{lem:H2_pc_LMI} allows us to use elements of the set
$G_2^d$ to represent the subset of $\tilde{H}_2$ of finite-rank
operators with piecewise-continuous eigenvectors. Since the
transformation from $G_2^d$ to $\tilde{G}_2^d$ is affine, this
allows us to replace $R \in \tilde{H}_3$ with a semidefinite
programming constraint.

\subsection{Parameter Dependent Positive Finite-Rank Operators }

We now briefly discuss the extension of the results of the previous
two subsections when the dynamics contain parametric uncertainty.
\begin{lem}
For a given $d\in \Z^+$, suppose there exists a $P\in
\bar{\Sigma}^d_s$ such that
\begin{equation*}
M(\theta,\omega,\alpha)=\bar{Z}^d_n[\theta]^T P(\alpha)
\bar{Z}^d_n[\omega]
\end{equation*}
Then $M(\cdot,\cdot,\alpha)\in H_2^+$ for all $\alpha \in \R^n $
\end{lem}
Define $K_p:=\{\alpha: p_i(\alpha) \ge 0 , i=1 ,\ldots, l\}$.

\begin{lem} For a given $d \in \Z^+$, suppose there exist $P\in \bar{\Upsilon}^{K_p}_{n(d+1),d}$ such that
\begin{equation*}
M(\theta,\omega,\alpha)=\bar{Z}^d_n[\theta]^T P(\alpha)
\bar{Z}^d_n[\omega].
\end{equation*}
Then $M(\cdot,\cdot,\alpha) \in H_2^+$ for all $\alpha\in K_p$.
\end{lem}
These lemmas allow us to search for parameter dependent positive
integral operators using the set of positive semidefinite matrix
functions. The extension of these lemmas to the set $\tilde{H}_2^+$
should be clear and is not explicitly discussed.

\section{Results}

We now turn our attention to the following three classes of
time-delay systems which we will explicitly consider.
\begin{enumerate}
\item \quad $\dot{x}(t)=A_0 x(t)+A_1 x(t-\tau)$\label{num:single}
\item \quad $\dot{x}(t)=A_0 x(t)+\sum_{i=1}^n A_i
x(t-\tau_i)$\label{num:multiple}
\item \quad $\dot{x}(t)=A_0 x(t) + \int_{-\tau}^{0}A(\theta) x(t+\theta) d
\theta$ \label{num:distributed}
\end{enumerate}
Here $x(t) \in \R^n$, $A_i\in \R^{n \times n}$ and $A\in\R^{n \times
n}[\theta]$ is a matrix of polynomials. These classes contain
constant, finite aftereffect. We refer to system~\ref{num:single})
as the case of single delay, system~\ref{num:multiple}) as the case
of multiple delays and system~\ref{num:distributed}) as the case of
distributed delay. In this section, we show that in each of these
cases, the derivative of the complete quadratic functional along
trajectories of the system can be represented by a quadratic
functional defined by the integral and multiplier operators
parameterized in the previous section. Moreover, the transformation
from the coefficients of the polynomials defining the complete
quadratic functional to those defining its derivative is affine.

\subsection{The Single Delay Case}

In this section, we consider the first subclass of time-delay
systems which are given in the following form for matrices $A,B \in
\R^{n \times n}$
\begin{equation}\dot{x}(t)=Ax(t)+Bx(t-\tau)
\label{eqn:single_delay}
\end{equation}
The following theorem gives a condition for stability of the system.
\begin{thm}
For integer $d\ge0$, the solution map, $G$, defined by
Equation~\eqref{eqn:single_delay} is asymptotically stable if there
exists a constant $\epsilon>0$, matrix $P \in \S^n$, and continuous
matrix-valued functions $S: \R \rightarrow \S^n$, $Q:\R \rightarrow
\R^{n\times n}$ and $R:\R^2 \rightarrow \R^{n \times n}$ where
$R(\theta,\eta)=R(\eta,\theta)^T$ and such that
\begin{align*}
&\bmat{P -\epsilon I & \tau Q \\ \tau Q^T & \tau S}\in G_1^d \quad
&R \in G_2^d\\
&-D_1 \in G_3^d & M \in G_2^d
\end{align*}
\begin{align*}
&D_1(\theta) =\bmat{D_{11} & PB - Q(-\tau)
& \tau (A^T Q(\theta) - \dot{Q}(\theta)+ R(0,\theta))\\
*^T & -S(-\tau) & \tau (B^T Q(\theta) - R(-\tau,\theta))\\
*^T & *^T& - \tau \dot{S}(\theta)}\\
&M(\theta,\omega)=\frac{d}{d \theta} R(\theta,\omega)+\frac{d}{d
\omega} R(\theta,\omega)\\
&D_{11}=P A+A^T P + Q(0)+Q(0)^T + S(0) + \epsilon I
\end{align*}\label{thm:single_delay}
\end{thm}
\begin{proof}
We use the following complete quadratic functional.
\begin{align*}
&V(\phi)=\phi(0)^T P \phi(0) + 2 \phi(0)^T
\int_{-\tau}^{0} Q(\theta) \phi(\theta) d \theta \\
&+ \int_{-\tau}^{0}\phi(\theta)^T S(\theta) \phi(\theta) d \theta +
\int_{-\tau}^{0} \int_{-\tau}^{0}\phi(\theta)^T R(\theta,
\omega)\phi(\omega) d \theta d \omega\\
&= \frac{1}{\tau} \int_{-\tau}^0 \bmat{\phi(0) \\ \phi(\theta)}^T
\bmat{P -\epsilon I& \tau Q(\theta)\\ \tau Q(\theta)^T & \tau S(\theta)} \bmat{\phi(0) \\
\phi(\theta) } d \theta\\
&+ \int_{-\tau}^0 \int_{-\tau}^0 \phi(\theta)^T R(\theta,\omega)
\phi(\omega) d \theta d \omega+\epsilon \norm{\phi(0)}^2 \ge
\epsilon \norm{\phi(0)}^2
\end{align*}
All that is required for asymptotic stability is strict negativity
of the the derivative. The Lie derivative of this functional along
trajectories of the system is given by the following.
\begin{align*}
&\dot{V}(\Gamma(\phi,0))=\phi(0)^T P(A \phi(0)+B \phi(-\tau))+(
\phi(0)^T A^T+\phi(-\tau)^T B^T)P \phi(0)\\
&+2\int_{-\tau}^0 \phi(\theta)^T Q(\theta)^T
(A\phi(0)+B\phi(-\tau))d \theta\\
&+2 \phi(0)^T \left(Q(0)\phi(0)-Q(-\tau)\phi(-\tau)-\int_{-\tau}^0
\dot{Q}(\theta) \phi(\theta)d \theta\right)\\
&+\phi(0)^T S(0) \phi(0) - \phi(-\tau)^T S(-\tau)
\phi(-\tau)-\int_{-\tau}^0 \phi(\theta)^T \dot{S}(\theta)
\phi(\theta) d \theta \\
& +\int_{-\tau}^0 \left(\phi(0)^T R(0,\theta)
\phi(\theta)-\phi(-\tau)^T R(-\tau,\theta) \phi(\theta)\right)d
\theta \\
&+\int_{-\tau}^0 \left(\phi(\theta)^T R(\theta,0) \phi(0) -
\phi(\theta)^T R(\theta,-\tau) \phi(-\tau) \right) d
\theta \\
&-\int_{-\tau}^0 \int_{-\tau}^0 \phi(\theta)\left(\frac{d}{d \theta}
R(\theta,\omega) +\frac{d}{d \omega} R(\theta,\omega)
\right)\phi(\omega) d
\theta d \omega\\
&=\frac{1}{\tau}\int_{-\tau}^0
\bmat{\phi(0)\\\phi(-\tau)\\\phi(\theta)}^T D_1(\theta)
\bmat{\phi(0)\\\phi(-\tau)\\\phi(\theta)}-\int_{-\tau}^0
\int_{-\tau}^0 \phi(\theta)M(\theta,\omega)\phi(\omega) d \theta d
\omega-\epsilon
\norm{\phi(0)}^2 \\
&\le-\epsilon \norm{\phi(0)}^2
\end{align*}
Therefore, if the conditions of the theorem hold, then the
functional is strictly decreasing along trajectories of the system,
proving asymptotic stability.
\end{proof}

\subsection{Procedure for the Case of a Single Delay}
In this subsection, we explicitly show how the conditions associated
with Theorem~\ref{thm:single_delay} can be tested using semidefinite
programming. As a caveat, we note that direct implementation of the
procedure presented in this subsection may not be the the most
numerically efficient method possible. However, from a conceptual
standpoint, it is the most direct and, in terms of accuracy, will
deliver the results contained in the numerical examples given later
in the text.

\paragraph{Problem:} Suppose we want to determine whether the
following delay-differential equation is stable for a certain value
of delay, $\tau$, using polynomial functions of degree $d=4$.
\[\dot{x}(t)=ax(t)+bx(t-\tau)
\]
\paragraph{Defining Polynomial Variables:} Our first step is to
construct the variables which correspond to $P$ and polynomials $Q$
and $S$. Declaring the scalar variable $P$ should be clear. To
declare the polynomial matrices $Q$ and $S$ requires that we first
construct the vector of monomials which are contained in these
polynomials. Note that these vectors are not variables, rather
simply structures which will assist in constructing equality
constraints. Specifically, let
\[Z^1(\theta)=\bmat{1& \theta &\theta^2& \theta^3 & \theta^4}^T.
\]
Then $Q$ and $S$ can be represented using vector variables $v_Q$ and
$v_S$ as $Q(\theta)=v_Q^T Z^1(\theta)$ and $S(\theta)=v_S^T
Z^1(\theta)$. For convenience, we also put $P$ in vector form, so
that $P=v_P^T Z^1(\theta)$, where
\[v_P=\bmat{P&0&0&0&0}.
\]

\paragraph{Constructing elements of $G_1^4$:} For $G_1^4$, we
declare SOS variable $M$ and polynomial variable $T$. We first
declare vector variable $v_T$ and define $T(\theta)=v_T^T
Z^1(\theta)$. We then impose the constraint $\int_{-\tau}^0
T(\theta) d\theta=0$ as $v_T^T c^i=0$, where
\[c^i=\bmat{-\tau& 1/2\tau^2 &-1/3\tau^3 &1/4\tau^4 &-1/5\tau^5}
\]
Since
\begin{align*}
\int_{-1}^0 T(\theta) d\theta&=\int_{-1}^0 v_T^T Z^1(\theta)
d\theta=v_T^T \int_{-1}^0 Z^1(\theta) d\theta = v_T^T c^i =0
\end{align*}

To declare SOS variable $M$, we declare the positive semidefinite
matrix variable $m_M \ge 0$ and construct the vector of monomials
$Z^2$.
\[Z^2(\theta)=\bmat{Z^2_s(\theta)&0\\0&Z_s^1(\theta)}^T
\]
where
\[Z^2_s(\theta)=\bmat{1&\theta&\theta^2}
\]
We then define $M(\theta)=Z^2(\theta)^T m_M Z^2(\theta)$. By
construction, $M \in \Sigma_s$. For future use, it will be
convenient to write $M$ in vector form. Specifically, let
\[M(\theta)=\bmat{M_{11}(\theta) &
M_{12}(\theta)\\M_{21}(\theta)&M_{22}(\theta)} \qquad
\text{and}\qquad m_M=\bmat{m_{M11} & m_{M12}\\m_{M21} & m_{M22}}.
\]

Now let $v_{Mij}$ be the vector of columns of $m_{Mij}$. Then
\[M_{ij}(\theta)=v_{Mij}^T C_1 Z^1(\theta)
\]
where
%\[C_1=\bmat{1&&&&\\&1&&&\\&&1&&\\&1&&&\\&&1&&\\&&&1&\\&&1&&\\&&&1&\\&&&&1}
%\]
\[C_1=\bmat{1&&&&&&&&\\&1&&1&&&&&\\&&1&&1&&1&&\\&&&&&1&&1&\\&&&&&&&&1}^T.
\]

\textbf{Note:} Theorem~\ref{lem:H1_LMI} only requires $M(\theta) \ge
0$ for $\theta \in [-\tau,0]$. However, by imposing the constraint
$M \in \Sigma_s$, we are constraining $M(\theta)\ge 0$ for $\theta
\in \R$. Although this difference should not matter in a technical
sense since we only consider the interval $[-\tau,0]$, in practical
terms, this constraint severely restricts the set of admissable
polynomial functions. For this reason, in implementation, we
actually use the set $\Gamma_{2,4}^{K_p}$, where
$p(\theta)=-\theta(\theta+\tau)$. In this case, we would declare an
auxiliary positive semidefinite matrix variable $m_{S} \ge 0$ and
define $M(\theta)=Z^2(\theta)^T m_M Z^2(\theta)-\theta(\theta+1)
Z^2(\theta)^T m_{S} Z^2(\theta)$. For simplicity, however, we will
not explicitly include $m_S$ in our discussion.

\paragraph{Equality constraints for $G_1^4$:}
We can now impose the constraint
\[\bmat{P -\epsilon I & \tau Q \\ \tau Q^T & \tau S}\in G_1^d.
\]
This is easily done using the vector form of $M$. Specifically, we
include the constraints
\begin{align*}
v_{P-\epsilon}-v_T=v_{M11}C_1,\qquad v_{Q}=v_{M12}C_1,\quad \text{
and }\quad v_{S}=v_{M22}C_1.
\end{align*}

\paragraph{Constructing elements of $G_2^4$:} To construct an
element of $G_2^4$, we simply declare a positive semidefinite matrix
variable $m_R$ and define $R(\theta,\omega)=Z^2(\theta)m_R
Z^2(\omega)$. For convenience, it will be helpful to represent $R$
in vector form as $R(\theta,\omega)=v_R^T Z^3(\theta,\omega)$, where
$v_R$ is the vector of column vectors of $m_R$ and
\[Z^3(\theta,\omega):=\bmat{1&\theta&\theta^2&\omega & \theta \omega & \theta^2 \omega& \omega^2 & \theta \omega^2& \theta^2
\omega^2}^T
\]

\paragraph{Constructing the Derivative:}
Recall the derivative is defined $D$ and $G$, which are given by the
following.
\begin{align*}
&D(\theta) =\\
&\bmat{P A+A^T P + Q(0)+Q(0)^T + S(0) & PB - Q(-\tau)
& \tau (A^T Q(\theta) - \dot{Q}(\theta)+ R(0,\theta))\\
*^T & -S(-\tau) & \tau (B^T Q(\theta) - R(-\tau,\theta))\\
*^T & *^T& - \tau \dot{S}(\theta)}\\
&G(\theta,\omega)=\frac{d}{d \theta} R(\theta,\omega)+\frac{d}{d
\omega} R(\theta,\omega)
\end{align*}

Constructing $D$ and $G$ is relatively easy. Decompose $D(\theta)$
as
\[D(\theta)=\bmat{D_{11}&D_{12}&D_{13}(\theta)\\ *^T&D_{22}&D_{23}(\theta)\\ *^T&
*^T&D_{33}(\theta)}.
\]
Then
\begin{align*}
D_{11}&=2Pa+(2 v_Q^T +v_S) C_0\\
D_{12}&=Pb-v_Q C_\tau\\
D_{22}&=-v_S C_\tau\\
D_{13}(\theta)&=v_{D13}^T Z^1(\theta)=\tau(a v_Q^T - v_Q^T C_d +v_R C_{0,1})Z^1(\theta)\\
D_{23}(\theta)&=v_{D23}^T Z^1(\theta)=\tau(b v_Q-v_R
C_{\tau,1})Z^1(\theta)\\
D_{33}(\theta)&=v_{D23}^T Z^1(\theta)=-\tau(v_S C_d) Z^1(\theta)
\end{align*}
where
\begin{align*}
C_0&=Z^1(0)=\bmat{1&0&0&0&0}\\
C_\tau&=Z^1(-\tau)=\bmat{1& -\tau &\tau^2& -\tau^3 & \tau^4}^T\\
C_d&=\bmat{&&&&\\1&&&&\\&2&&&\\&&3&&\\&&&4&}\qquad
C_{0,1}=\bmat{1&&&&\\&&&&\\&&&&\\&1&&&\\&&&&\\&&&&\\&&1&&\\&&&&\\&&&&}\qquad
C_{\tau,1}=\bmat{1&&&&\\-\tau&&&&\\\tau^2&&&&\\&1&&&\\&-\tau&&&\\&-\tau^2&&&\\&&1&&\\&&-\tau&&\\&&-\tau^2&&}.
%=\bmat{1&&&&&&&&\\&&&1&&&&&\\&&&&&&1&&\\&&&&&&&&\\&&&&&&&&}
\end{align*}
Furthermore,
$G(\theta,\omega)=v_R^T(C_{rd1}+C_{rd2})Z^3(\theta,\omega)=v_G^T
Z^3(\theta,\omega)$, where
\[C_{rd1}=\bmat{&&&&&&&&\\1&&&&&&&&\\&2&&&&&&&\\&&&&&&&&\\&&&1&&&&&\\&&&&2&&&&\\&&&&&&&&\\&&&&&&1&&\\&&&&&&&2&}\qquad\text{and}\qquad
C_{rd2}=\bmat{&&&&&&&&\\&&&&&&&&\\&&&&&&&&\\1&&&&&&&&\\&1&&&&&&&\\&&1&&&&&&\\&&&2&&&&&\\&&&&2&&&&\\&&&&&2&&&}.
\]

\paragraph{Equality constraints for $G_3^4$:} To impose $-D \in
G_3^4$, we construct positive semidefinite variable $m_K\ge0$ with
and auxiliary variables $v_{U11}$, $v_{U12}$, and $v_{U22}$ with
constraints $v_{Uij}^T c^i=0$. We then represent
\[m_K = \bmat{m_{K11}&m_{K12}&m_{K13}\\ *^T&m_{K22}&m_{K23}\\ *^T&
*^T&m_{K33}}.
\]
Furthermore, let $v_{Kij}$ be the vector form of $m_{Kij}$. Then we
impose the following constraints.
\begin{align*}
&-v_{D11}=v_{K11}C_1-v_{U11} \qquad &-v_{D12}=v_{K12}C_1-v_{U12} \qquad &-v_{D22}=v_{K22}C_1-v_{U22}\\
&-v_{D13}=v_{K13}C_1  &-v_{D23}=v_{K23}C_1\qquad
&-v_{D33}=v_{K33}C_1,
\end{align*}
where
\begin{align*}
v_{D11}&=\bmat{D_{11}&0&0&0&0}\\
v_{D12}&=\bmat{D_{12}&0&0&0&0}\\
v_{D22}&=\bmat{D_{22}&0&0&0&0}.
\end{align*}

\paragraph{Equality constraints for $G_2^4$:}
Finally, we constrain $G \in G_2^4$. For this, we declare a new
positive semidefinite matrix variable $m_L\ge 0$ with vector form
$v_L$ and impose the following constraint.
\[v_G=v_L
\]

\paragraph{Conclusion} Once all the variables and constraints have
been declared, the problem can be submitted to SDP solvers such as
SeDuMi for evaluation. Although the problem presented in this
subsection is limited in scope, it will hopefully provide guidance
for those who are interested in understanding the mechanics of the
problem. For those solely interested in implementing the algorithm
at a higher level, there exist a number of alternatives. In
particular, SOSTools~\cite{sostools_web} provides a high-level
interface to solving sum-of-squares problems. An advantage of this
software package is the ability to declare polynomial variables and
implement constraints using Matlab's symbolic math toolbox without
creating the explicit matrix maps of the form $C$ which were used so
prominently in this example. The only drawback to using this
software lies in the fact that current implementations are somewhat
inefficient in constructing these constraints and so generally the
problem will take much longer computationally to set up than to
actually solve. For those who would like to simply type in a
dynamical system and receive an answer, we also suggest the Matlab
toolbox we have created for solving these types of problems and
which can be found online at
{\tt{\verb+http://www.stanford.edu/~mmpeet/software+}}.

\subsection{Examples of a Single Delay}

\textbf{Example 1:} In this example, we compare our results with the
discretized Lyapunov functional approach used by Gu et al.
in~\cite{gu_2003} in the case of a system with a single delay.
Although numerous other papers have also given sufficient conditions
for stability of time-delay systems, e.g.
\cite{fridman_2002,moon_2001,lu_2002}, we use the approach
introduced by Gu since it has demonstrated a particularly high level
of precision. When we are comparing with the piecewise linear
approach here and throughout this chapter, we will only consider
examples which have been presented in the work~\cite{gu_2003} and we
will compare our results with the numbers that are cited therein. We
use SOSTools~\cite{sostools_web} and SeDuMi~\cite{sturm_1999} for
solution of all semidefinite programming problems. Now consider the
following system of delay differential equations.
\[\dot{x}(t)=\bmat{0 & 1\\ -2 & 0.1}x(t)+\bmat{0 & 0\\ 1 & 0}x(t-\tau)
\]
The problem is to estimate the range of $\tau$ for which the
differential equation remains stable. Using the method presented in
this chapter and by sweeping $\tau$ in increments of $.1$, we
estimate the range of stability to be an interval. We then use a
bisection method to find the minimum and maximum stable delay. Our
results are also summarized in Table~\ref{Table1} and are compared
to the analytical limit and the results obtained in~\cite{gu_2003}.
\begin{table}[!ht]
\centering
\begin{tabular}[!]{|c|c|c|}
\hline
\multicolumn{3}{|c|}{Our results}\\
\hline
$d$ & $\tau_{\min}$ & $\tau_{\max}$\\
\hline
2 & .10017 & 1.6249\\
\hline
4 & .10017 & 1.7172\\
\hline
6 & .10017 & 1.71785\\
\hline
Analytic & .10017 & 1.71785\\
\hline
\end{tabular}
\begin{tabular}[!]{|c|c|c|}
\hline
\multicolumn{3}{|c|}{Piecewise Functional}\\
\hline
$N_2$& $\tau_{\min}$ & $\tau_{\max}$\\
\hline
1 & .1006 & 1.4272\\
\hline
2 & .1003 & 1.6921\\
\hline
3 & .1003 & 1.7161\\
\hline
 &  & \\
\hline
\end{tabular}
\vspace{0.1cm} \caption{$\tau_{\max}$ and $\tau_{\min}$ for
discretization level $N_2$ using the piecewise-linear Lyapunov
functional and for degree $d$ using our approach and compared to the
analytical limit} \label{Table1}
\end{table}

Clearly, the results for this test case illustrate a high rate of
convergence to the analytical limit. However, although the results
presented here give reasonable estimates for the interval of
stability, they do not prove stability over any interval, but rather
only at the specific values of $\tau$ for which the algorithm was
tested. To provide a more rigorous analysis, we now include $\tau$
as an uncertain parameter and search for parameter dependent
Lyapunov functionals which prove stability over an interval. The
results from this test are given in Table~\ref{Table3}.
\begin{table}[!ht]
\centering
\begin{tabular}[!]{|c|c|c|c|}
\hline
$d$ in $\tau$ & $d$ in $\theta$ & $\tau_{\min}$ & $\tau_{\max}$\\
\hline
2&2 & .1002 & 1.6246\\
\hline
2&4 & .1002 & 1.717\\
\hline
Analytic && .10017 & 1.71785\\
\hline
\end{tabular}\vspace{0.1cm} \caption{Stability on the interval
$[\tau_{\min}, \tau_{\max}]$ vs. degree using a parameter-dependent
functional} \label{Table3}
\end{table}

\begin{figure}[htb]
  \centerline{\includegraphics[width=0.3\textwidth]{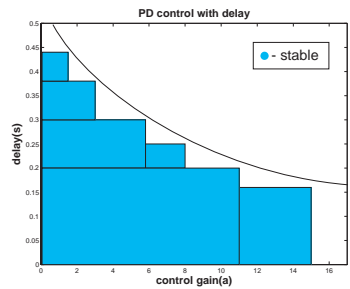}}
  \caption{Regions of Stability for Example 2}
  \label{fig:regions}
\end{figure}

\textbf{Example 2:} In this example, we illustrate the flexibility
of our algorithm through a simplistic control design and analysis
problem. Suppose we wish to control a simple inertial mass remotely
using a PD controller. Now suppose that the derivative control is
half of the proportional control. Then we have the following
dynamical system.
\[\ddot{x}(t)=-a x(t) -\frac{a}{2}\dot{x}(t)
\]
It is easy to show that this system is stable for all positive
values of $a$. However, because we are controlling the mass
remotely, some delay may be introduced due to, for example, the
fixed speed of light. We assume that this delay is known and changes
sufficiently slowly so that for the purposes of analysis, it may be
taken to be fixed. Now we have the following delay-differential
equation with uncertain, time-invariant parameters $a$ and $\tau$.
\[\ddot{x}(t)=-a x(t-\tau) -\frac{a}{2}\dot{x}(t-\tau)
\]
Whereas before the system was stable for all positive values of $a$,
now, for any fixed value of $a$, there exists a $\tau$ for which the
system will be unstable. In order to determine which values of $a$
are stable for any fixed value of $\tau$, we divide the parameter
space into regions of the form $a \in [a_{\min},a_{\max}]$ and $\tau
\in [\tau_{\min},\tau_{\max}]$. This type of region is compact and
can be represented as a semi-algebraic set using the polynomials
$p_1(a)=(a - a_{\min})(a-a_{\max})$ and $p_2(\tau)=(\tau -
\tau_{\min})(\tau-\tau_{\max})$. By using these polynomials, we are
able to construct parameter dependent Lyapunov functionals which
prove stability over a number of parameter regions. These regions
are illustrated in Figure~\ref{fig:regions}.

\subsection{Multiple Delay Case}

We now consider the case of multiple delays. The system is now
defined by the following functional differential equation where $A_i
\in \R^{n \times n}$, $\tau_{i}>\tau_{i-1}$ for $i=1,\ldots,K$, and
$\tau_0=0$.
\begin{equation} \dot{x}(t)=\sum_{i=0}^K A_i x(t-\tau_i)
\label{eqn:multiple_delay}
\end{equation}
The difference between the analysis of the case of a single delay
and that of multiple delays is that the matrix-valued functions
defining the complete quadratic functional necessary for stability
may now contain discontinuities at discrete points given by the
values of the delay. In this case, we express the the complete
quadratic functional in the following form where $Q_i:\R \rightarrow
\R^{n \times n}$, $S_i: \R \rightarrow \S^n$ and $R_{ij}: \R^2
\rightarrow \R^{n \times n}$ for $i,j = 1 ,\ldots, K$ are continuous
matrix-valued functions, $P \in \S^n$ and
$R_{ij}(\theta,\omega)=R_{ji}(\omega,\theta)^T$.
\begin{align*}
&V(\phi)=\phi(0)^T P \phi(0) +2 \sum_{j=1}^K \phi(0)^T
\int_{-\tau_i}^{-\tau_{i-1}} Q_i(\theta) \phi(\theta) d \theta \\
&+ \sum_{j=1}^K\int_{-\tau_i}^{-\tau_{i-1}}\phi(\theta)^T
S_i(\theta) \phi(\theta) d \theta +
\sum_{i=1}^{K}\sum_{j=1}^{K}\int_{-\tau_i}^{-\tau_{i-1}}
\int_{-\tau_j}^{-\tau_{j-1}}\phi(\theta)^T R_{ij}(\theta,
\omega)\phi(\omega) d \theta d \omega
\end{align*}

\begin{thm}
For integer $d\ge0$, suppose $P\in \S^n$, $\eta
>0$, $Q_i: \R\rightarrow \R^{n
\times n}$, $S_i:\R \rightarrow \S^n$ and $R_{ij}:\R^2 \rightarrow
\R^{n \times n}$ where
$R_{ij}(\theta,\omega)=R_{ji}(\omega,\theta)^T$ are continuous for
$i,j=1 ,\ldots, K$. Let $R(\theta,\omega)=R_{ij}(\theta,\omega)$ for
$\theta\in I_i$, $\omega \in I_j$ where $I_i=[-\tau_i,-\tau_{i-1}]$.
Then the solution map defined by equation~\eqref{eqn:multiple_delay}
is asymptotically stable if the following holds.
\begin{align*}M \in \tilde{G}_1^d & \qquad R \in \tilde{G}_2^d\\
-D \in \tilde{G}_3^d & \qquad L \in \tilde{G}_2^d
\end{align*}
\begin{align*}
&M(\theta)=\bmat{P-\eta I & \tau Q_i(\theta)\\\tau Q_i(\theta)^T
&\tau S_i(\theta)}
\qquad \theta \in I_i\\
&L(\theta,\omega)=\frac{\delta}{\delta
\theta}R_{ij}(\theta,\omega)+\frac{\delta}{\delta
\theta}R_{ij}(\theta,\omega)\qquad \theta \in I_i,\omega \in I_j\\
&D(\theta)=\bmat{D11 & \tau D12_{i}(\theta)\\
\tau D12_{i}(\theta)^T & \tau D22_{i}(\theta)} \qquad \text{for } \theta \in I_i\\
&D11_{11}=P A_0+A_0^T P
+Q_1(0)+Q_1(0)^T+S_1(0)-\eta I\\
&D11_{ij}=\\
&\begin{cases}
P A_{i-1}-Q_{i-1}(-\tau_{i-1})+Q_{i}(-\tau_{i-1})&i,j=1,2\ldots K\\
S_{i}(-\tau_{i-1})-S_{i-1}(-\tau_{i-1}) & i=j=2\ldots K\\
P A_K -Q_K(-\tau_K) & i,j=1,K+1\\
-S_K(-\tau_K)& i=j=K+1\\
0 & \text{otherwise}
\end{cases}\\
&D12_j(\theta)=\\
&\bmat{ R_{1j}(0,\theta)+A_0^T Q_j(\theta)-\dot{Q}_j(\theta)\\
R_{2j}(-\tau_1,\theta)-R_{1j}(-\tau_1,\theta)+A_1^T Q_j(\theta) \\ \vdots \\
R_{(K)j}(-\tau_{K-1},\theta)-R_{(K-1)j}(-\tau_{K-1},\theta)+A_{K-1}^T Q_j(\theta) \\
A_K^T Q_j(\theta)-R_{Kj}(-\tau_K,\theta)}\\
&D22_i(\theta)=-\dot{S}_i(\theta)
\end{align*}
\end{thm}
\begin{proof}
We use the following complete quadratic functional.
\begin{align*}
&V(\phi)=\frac{1}{\tau}\int_{-\tau_K}^0\bmat{\phi(0)\\\phi(\theta)}^T
M(\theta)
\bmat{\phi(0)\\\phi(\theta)}\\
&+\int_{-\tau_K}^0\int_{-\tau_K}^0 \phi(\theta)^T
R(\theta,\omega)\phi(\omega)d \theta d \omega + \eta
\norm{\phi(0)}^2 \ge \eta \norm{\phi(0)}^2
\end{align*}

The system is asymptotically stable if the derivative of the
functional is strictly negative. The Lie derivative of the
functional along trajectories of the systems is given by the
following.
\begin{align*}
&\dot{V}(\Gamma(\phi,0))
=\frac{1}{\tau}\int_{-\tau_K}^0\bmat{\phi(0)\\\vdots \\
\phi(-\tau_K)\\\phi(\theta)}^T D(\theta)\bmat{\phi(0)\\\vdots \\
\phi(-\tau_K)\\\phi(\theta)}d \theta\\
&+\int_{-\tau_K}^0 \phi(\theta)^T L(\theta,\omega)\phi(\omega)d
\theta d \omega-\eta\norm{\phi(0)}^2 \le -\eta \norm{\phi(0)}^2
\end{align*}
Thus if the conditions of the theorem hold, then the derivative of
the functional is strictly negative, which implies asymptotic
stability.
\end{proof}

\subsection{Examples of Multiple Delay}

\textbf{Example 3:} Consider the following system of
delay-differential equations.
%\[\dot{x}(t)=\bmat{-2 & 0\\ 0 & -0.9}x(t)+\bmat{-1 & 0\\ -1 & -1}[0.05x(t-0.5\tau)+0.95 x(t-\tau)]
%\]
\[\dot{x}(t)=\bmat{-2 & 0\\ 0 & -\frac{9}{10}}x(t)
+\bmat{-1 & 0\\ -1 &
-1}\left[\frac{1}{20}x(t-\frac{\tau}{2})+\frac{19}{20}
x(t-\tau)\right]
\]
Again, the system is stable when $\tau$ lies on some interval. The
problem is to search for the minimum and maximum value of $\tau$ for
which the system remains stable. In applying the methods of this
chapter, we again use a bisection method to find the minimum and
maximum value of $\tau$ for which the system remains stable. Our
results are summarized in Table~\ref{Table3} and are compared to the
analytical limit as well as piecewise-linear functional method. For
the piecewise functional method, $N_2$ is the level of both
discretization and subdiscretization.

\begin{table}[!ht]
\centering
\begin{tabular}[!]{|c|c|c|}
\hline
\multicolumn{3}{|c|}{Our Approach}\\
\hline
$d$ & $\tau_{\min}$ & $\tau_{\max}$\\
\hline
2 & .20247 & 1.354\\
\hline
4 & .20247 & 1.3722\\
\hline
Analytic & .20246 & 1.3723\\
\hline
\end{tabular}
\begin{tabular}[!]{|c|c|c|}
\hline
\multicolumn{3}{|c|}{Piecewise Functional}\\
\hline
$N_2$ & $\tau_{\min}$ & $\tau_{\max}$\\
\hline
1 & .204 & 1.35\\
\hline
2 & .203 & 1.372\\
\hline
 &  & \\
\hline
\end{tabular}
\vspace{0.1cm} \caption{$\tau_{max}$ and $\tau_{min}$ using the
piecewise-linear Lyapunov functional of Gu et al. and our approach
and compared to the analytical limit} \label{Table3}
\end{table}

\subsection{Distributed Delay Case}

We now consider the case of distributed delay where the dynamics are
given by the following functional differential equation where
$A_0\in \R^{n \times n}$ and $A\in \R^{n \times n}[\theta]$.
\begin{equation}
\dot{x}(t)= A_0 x(t)+\int_{- \tau}^{0}A(\theta)x(t+\theta)d \theta
\label{eqn:distributed}
\end{equation}
We can again assume that the matrix-valued functions defining the
complete quadratic functional are continuous. This leads to the
following theorem.
\begin{thm}
For integer $d\ge0$, the solution map, $G$, defined by
equation~\eqref{eqn:distributed} is asymptotically stable if there
exists a constant $\eta
> 0$, a matrix $P
\in \S^n$ and matrix functions $Q:\R \rightarrow \R^{n \times n}$,
$S: \R \rightarrow \S^{n}$, $R: \R^2 \rightarrow \R^{n \times n}$
where $R(\theta,\omega)=R(\omega,\theta)^T$ and such that the
following holds
\begin{align*}
&\bmat{P-\eta I & \tau Q\\ \tau Q^T & \tau S} \in G_1^d
&-D \in G_{3}^d\\
&R \in G_{2}^d & M \in G_2^d
\end{align*}
where
\begin{align*}
&D(\theta)=\bmat{D_{11} + \eta I& \tau D_{12}(\theta)\\ \tau D_{12}(\theta)^T & -\tau \dot{S}(\theta)}\\
&D_{11}= \bmat{P A_0+ A_0^T P +Q(0)+Q(0)^T+S(0)
& -Q(-\tau) \\ -Q(-\tau)^T & -S(-\tau)}\\
&D_{12}(\theta)=\bmat{A_0^T Q(\theta) +P A(\theta) -\frac{d}{d \theta}Q(\theta)+R(0,\theta)\\ -R(-\tau,\theta)}\\
&M(\theta,\omega)=\frac{d}{d \theta}R(\theta,\omega)+\frac{d}{d
\omega} R(\theta,\omega)-A(\theta)^T Q(\omega)-Q(\theta)^T A(\omega)
\end{align*}
\end{thm}
\begin{proof}
We consider the complete quadratic functional
\begin{align*}&V(\phi)
=\frac{1}{\tau}\int_{-\tau}^0\bmat{\phi(0)\\ \phi(\theta)}^T
\bmat{P-\eta I & \tau
Q(\theta)\\ \tau Q(\theta)^T & \tau S(\theta)} \bmat{\phi(0)\\
\phi(\theta)}^T d \theta \\
&+\int_{-\tau}^0 \int_{-\tau}^0 \phi(\theta)^T
R(\theta,\omega)\phi(\omega)d \theta d \omega+ \eta \norm{\phi(0)}^2
\ge \eta \norm{\phi(0)}^2
\end{align*}
The system is asymptotically stable if the Lie derivative is
strictly negative. The Lie derivative of the functional along
trajectories of the system is given by
\begin{align*}
&\dot{V}(\Gamma(\phi,0)) =\int_{-\tau}^{0}
\bmat{\phi(0)\\\phi(-\tau) \\ \phi(\theta) }^T
D(\theta)\bmat{\phi(0)\\\phi(-\tau)\\ \phi(\theta) } d \theta\\
&-\int_{-\tau}^{0} \int_{-\tau}^{0}  \phi(\theta)^T
M(\theta,\omega)\phi(\omega) d \omega d \theta-\eta \norm{\phi(0)}^2
\le -\eta \norm{\phi(0)}^2
\end{align*}
Thus, if the conditions of the theorem hold, then the system is
asymptotically stable.
\end{proof}

\section{Conclusion}

In this chapter, we have shown how to compute solutions to an
operator-theoretic version of the Lyapunov inequality. Our approach
is to combine results from real algebraic geometry with functional
analysis in order to parameterize certain convex cones of positive
operators using the convex cone of positive semidefinite matrices.
We then show that the operator inequalities can be expressed using
affine constraints on these matrices. This allows us to compute
solutions using semidefinite programming, for which there exist
efficient numerical algorithms. We have further extended our results
to the case when the dynamics contain parametric uncertainty using a
construction based on certain results of Putinar and others. The
numerical examples given in this Chapter demonstrate a quick
convergence to the analytic limit of stability.

We conclude by mentioning that the methods of this chapter can be
used to construct full-rank solutions of the Lyapunov inequality.
Furthermore, we have developed methods for constructing the inverses
of these full-rank operators. By computing solutions to the Lyapunov
inequality for the adjoint system constructed by Delfour and
Mitter~\cite{delfour_72a}, this invertibility result seems to imply
that one can construct stabilizing controllers for linear time-delay
systems. This work is ongoing, however, and is therefore not
detailed in this thesis.

\chapter{Stability of Nonlinear Time-Delay
Systems}\label{chp:nonlinearcase}
\section{Introduction}
The question of stability of nonlinear functional differential
equations is complicated by the lack of a complete Lyapunov
functional structure whose existence is necessary for stability of a
general nonlinear time-delay system. For this reason all results
given in this section will only be sufficient. Furthermore, we
cannot claim, as was done in the previous chapter, that these
conditions will generally approach necessity in any sense. Given
these limitations, we felt that the best approach was to construct a
sequence of conditions that would approach necessity in at least one
special case, that of linear systems. For this reason, the
conditions presented in this chapter are a generalization of those
in the previous one and the Lyapunov functionals reduce to the
complete quadratic functional in the case where the degree of $x$ is
restricted to $2$. This chapter differs from the presentation of the
linear case in that it includes a section on delay-independent
stability, a topic not previously discussed.

\section{Delay-Dependent Stability}\label{sec:delay_dep}
Consider a nonlinear discrete time-delay system defined by a
functional of the following form for $\tau_i<\tau_{i-1}$ for
$i=1,\ldots,K$ and $\tau_0=0$.
\begin{align*}
f(x_t)&= f(x(t),x(t-\tau_1),\cdots ,x(t-\tau_K))%\\
%z(x_t)&=\left[x_t(-\tau_0)^T,x_t(-\tau_1)^T,\cdots,x_t(-\tau_{K})^T
%\right]^T\\
%&=\left[x(t)^T,x(t-\tau_1)^T,\cdots,x(t-\tau_{K})^T \right]^T
\end{align*}

In this section, we present conditions for stability based on the
use of a generalization of the compete quadratic functional. The
generalized Lyapunov functional has the following form.

\begin{align*}
 V(\phi):=&\int_{-\tau_K}^0 Z_d[\phi(0),\phi(\theta)]^T M(\theta) Z_d[\phi(0),\phi(\theta)] d \theta \\
 &+ \int_{-\tau_K}^0 \int_{-\tau_K}^0
 Z_d[\phi(\theta)]^T R(\theta,\omega) Z_d[\phi(\omega)]d \theta d
 \omega
\end{align*}

\subsection{Single Delay Case}
In this subsection, we consider the special case of a single delay.
\[\dot{x}(t)=f(x(t),x(t-\tau))
\]

Here we assume $x(t) \in \R^n$ and $f$ is continuous. We first
introduce some notation.

\begin{defn}We say a continuous function $f : \R^{2n} \times \R \rightarrow \R$ satisfies
$f \in K_1$ if there exists some $\alpha>0$ such that the following
holds for all $\phi \in \mathcal{C}_\tau$.
\[\int_{-\tau}^0 f(\phi(0),\phi(\theta),\theta) d\theta \ge \alpha
\norm{\phi(0)}^2
\]
\end{defn}

\begin{defn}We say a continuous function $f : \R^{3n} \times \R^2 \rightarrow \R$ satisfies
$f \in K_2$ if the following holds for all $\phi \in
\mathcal{C}_\tau$.
\[\int_{-\tau}^0 \int_{-\tau}^0
f(\phi(\theta),\phi(\omega),\theta,\omega) d\theta d \omega \ge 0
\]
\end{defn}

\begin{defn}We say a continuous function $f : \R^{3n} \times \R \rightarrow \R$ satisfies
$f \in K_3$ if there exists some $\alpha>0$ such that the following
holds for all $\phi \in \mathcal{C}_\tau$.
\[\int_{-\tau}^0 f(\phi(0),\phi(-\tau),\phi(\theta),\theta) d\theta \ge \alpha
\norm{\phi(0)}^2
\]
\end{defn}

The sets $K_1$, $K_2$, and $K_3$ are used to define a positive
Lyapunov functional and some of its derivative forms. For a given
degree bound, these sets can be parameterized by the space of
positive semidefinite matrices using the following subsets.

%\begin{defn} For variables $x$ and integer $m\ge 0$, define $Z_m(x)$ to be
%the vector of monomials, lexigraphically ordered, in variables $x$
%of degree less than or equal to $m$.
%\end{defn}
%
%For $n$ variables, $Z_m(x)$ is of length $n+m \choose m$.

\begin{defn} Denote $f \in \Xi^d_1 \subset K_1$ if $f$ is a polynomial of degree
$d$ or less such that there exists a polynomial function
$t:\R^{n}\times \R \rightarrow \R$, of degree $d$ or less, and
$\alpha>0$ such that the following holds.
\begin{align*}
f(x_1,x_2,\theta)-t(x_1,\theta)-\alpha\norm{x_1}^2 &\in \Sigma_s \\
%s(x_1,x_2,\theta)&\ge 0 \quad \forall \theta \in [-\tau,0]\\
\int_{-\tau}^0 t(x_1,\theta) d \theta &=0
\end{align*}
\end{defn}

\begin{defn}
Denote $f \in \Xi^{2d}_2 \subset K_2$ if $f$ is a polynomial of
degree $2d$ or less and there exists a polynomial matrix function
$R:\R^2\rightarrow \R^{{n+d \choose d}\times {n+d \choose d}}$ such
that $R \in G_2^{2d}$(defined in Chapter~\ref{chp:LinearCase}) and
the following holds.
\begin{align}
f(x_1,x_1,\theta,\omega) =Z_d(x_1)^T R(\theta,\omega) Z_d(x_2)
\end{align}
\end{defn}

\begin{defn}
Denote $f \in \Xi_3^d \subset K_3$ if $f$ is a polynomial of degree
$d$ or less such that there exist polynomial function
$t:\R^{2n}\times \R \rightarrow \R$, of degree $d$ or less, and
$\alpha>0$ such that the following holds.
\begin{align*}
f(x_1,x_2,x_3,\theta)-t(x_1,x_2,\theta)-\alpha\norm{x_1}^2&\in \Sigma_s \\
%s(x_1,x_2,x_3,\theta)&\ge 0 \quad \forall \theta \in [-\tau,0]\\
\int_{-\tau}^0 t(x_1,x_2,\theta) d \theta &=0
\end{align*}
\end{defn}

We now state the stability theorem.

\begin{thm}\label{thm:nonlinear_single}
For a given integer $d\ge 0$, suppose there exist functions $g\in
\Xi_1^d$, $h,\hat{h} \in \Xi^d_2$, and $-\hat{g} \in \Xi^d_3$ such
that the following conditions hold.
\begin{align*}
&\hat{g}(x_t(0),x_t(-\tau),x_t(\theta),\theta)\\
&\qquad=g(x_t(0),x_t(0),0) -g(x_t(0),x_t(-\tau),-\tau)\\
&\qquad+\tau \nabla_{x_t(0)} g(x_t(0),x_t(\theta),\theta)^T
f(x_t(0),x_t(-\tau))-\tau \frac{\delta}{\delta
\theta}g(x_t(0),x_t(\theta),\theta)\\
&\qquad+\tau h(x_t(0),x_t(\theta),0,\theta)-\tau h(x_t(-\tau),x_t(\theta),-\tau,\theta)\\
&\qquad+\tau h(x_t(\theta), x_t(0), \theta, 0)-\tau
h(x_t(\theta),x_t(-\tau),\theta,-\tau)
\\
&\hat{h}(x_t(\theta),x_t(\omega),\theta,\omega)=\frac{\delta}{\delta
\theta}h(x_t(\theta),x_t(\omega),\theta,\omega)+\frac{\delta}{\delta
\omega}h(x_t(\theta),x_t(\omega),\theta,\omega)
\end{align*}
Then the time-delay system defined by $f$ is globally asymptotically
stable.
\end{thm}
\begin{proof}
Consider the Lyapunov functional defined as follows
\begin{align*}
V(x_t):=\int_{-\tau}^0 g(x_t(0),x_t(\theta),\theta) d\theta +
\int_{-\tau}^0 \int_{-\tau}^0
h(x_t(\theta),x_t(\omega),\theta,\omega) d\theta d\omega \ge \alpha
\norm{x_t(0)}^2
\end{align*}
The derivative of this functional along trajectories of the systems
can be expressed as follows.
\begin{align*}
\dot{V}(x_t):= &g(x_t(0),x_t(0),0) -g(x_t(0),x_t(-\tau),-\tau)\\
&+\int_{-\tau}^0 \bbl( \nabla_{x_t(0)}
g(x_t(0),x_t(\theta),\theta)^T
p(x_t(0),x_t(-\tau))-\frac{\delta}{\delta
\theta}g(x_t(0),x_t(\theta),\theta)\bbr) d\theta\\
&+\int_{-\tau}^0\bbl(
h(x_t(0),x_t(\theta),0,\theta)-h(x_t(-\tau),x_t(\theta),-\tau,\theta)\bbr)d
\theta\\
&+\int_{-\tau}^0\bbl( h(x_t(\theta), x_t(0), \theta,
0)-h(x_t(\theta),x_t(-\tau),\theta,-\tau)\bbr)d
\theta\\
&-\int_{-\tau}^0 \int_{-\tau}^0 \bbl( \frac{\delta}{\delta
\theta}h(x_t(\theta),x_t(\omega),\theta,\omega)+\frac{\delta}{\delta
\omega}h(x_t(\theta),x_t(\omega),\theta,\omega) \bbr) d\theta
d\omega\\
=&\frac{1}{\tau}\int_{-\tau}^0
\hat{g}(x_t(0),x_t(-\tau),x_t(\theta),\theta)d\theta+\int_{-\tau}^0
\int_{-\tau}^0 \hat{h}(x_t(\theta),x_t(\omega),\theta,\omega)
d\theta d \omega\\
 \le& -\alpha \norm{x_t(0)}^2
\end{align*}
Thus, we have that if the conditions of the theorem are satisfied,
then the derivative of the Lyapunov functional is negative definite,
proving global asymptotic stability of the functional differential
equation defined by $f$.
\end{proof}

As with the linear case, the conditions contained in
Theorem~\ref{thm:nonlinear_single} can be expressed as a
semidefinite program.

\subsection{Multiple Delays}
In this section, we consider stability of nonlinear time-delay
systems defined by functionals of the following form for
$\tau_i<\tau_{i-1}$ for $i=1,\ldots,K$ and $\tau_0=0$.
\[\dot{x}(t)=f(x(t),x(t-\tau_1),\cdots ,x(t-\tau_K))
\]

As before, $x(t)\in \R^n$, $f$ is continuous, and we introduce some
notation.

\begin{defn}We say a function $f : \R^{2n} \times \R \rightarrow \R$ satisfies
$f \in \tilde{K}_1$ if $f(x_1,x_2,\theta)$ is continuous except
possibly at points $\theta \in \{-\tau_i\}_{i=1}^{K-1}$ and there
exists some $\alpha>0$ such that the following holds for all $\phi
\in \mathcal{C}_\tau$.
\[\int_{-\tau_K}^0 f(\phi(0),\phi(\theta),\theta) d\theta \ge \alpha
\norm{\phi(0)}^2
\]
\end{defn}

\begin{defn}We say a function $f : \R^{n} \times \R^2 \rightarrow \R$ satisfies
$f \in \tilde{K}_2$ if $f(x_1,x_2,\theta,\omega)$ is continuous
except possibly at points $\theta,\omega \in
\{-\tau_i\}_{i=1}^{K-1}$ and the following holds for all $\phi \in
\mathcal{C}_\tau$.
\[\int_{-\tau_K}^0 \int_{-\tau_K}^0
f(\phi(\theta),\phi(\omega),\theta,\omega) d\theta d \omega \ge 0
\]
\end{defn}

\begin{defn}We say a function $f : \R^{n(K+2)} \times \R \rightarrow \R$ satisfies
$f \in \tilde{K}_3$ if $f$ is continuous except possibly at points
$\theta \in \{-\tau_i\}_{i=1}^K $ and there exists some $\alpha>0$
such that the following holds for all $\phi \in \mathcal{C}_\tau$.
\[\int_{-\tau_K}^0 f(\phi(0),\phi(-\tau_1),\cdots,\phi(-\tau_K),\phi(\theta),\theta) d\theta \ge \alpha
\norm{\phi(0)}^2
\]
\end{defn}

The sets $\tilde{K}_1$, $\tilde{K}_2$, and $\tilde{K}_3$ are used to
define positive Lyapunov functionals and some of its derivative
forms in the case of multiple delays. For a given degree bound,
these sets can be parameterized by the space of positive
semidefinite matrices using the following subsets.

\begin{defn} Denote $f \in \tilde{\Xi}_1^d \subset \tilde{K}_1$ if there exist polynomial functions
$f_i:\R^{2n}\times \R \rightarrow \R$, $t_i:\R^{n}\times \R
\rightarrow \R$, $i=1 \cdots K $ of degree $d$ or less and
$\alpha>0$ such that the following holds for all $x_i \in \R^n$.
\begin{align*}
&f_i(x_1,x_2,\theta)-t_i(x_1,\theta)-\alpha\norm{x_1}^2 \in \Sigma_s \\
&f(x_1,x_2,\theta)=f_i(x_1,x_2,\theta) \quad \forall \theta \in
[-\tau_i,-\tau_{i-1}]\\
&\sum_{i=1}^{K} \int_{-\tau_i}^{-\tau_{i-1}} t_i(x_1,\theta) d
\theta =0
\end{align*}
\end{defn}

\begin{defn}
Denote $f \in \tilde{\Xi}_2^{2d} \subset \tilde{K}_2$ if $f$ is a
polynomial of degree $2d$ or less such there exists a matrix
function $R:\R^2 \rightarrow \R^{{n+d \choose d} \times {n+d \choose
d}}$ such that $R \in \tilde{G}_2^{2d}$ and the following holds.
\begin{align}
f(x_1,x_1,\theta,\omega) =Z_d(x_1)^T R(\theta,\omega) Z_d(x_2)
\end{align}
\end{defn}

\begin{defn}
Denote $f \in \tilde{\Xi}_3^{d} \subset \tilde{K}_3$ if there exist
functions $f_i:\R^{n(K+2)}\times \R \rightarrow \R$,
$t_i:\R^{n(K+2)}\times \R \rightarrow \R$, $i=1 \cdots K$, of degree
$d$ or less, and $\alpha>0$ such that the following holds.
\begin{align*}
&f_i(x_1,\cdots,x_{K+1},x_{K+2},\theta)-t_i(x_1,\cdots,x_{K+1},\theta)-\alpha\norm{x_1}^2
\in \Sigma_s \\
&f(x_1,\cdots,x_{K+1},x_{K+2},\theta)=f_i(x_1,\cdots,x_{K+1},x_{K+2},\theta)
\quad \forall \theta
\in [-\tau_i,-\tau_{i-1}]\\
&\sum_{i=1}^K \int_{-\tau_{i}}^{-\tau_{i-1}}
t_i(x_1,\cdots,x_{K+1},\theta) d \theta =0
\end{align*}
\end{defn}

We now state the stability theorem in the case of multiple delays.

\begin{thm} \label{thm:nonlinear_multiple}
For a given integer $d\ge0$, suppose there exist functions $h\in
\tilde{\Xi}^d_1$, $g,\hat{g} \in \tilde{\Xi}^d_2$, and $\hat{h} \in
\tilde{\Xi}^d_3$ such that the following conditions hold where
$I_i:=[-\tau_i,-\tau_{i-1}]$ and $\Delta
\tau_{i}:=\tau_i-\tau_{i-1}$.
\begin{align*}
&\hat{g}_i(x_t(0),\cdots,x_t(-\tau_K),x_t(\theta),\theta)\\
&\qquad=g_i(x_t(0),x_t(-\tau_{i-1}),-\tau_{i-1})-g_i(x_t(0),x_t(-\tau_{i}),-\tau_{i})\\
&\qquad+\Delta \tau_{i}
\nabla_{x_t(0)}g_i(x_t(0),x_t(\theta),\theta)^T
f(x_t(0),\cdots,x_t(-\tau_K)) - \Delta
\tau_{i}\frac{\delta}{\delta \theta}g_i(x_t(0),x_t(\theta),\theta)\\
&\qquad+ \sum_{j=1}^K \bbl(
h_{ij}(x_t(-\tau_{i-1},x_t(\theta),-\tau_{i-1},\theta)
-h_{ij}(x_t(-\tau_{i}),x_t(\theta),-\tau_{i},\theta) \bbr)\\
&\qquad + \sum_{j=1}^K
\bbl(h_{ji}(x_t(\theta),x_t(-\tau_{j-1}),\theta,-\tau_{j-1})-h_{ji}(x_t(\theta),x_t(-\tau_j),\theta,-\tau_j)\bbr)
\\
&\hat{h}_{ij}(x_t(\theta),x_t(\omega),\theta,\omega)=\frac{\delta}{\delta
\theta}
h_{ij}(x_t(\theta),x_t(\omega),\theta,\omega)+\frac{\delta}{\delta
\omega} h_{ij}(x_t(\theta),x_t(\omega),\theta,\omega)
\end{align*}
Where
\begin{align*}
&\Delta \tau_{i}
g_i(x_t(0),x_t(\theta),\theta):=g(x_t(0),x_t(\theta),\theta)
\quad \forall \theta \in I_i\\
&\hat{g}_i(x_t(0),\cdots,x_t(-\tau_K),x_t(\theta),\theta):=\hat{g}(x_t(0),\cdots,x_t(-\tau_K),x_t(\theta),\theta)
\quad \forall \theta \in I_i\\
&h_{ij}(x_t(\theta),x_t(\omega),\theta,\omega):=h(x_t(\theta),x_t(\omega),\theta,\omega)\quad \forall \theta,\omega \in I_i\\
&\hat{h}_{ij}(x_t(\theta),x_t(\omega),\theta,\omega):=\hat{h}(x_t(\theta),x_t(\omega),\theta,\omega)\quad \forall \theta,\omega \in I_i\\
\end{align*}

Then the time-delay system defined by $f$ is globally asymptotically
stable.
\end{thm}
\begin{proof}
Consider the Lyapunov functional defined as follows.
\begin{align*}
&V(x_t):=\int_{-\tau_K}^0 g(x_t(0),x_t(\theta),\theta) d\theta +
\int_{-\tau_K}^0 \int_{-\tau_K}^0
h(x_t(\theta),x_t(\omega),\theta,\omega) d\theta d\omega\\
&=\sum_{i=1}^K \int_{-\tau_i}^{-\tau_{i-1}} \Delta
\tau_{i}g_i(x_t(0),x_t(\theta),\theta) d\theta + \sum_{i,j=1}^K
\int_{-\tau_i}^{-\tau_{i-1}} \int_{-\tau_j}^{-\tau_{j-1}}
h_{ij}(x_t(\theta),x_t(\omega),\theta,\omega) d\theta d\omega \\
&\ge \alpha \norm{x_t(0)}^2
\end{align*}
The inequality holds for some $\alpha>0$ by definition of the set
$\tilde{K}_1$. The derivative of this functional along trajectories
of the system can be expressed as follows.
\begin{align*}
&\dot{V}(x_t):=\sum_{i=1}^K \Delta
\tau_{i} \bbl(g_i(x_t(0),x_t(-\tau_{i-1}),-\tau_{i-1})-g_i(x_t(0),x_t(-\tau_{i}),-\tau_{i}) \bbr)\\
&+\sum_{i=1}^K \int_{-\tau_i}^{-\tau_{i-1}} \Delta \tau_{i} \bbl(
\nabla_{x_t(0)} g_i(x_t(0),x_t(\theta),\theta)^T
f(x_t(0),\cdots,x_t(-\tau_K)) - \frac{\delta}{\delta \theta} g_i(x_t(0),x_t(\theta),\theta) \bbr) d\theta\\
&+\sum_{i,j=1}^K \int_{-\tau_i}^{-\tau_{i-1}}
\bbl( h_{ij}(x_t(-\tau_{j-1},x_t(\omega),-\tau_{j-1},\omega) -h_{ij}(x_t(-\tau_{j}),x_t(\omega),-\tau_{j},\omega) \bbr)d\omega\\
&+\sum_{i,j=1}^K \int_{-\tau_j}^{-\tau_{j-1}}
\bbl( h_{ij}(x_t(\theta),x_t(-\tau_{i-1}),\theta,-\tau_{i-1})-h_{ij}(x_t(\theta),x_t(-\tau_i),\theta,-\tau_i) \bbr)d\theta\\
&- \sum_{i,j=1}^K \int_{-\tau_i}^{-\tau_{i-1}}
\int_{-\tau_j}^{-\tau_{j-1}} \bbl(\frac{\delta}{\delta \theta}
h_{ij}(x_t(\theta),x_t(\omega),\theta,\omega)+\frac{\delta}{\delta
\omega}
h_{ij}(x_t(\theta),x_t(\omega),\theta,\omega)\bbr) d\theta d\omega\\
&=\frac{1}{\tau_K}\sum_{i=1}^K \int_{-\tau_i}^{-\tau_{i-1}}
\hat{g}_i(x_t(0),\cdots,x_t(-\tau_K),x_t(\theta),\theta) d\theta\\
&+ \sum_{i,j=1}^K \int_{-\tau_i}^{-\tau_{i-1}}
\int_{-\tau_j}^{-\tau_{j-1}}
\hat{h}_{ij}(x_t(\theta),x_t(\omega),\theta,\omega) d\theta d\omega\\
&=\int_{-\tau}^0 \hat{g}(x_t(0),\cdots,x_t(-\tau_K),\theta) d\theta
+ \int_{-\tau}^0 \int_{-\tau}^0
\hat{h}(x_t(\theta),x_t(\omega),\theta,\omega) d\theta d\omega \le
-\alpha \norm{x_t(0)}^2
\end{align*}
The final inequality follows from the definition of $\tilde{K}_2$
and $\tilde{K}_3$. Thus the conditions of the theorem imply global
asymptotic stability of the time-delay system defined by $f$.
\end{proof}

Similar to the case of a single delay, the conditions of
Theorem~\ref{thm:nonlinear_multiple} can be expressed as a
semidefinite program.

\section{Delay-Independent Stability}\label{sec:delay_ind}
We now briefly discuss conditions for delay-independent stability of
time-delay systems. Recall that a system defined by a functional of
the following form is stable independent of delay if it is stable
for any finite values of $\{\tau_i\}_{i=1}^K$.
\[\hat{f}(x_t):=f(x_t(0),x_t(-\tau_1),\cdots,x_t(-\tau_K))
\]

Here we assume $f$ is continuous and $x(t)\in \R^n$. The theorems of
the previous section can be applied directly to the problem of
delay-independent stability by considering $\tau_i$ to be an
uncertain parameter and searching for a parameter-dependent Lyapunov
functional. However, a less computationally intensive approach is to
consider specific Lyapunov functional forms where $\tau_i$ does not
explicitly appear in the stability conditions. One such Lyapunov
functional form is given as follows.

\[
V(\phi)=p_0(\phi(0))+\sum_{i=1}^K \int_{-\tau_i}^{0}
p_i(\phi(\theta)) \,d \theta
\]
This functional has an upper Lie derivative defined as follows.
\[
\dot{V}(\phi)= \nabla p_0(\phi(0))^T
f(x_t(0),x_t(-\tau_1),\cdots,x_t(-\tau_K))+\sum_{i=1}^K
p_i(\phi(0))-p_2(\phi(-\tau_i))
\]

Since $\tau_i$ does not appear explicitly in the derivative of the
functional, we can prove stability for arbitrary $\tau_i$ using the
following theorem.

\begin{thm}\label{thm:nonlinear_multiple_ind}
For a given integer $d$, suppose there exist polynomials
$\{p_i\}_{i=0}^K$ and constant $\alpha>0$ such that $p_i(0)=0$ for
$i=0,\ldots, K$ and the following conditions hold.
\begin{align*}
&p_0(x_0)-\alpha \norm{x_0}^2 \in \Sigma_s^d\\
&p_i(x_i)\in \Sigma_s^d \qquad i=1,\ldots,K\\
&-\nabla p_0(x_0)^T f(x_0,x_1,\cdots,x_K)-\sum_{i=1}^K
\bbl(p_i(x_0)-p_2(x_i)\bbr) - \alpha \norm{x_0}^2 \in \Sigma_s
\end{align*}
Then the system defined by $f$ is globally asymptotically stable for
any finite values of $\tau_i$.
\end{thm}
\begin{proof}
Consider the following Lyapunov functional.
\[
V(x_t)=p_0(x_t(0))+\sum_{i=1}^K \int_{-\tau_i}^{0} p_i(x_t(\theta))
\,d \theta \ge \alpha \norm{x_t(0)}^2
\]
The Lyapunov functional is positive definite by the conditions of
the theorem. The derivative of the functional along trajectories of
the system is given as follows.
\begin{align*}
\dot{V}(\phi)&= \nabla p_0(\phi(0))^T
f(x_t(0),x_t(-\tau_1),\cdots,x_t(-\tau_K))+\sum_{i=1}^K
p_i(\phi(0))-p_2(\phi(-\tau_i))\\
&\le -\alpha \norm{x_t(0)}^2
\end{align*}
Therefore the derivative is negative definite by the conditions of
the theorem. Therefore we have that the system defined by $f$ is
globally asymptotically stable for any finite values of $\tau_i$.
\end{proof}

The computational burden associated with the conditions expressed in
Theorem~\ref{thm:nonlinear_multiple_ind} is significantly less than
that associated with proving delay-independent stability using
Theorem~\ref{thm:nonlinear_multiple}. This reduction occurs because
in this case, we have deliberately chosen a Lyapunov functional
structure where the delays $\tau_i$ do not explicitly appear in the
derivative. Naturally, the use of a more specialized Lyapunov
functional structure may result in an increased level of
conservativity. For this reason, when testing delay-independent
stability, it is recommended that
Theorem~\ref{thm:nonlinear_multiple_ind} be used as a first
approximation and that Theorem~\ref{thm:nonlinear_multiple} be used
if additional accuracy is required.

\subsection{Numerical Examples}\label{sec:example}
The contribution of this chapter is to show that Lyapunov functions
to prove stability of arbitrary non-linear polynomial systems may be
computed algorithmically. To illustrate this, we compare our results
with those derived through less algorithmic analysis.

\textbf{Example 1: Reflection Dynamics}

To begin, we compare our results with a solution given by
Hale~\cite{hale_1993} and attributed to LaSalle for the following
dynamics.
\begin{equation}
\dot{x}(t)=a x^3(t)+ b x^3(t-\tau)
\end{equation}
The following Lyapunov function is used by Hale to show that the
system is stable for any $a<0$, $|b|<|a|$.
\begin{equation}
  \label{eqn:Hale_functional}
  V(x_t)=-\frac{x_t^4(0)}{2 a}
  + \int_{-\tau}^{0}
  x_t^6(\theta)\, d \theta
\end{equation}

Our own computation of various points give similar results. For
example, for $a=-1$, $b=.9$, we can prove asymptotic stability using
the following function
\[
V(x_t)=-3.92 x_t^2(0)+1.91 x_t^4(0) + \int_{-\tau}^{0}
5.35x_t^4(\theta)+ 3.72 x_t^6(\theta) \,d \theta
\]

\begin{figure}
\centerline{\includegraphics[width=0.3\textwidth]{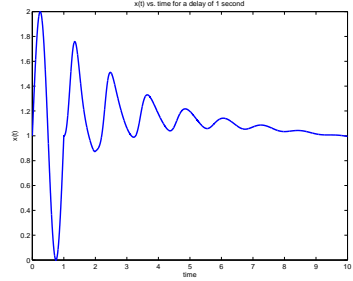}}
\caption{Example of a trajectory of $x(t)$ for $a=-1$, $b=.9$}
\end{figure}
\begin{figure}
\centerline{\includegraphics[width=0.3\textwidth]{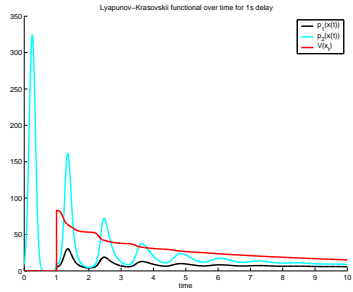}}
\caption{Example of a trajectory of $V(x_t)$ for $a=-1$, $b=.9$}
\end{figure}

%\begin{figure}[h]
%\centerline{\includegraphics[width=0.5\textwidth]{numexample1_1.eps}}
%\caption{Trajectory of $x(t)$}
%\end{figure}
%
%\begin{figure}[h]
%\centerline{\includegraphics[width=0.5\textwidth]{numexample1_2.eps}}
%\caption{Trajectory of $V(x_t)$}
%\end{figure}

The purpose here is not to construct this particular Lyapunov
function, but rather to demonstrate that such functions can be
numerically computed without the insight needed to manually derive
the results. To illustrate this expanded flexibility, suppose we
introduce a cross-term, such as might arise due to interference,
then we have
\begin{equation}
\dot{x}(t)=a x^3(t)+ c (x(t) x(t-\tau))^2 +b x^3(t-\tau)
\end{equation}
We find the following Lyapunov function in the case when $a=-1$,
$b=0.5$, $c=0.2$.
\begin{align*}
V(x_t)=&4.233 x_t^2(0)-0.2147 x_t^3(0) \\
&+ 1.856 x_t^4(0)
 \int_{-\tau}^{0}
4.107 x_t^4(\theta)+ 0.1525 x_t^5(\theta) + 3.454 x_t^6(\theta) \,d
\theta
\end{align*}

\textbf{Example 2: Epidemiology}

We now consider the following epidemiological model. Consider a
human population subject to non-lethal infection by a virus. Assume
the disease has incubation period $\tau$. Cooke~\cite{cooke_1979}
models the percentage of infected population, $y(t)$ using the
following dynamics.
\[\dot{y}(t)=-a y(t)+b y(t-\tau)\br(1-y(t)\bl)
\]
Where
\begin{itemize}
\item $a$ is the recovery rate for the infected population.
\item $b$ is the rate of infection for those exposed to the virus.
\end{itemize}

The model is nonlinear and equilibria exist at $y^*=0$ and
$y^*=(b-a)/b$. Cooke used the following Lyapunov functional to prove
delay-independent stability of the $0$ equilibrium for $a>b>0$.

\[V(\phi)=\frac{1}{2}\phi(0)^2+\frac{1}{2}\int_{-\tau}^0 a
\phi(\theta)^2 d\theta
\]

Using semidefinite programming, we were also able to prove
delay-independent stability for $a>b>0$ using the following
functional.
\[V(\phi)=1.75 \phi(0)^2+\int_{-\tau}^0(1.47 a + .28
b)\phi(\theta)^2 d \theta
\]

The point of these examples is obviously not to trivialize the work
of Cooke and Hale, but rather to show the broad applicability and
utility of the algorithms presented in this thesis.

\section{Conclusion}\label{sec:conclusion}

This chapter provides a generalization of the stability analysis
results for linear time-delay systems presented in
Chapter~\ref{chp:LinearCase} to the case of nonlinear dynamics. To
this end, we have proposed a sequence of sufficient conditions for
both delay-dependent and delay independent stability of nonlinear
systems with multiple delays. The numerical examples presented
illustrate that many previous results on stability analysis can be
derived using our algorithmic framework. Unfortunately, unlike in
the linear stability case, we are not optimistic that the results of
this chapter will lead to algorithms for the construction of
stabilizing controllers for any broad class of nonlinear time-delay
systems. This is because the results of this chapter to not fit as
neatly into the operator-theoretic framework of the linear case.

\chapter{Conclusion}
The purpose of this thesis has been to find ways of addressing the
question of stability of complex systems. In particular, we have
sought ways of using results from functional analysis, real
algebraic geometry and computer science to provide new approaches to
answering the question of stability of nonlinear differential
equations which contain delay. The results can be divided into two
different approaches.

\section{Generalized Passivity}
In Chapter~\ref{chp:IQCresult}, we showed how one can restructure
the stability question of a nonlinear, discontinuous model of
internet congestion control which contains delay. We showed how the
question of stability of the model could be recast as the question
of passivity of each of two interconnected operators. We then showed
that by choosing to decompose the system into simpler subsystems,
one of which is nonlinear and discontinuous but delay free and the
other which is linear, but contains delay, one can use more
specialized tools on each subsystem. By using more specialized
tools, we were able to obtain less conservative results for each
subsystem. We then combined these improved bounds using the
generalized passivity framework to obtain less conservative results
for the combined system.

\paragraph{Practical Impact} There are a number of factors to
consider when evaluating the practical impact of this work on
stability of Internet congestion control. As mentioned earlier, it
has been shown that current internet protocols will destabilize as
transmission rates and link capacities continue to increase. Finding
a stable alternative to current protocols which converges to a
global optimum is therefore a high priority for researchers. In
addition, because of the broad geographic reach of the internet, one
should expect to find that transmission delay will be a major
consideration for any future internet design. In particular, a
simple calculation reveals that the speed of light alone can account
for a round-trip delay of .1 seconds or more. It is well-known to
the practicing engineer that such delays will have a destabilizing
effect on most control systems at high values of gain. An
alternative version of internet congestion control which
specifically attempts to account for these delays is the FAST
protocol which has presented in this paper. It has been shown that
the FAST protocol is stable about the equilibrium for arbitrary
topology and delay.

The work in this paper has shown that the FAST protocols are also
stable globally for a more limited network topology. This work is
important in that it shows, for the first time, that the protocols
developed for the linearized system also converge in the nonlinear
setting without additional constraints on the parameters. From a
practical standpoint, the improvement of 57\% in the gain parameter
$\alpha$ over the best available previous bound should allow a
substantial increase in the rate of convergence to the global
optimal. Since the internet is generally a dynamic environment,
quick convergence of the system to equilibrium is of particular
importance. In addition, the increase in $\alpha$ also results in a
$37\%$ decrease in the equilibrium value of the price. Since price
tracks queue length in most implementations, this result gives a
$37\%$ decrease in size of the equilibrium queue. This results in a
decrease of the storage space needed at the queues and therefore
cost savings for providers of routers. However, we note that for
this increase in the gain parameter to be implemented throughout the
internet, one would like to see an extension of our results to the
general network topology.

\paragraph{Research Directions} As we mentioned in
Chapter~\ref{chp:IQCresult}, the proof we have given in the limited
topology framework is not easily extended to the more general case.
The reason for this is that the dynamics at the link are nonlinear
due to the discontinuity due to the positive projection. The
approach in this thesis was to separate the nonlinearity from the
delay. Ideally in the multiple-user/multiple-link framework, we
would also separate the nonlinearity and delay so that we would be
left with the feedback interconnection of multiple nonlinear
subsystems and multiple linear, but delayed subsystems. However, the
presence of nonlinearity at the link makes this approach
problematic. To date, we have not overcome this difficulty.

Another extension of this work is the use of our methods to consider
alternative versions of internet congestion control. Almost any
protocol should be amenable to this type of analysis provided that
the nonlinearity satisfies a sector bound. In addition, if one
analyzes a linear source control with a nonlinear, discontinuous
link control, then the problems discussed earlier are lessened and
extensions to the general network topology are more likely to be
found.

\section{Algorithms for Stability Analysis}
In Chapter~\ref{chp:LinearCase}, we showed how the question of
stability for linear time-delay systems could be posed as the search
for a positive complete quadratic Lyapunov functional with negative
derivative. We then explained how the search for such a functional
could be expressed as a convex feasibility problem over certain
convex cones of positive linear operators. We then demonstrated that
if one considers the subset of such operators which are defined by
polynomials of bounded degree, then the search for such operators
can be expressed using semidefinite programming constraints. This
allowed us to express the question of stability of linear time-delay
systems as a semidefinite programming problem which is amenable to
solution using recently developed interior-point methods. In
addition, we were able to generalize this approach to linear systems
with parametric uncertainty. In Chapter~\ref{chp:nonlinearcase}, we
were able to further generalize this approach to nonlinear
time-delay systems through the use of non-quadratic Lyapunov
functionals. We also dealt with the particular case of
delay-independent stability using a particular form of Lyapunov
functional for which the resulting semidefinite program was of
significantly decreased complexity.

\paragraph{Practical Impact}
The practical impact of Chapters~\ref{chp:LinearCase}
and~\ref{chp:nonlinearcase} is to show how stability proofs for a
broad class of time-delay systems can be computed quickly and
efficiently using existing numerical algorithms. In particular, by
using Matlab packages such as SOSTools, the equality constraints
implied in the stability theorems can be implemented quickly using a
minimal amount of coding. The semidefinite program can then be
solved efficiently using interior-point algorithms such as that in
SeDuMi and its Matlab interface.

In addition, for the linear case at least, the size of the
semidefinite programs associated with the various stability theorems
scales well in the size of the problem. In particular, for a
$n$-dimensional single-delay system using polynomial functions of
degree $2d$, the number of variables scales as $(n \cdot d)^2$. For
a $n$-dimensional system with $K$-delays using polynomial functions
of degree $2d$, the number of variables scales as $(n\cdot d \cdot
K)^2$. When parametric uncertainty or nonlinearity are present,
however, the scaling becomes much worse. In particular, for the case
of a linear $n$-dimensional system with $K$-delays using polynomial
functions of degree $2d$, and $m$ uncertain parameters $\alpha_i$
using polynomials of degree $2d_2$ in $\alpha_i$, the number of
variables scales as $(n\cdot d \cdot K \cdot {m+d_2 \choose
d_2})^2$.

In engineering practice, of course, there are many ways that one can
deal with the presence of delays. For example, in designing
controllers, it is standard to use a Pad\'{e} approximation to
replace the transcendental term in the characteristic equation with
a rational polynomial. This rational transfer function can then be
controlled using standard techniques. In general, however, such
controlled systems are not guaranteed to be stable unless the system
also satisfies some additional restrictions including bounds on size
of the delay and bandwidth restrictions on the input. For analysis,
the advantage of our method over direct use of the Pad\'{e}
approximation lies in the fact that our method guarantees stability
regardless of dynamic properties of the system. Moreover, if one
considers the question of control, then our approach can be seen
more as a complement to use of the Pad\'{e} approximation, since the
one method can be used to efficiently design controllers while the
other can be used to verify stability of such controllers.

In addition to direct use of the Pad\'{e} approximation, there exist
methods which attempt to use a rational approximation for the delay
and combine this with some parametric uncertainty to give conditions
which are guaranteed to prove stablity. Furthermore, there exist
bounds~\cite{zhang_2003} on how well this Pad\'{e} technique
approximates the original system as a function of the degree of the
approximation. Unfortunately, the conditions associated with this
approach require the user to solve certain parameter-dependent LMI
problems. The advantage of our approach in this case is to give an
explicit numerical procedure for computing stability which involves
a minimum amount of effort by the user.

Finally, we conclude this section by mentioning perhaps the most
important advantage of our algorithm over previous existing methods.
This advantage lies in the ability of the algorithm to adapt to
numerous forms of generalization beyond the linear fixed-delay
result. In particular, our methods can be used to construct a broad
class of Lyapunov functionals for analysis of nonlinear and
uncertain time-delay systems. For such systems there exist very few
analysis tools, making our algorithm stand out simply due to the
lack of competition.

\paragraph{Research Directions}
The work in this thesis on computation of Lyapunov functionals for
time-delay systems provides a number of exciting directions for
further research. If one views the results in the linear case as an
attempt to generalize prior work on finite-dimensional linear
operator theory to infinite dimensions, then perhaps the most
obvious next step is to derive conditions for the synthesis of
stabilizing controllers. In particular, we have recent results that
show how to construct the inverse operators to the combined
multiplier and integral operators which define the complete
quadratic functional. If we rely on the analogy of linear
finite-dimensional systems, then this result should allow us to
directly synthesize stabilizing controllers. Unfortunately for
linear time-delay systems, however, the dual to the linear
time-delay system is not a linear time-delay system. Thus the most
obvious route for this research would seem closed to us.
Fortunately, however, the work by Delfour and
Mitter~\cite{delfour_72a}, among others, has shown us that instead
of analyzing stability of the dual system, we may consider stability
of the adjoint problem which is itself a linear time-delay system.
Furthermore, this adjoint system has stability properties that are
directly related to that of the original system. However, further
work in this direction is necessary before making any stronger
assertions. Apart from stabilization, other possible generalizations
of finite-dimensional linear theory include optimal control,
filtering and an operator-theoretic version of the KYP lemma for
transcendental polynomials.

Another direction for this research is to consider more general
classes of infinite-dimensional systems. In particular, some classes
of partial differential equations are also amenable to analysis
using Lyapunov functional methods. Thus it may be possible in some
cases to derive stability proofs for these types of systems. In
particular, the question of stability of the Navier-Stokes equation
is a subject of considerable interest. While we are far from
claiming that this question is a straightforward generalization of
our methods, we do not discount the possibility that it may
eventually fall within the scope of algorithmic analysis. Generally,
however, in order to study stability of a class of differential
equation, one must be able to structure the derivative of the
Lyapunov functional in such a way that its negativity can be
verified using semidefinite programming. In the case of linear
time-delay systems, the derivative of the Lyapunov functional could
be structured in a way which was similar to that of the functional
itself. For this reason, conditions for negativity of the derivative
were a simple generalization of the conditions for positivity of the
functional and therefore expressing the conditions for stability as
a semidefinite program was relatively easy. In general the same may
not be true for larger classes of system. We do have some reason to
be optimistic, however, since the case of linear time-delay systems
is actually structured in a rather complicated manner and showing
how to restructure the derivative of the functional was, in fact,
not trivial. Less structured types of systems may, in fact, be
easier to analyze. However, without further research we cannot
provide a more definitive answer to this line of inquiry.

Apart from further developments in theory, perhaps the most
important research direction is application of the algorithms
presented in this thesis to practical problems. Systems with
nonlinearity and delay arise in a number of different disciplines.
In particular, modeling of biological phenomena has show that many
aspects of the natural world are dominated by nonlinear and
infinite-dimensional behavior. Much recent research effort has been
put into modeling these type of system and as more models arise, one
can project the need for more efficient ways of analyzing them. Thus
we are considerably interested to see if the techniques proposed in
this thesis can find use in applications as diverse as determining
effectiveness of cancer therapy or the ability of the immune system
to compensate for the constantly changing genetic makeup of viruses
such as HIV.

\section{Final Word}
Recent developments have made the last 10 years an increasingly
exciting time to be studying control theory. In particular, the
advent of interior-point algorithms for semidefinite programming has
made a fundamental shift in how research is conducted. Instead of
asking whether we can solve a problem numerically, it has become
commonplace to instead ask whether the question can be reformulated
as a convex optimization problem. Thus the very definition of a
solution has become linked with the general availability of
efficient numerical algorithms. This interaction between control
theory and computer science has led to an ability to quickly analyze
systems which a generation ago would perhaps have been considered
prohibitively complex. Furthermore, as the definition of solution
continues to evolve, new opportunities will arise for adventurous
researchers to find ways to address problems which common sense
might now classify as being out of reach. We note that such
enthusiasm must be tempered with caution, however, lest mathematical
rigor be sacrificed in the process. As to our own work, we note that
none of the results presented would have been possible without the
solid theoretical underpinning of computational complexity, convex
optimization, real algebraic geometry, functional analysis, operator
theory, Lyapunov theory, and the input-output framework.
% And so we
%conclude with the hope that our results too will withstand rigorous
%scrutiny and may provide a solid basis for the possibility of future
%work of the type which we have so ambitiously listed in previous
%parts of this section.

% and the end material

\appendix

\chapter{Appendix to
Chapter~\ref{chp:FDEtheory}}\label{app:FDEtheory}

\paragraph{Additional Notation:}

Denote the complete metric space
\[V_{\alpha}:=\{x\in \mathcal{C}[-\tau,\infty): x(\theta)=\phi(\theta) \text{ for } \theta \in [-\tau,0],\;
\sup_{t \ge -\tau}\norm{x(t)}e^{-\alpha t} <\infty \},
\]
which has the following metric.
\[d(x,y)_{V_\alpha}=\sup_{t \ge -\tau}\norm{x(t)-y(t)}e^{-\alpha t}
\]
Define $\mathcal{C}_\tau$ to be the space of continuous functions
defined on $[-\tau,0]$ with norm
\[\norm{x}_{\mathcal{C}_\tau}=\sup_{t\in [-\tau,0]}\norm{x(t)}_2.
\]
%Define the operator $x_t:\mathcal{C}[-\tau, \infty)\times \R^+
%\rightarrow \mathcal{C_\tau}$ by
%\[x_t(\theta)=x(t+\theta) \qquad \text{for } \theta \in [-\tau,0].
%\]
For a given functional $f:\mathcal{C}_\tau \times \R^+ \rightarrow
\R^n$, let $\Lambda$ be defined by
\[(\Lambda x)(t)=\phi(0)+\int_{0}^t f(x_s,s)ds \qquad \text{for all }
t \ge 0
\]
and $(\Lambda x)(t)=\phi(t)$ for $t \in [-\tau,0)$.

\begin{defn}We say that a functional $f:\mathcal{C}_\tau \times \R^+ \rightarrow
\R^n$ satisfies \eemph{Assumptions 1} if
\begin{itemize}\item There exists a $K_1>0$ such that
\[\norm{f(x,t)-f(y,t)}_2 \le K_1 \norm{x-y}_{\mathcal{C}_\tau} \quad
\text{for all }x,y \in \mathcal{C}_\tau, t \ge 0
\]
\item There exists a $K_2>0$ such that
\[\norm{f(0,t)}\le K_2 \quad \text{for all }t \ge 0
\]
\item $f(x,t)$ is jointly continuous in $x$ and $t$. i.e. for every
$(x,t)$ and $\epsilon>0$, there exists a $\eta>0$ such that
\[\norm{x-y}_{\mathcal{C}_\tau}+\norm{t-s}\le \eta \; \Rightarrow \;
\norm{f(x,t)-f(y,s)}\le \epsilon.
\]
\end{itemize}
\end{defn}

\begin{defn}
Given a functional $f:\mathcal{C}_\tau \times \R^+ \rightarrow
\R^n$, we say that a function $x\in \mathcal{C}[-\tau,\infty)$ is a
\eemph{solution to Problem 1} with initial condition $\phi \in
\mathcal{C}_\tau$ if $x$ is differentiable for $t\ge0$,
$x(t)=\phi(t)$ for $t\in [-\tau,0]$ and
\[\dot{x}(t)=f(x_t,t)\qquad \text{for all }t \ge 0
\]
\end{defn}

\paragraph{Results:}

\begin{lem} \label{lem:continuity} If $f:\mathcal{C}_\tau \times \R^+ \rightarrow
\R^n$ satisfies Assumptions 1, and $x$ is a continuous function,
then $f(x_t,t)$ is a continuous function in $t$.
\end{lem}
\begin{proof}
If $f:\mathcal{C}_\tau \times \R^+ \rightarrow \R^n$ satisfies
Assumptions 1, then for any $t_1$, $\epsilon>0$, there exists a
$\eta_1>0$ such that
\[\norm{x_{t_1}-x_{t_2}}_{\mathcal{C}_\tau}+\norm{t_1-t_2}\le \eta_1 \; \Rightarrow \;
\norm{f(x_{t_1},t_1)-f(x_{t_2},t_2)}\le \epsilon.
\]
Now if $x(t)$ is a continuous function of $t$, then for any $t_1 \ge
0$, $x(t)$ is uniformly continuous on the interval
$I_{t_1}:=[\max\{-\tau,t_1-\tau-.1\},t_1+.1]$. Hence, for any
$\epsilon_2>0$, there exists a $\eta_2
>0$ such that for $s_1,s_2 \in I_{t_1}$, we
have that
\[|s_1-s_2| \le \eta_2 \; \Rightarrow \;\norm{x(s_1)-x(s_2)}_2 \le
\epsilon_2.
\]
Now for a given $t_1 \ge 0$ and $\epsilon>0$, let $\eta_1$ be
defined as above. Let $\epsilon_2=\eta_1/2$ define $\eta_2$. Then if
$\eta=\min\{.1,\eta_2,\eta_1/2\}$, we have that $|t_1-t_2|\le \eta$
implies that
\begin{align*}
\norm{x_{t_1}-x_{t_2}}_{\mathcal{C}_\tau}
&=\sup_{\theta \in [-\tau,0]}\norm{x_{t_1}(\theta)-x_{t_2}(\theta)}_2\\
&=\sup_{\theta \in [-\tau,0]}\norm{x(t_1+\theta)-x(t_2+\theta)}_2
\le \epsilon_2 =\eta_1/2.
\end{align*}
This holds since $\norm{(t_1+\theta)-(t_2+\theta)}=\norm{t_1-t_2}
\le \eta_2$ means $\norm{x(t_1+\theta) - x(t_2+\theta)}_2 \le
\epsilon_2$ since $t_1+\theta \in I_{t_1}$ for $\theta \in
[-\tau,0]$ and $t_2+\theta \in I_{t_1}$ for $\norm{t_1-t_2}<.1$,
$\theta \in [-\tau,0]$ and $t_2 \ge 0$. Therefore we have that
\[\norm{x_{t_1}-x_{t_2}}_{\mathcal{C}_\tau}+\norm{t_1-t_2}\le
\eta_1/2+\eta_1/2=\eta_1.
\]
This implies that
\[\norm{f(x_{t_1},t_1)-f(x_{t_2},t_2)}\le \epsilon.
\]
Thus we have continuity of $f(x_t,t)$.
\end{proof}
\begin{lem}\label{lem:map}
Suppose $f:\mathcal{C}_\tau \times \R^+ \rightarrow \R^n$ satisfies
Assumptions 1, then $\Lambda: V_\alpha \rightarrow V_\alpha$.
\end{lem}
\begin{proof}
Suppose $x \in V_\alpha$. Let $y=\Lambda x$. Then $y(t)=\phi(t)$ for
$t \in [-\tau,0]$, $y$ is continuous, and the following holds for $t
\ge 0$.
\begin{align*}
\norm{y(t)}_2&=\norm{\phi(0)+\int_{0}^t f(x_s,s)ds}_2\\
&\le\norm{\phi(0)}_2+\norm{\int_{0}^t f(x_s,s)ds}_2\\
&\le\norm{\phi(0)}_2+\int_{0}^t \norm{f(x_s,s)}_2 ds\\
&=\norm{\phi(0)}_2+\int_{0}^t \norm{f(x_s,s)-f(0,s)+f(0,s)}_2 ds\\
&\le \norm{\phi(0)}_2+\int_{0}^t \norm{f(x_s,s)-f(0,s)}_2+\norm{f(0,s)}_2 ds\\
&\le \norm{\phi(0)}_2+\int_{0}^t K_1\norm{x_s}_{\mathcal{C}_\tau}+K_2 ds\\
&= \norm{\phi(0)}_2+K_1 \int_{0}^t \sup_{\theta \in [-\tau,0]}\norm{x_s(\theta)}_2 ds + K_2 t\\
&= \norm{\phi(0)}_2+K_1 \int_{0}^t \sup_{\theta \in [-\tau,0]}\norm{x(s+\theta)}_2 ds + K_2 t\\
\end{align*}
Now, since $x\in V_\alpha$, there exists some $c>0$ such that
$\norm{x(t)}_2\le c e^{\alpha t}$ for all $t \ge -\tau$. Therefore,
$\sup_{\theta \in [-\tau,0]}\norm{x(t+\theta)}\le c e^{\alpha t}$
for all $t \ge 0$. Therefore, we have the following for $t \ge 0$.
\begin{align*}
\norm{y(t)}_2&\le \norm{\phi(0)}_2+K_1 \int_{0}^t \sup_{\theta \in [-\tau,0]}\norm{x(s+\theta)}_2 ds + K_2 t\\
&\le \norm{\phi(0)}_2+K_1 \int_{0}^t c e^{\alpha s} ds + K_2 t\\
&= \norm{\phi(0)}_2+\frac{c K_1}{\alpha} (e^{\alpha t}-1) + K_2 t
\end{align*}
Thus
\begin{align*}
\sup_{t \ge 0}\norm{y(t)}_2 e^{-\alpha t}&\le \sup_{t \ge 0}
\left(\norm{\phi(0)}_2e^{-\alpha t} +\frac{c K_1}{\alpha}
(1-e^{-\alpha t}) + K_2 t e^{-\alpha t}\right) \\
&\le \norm{\phi(0)}_2 +\frac{c K_1}{\alpha} + \frac{e
K_2}{\alpha}=c_1.
\end{align*}
Now let $c_2=\sup_{\theta \in [-\tau,0]}\norm{\phi(\theta)}$ and
$\tilde{c}=\max\{c_1,c_2\}$. Then $\sup_{t \ge -\tau}\norm{y(t)}_2
e^{-\alpha t} \le \tilde{c}$. Thus $y \in V_\alpha$.
\end{proof}

\begin{lem}\label{lem:contraction}
Suppose $f:\mathcal{C}_\tau \times \R^+ \rightarrow \R^n$ satisfies
Assumptions 1,  and $K_1 < \alpha$. Then $\Lambda$ is a contraction
on $V_\alpha$.
\end{lem}

\begin{proof}
Suppose $x,y\in V_\alpha$. Then we have the following for $t \ge 0$.
\begin{align*}
\norm{\Lambda x(t)-\Lambda y(t)}_2 &= \norm{\int_{0}^t
f(x_s,s)-f(y_s,s)ds}_2\\
&\le \int_{0}^t \norm{f(x_s,s)-f(y_s,s)}_2 ds\\
&\le \int_{0}^t K_1 \norm{x_s-y_s}_{\mathcal{C}_\tau} ds\\
&= K_1 \int_{0}^t \sup_{\theta \in [-\tau,0]}\norm{x_s(\theta)-y_s(\theta)}_2 ds\\
&= K_1 \int_{0}^t \sup_{\theta \in [-\tau,0]}\norm{x(s+\theta)-y(s+\theta)}_2 ds\\
&= K_1 \int_{0}^t \sup_{\theta \in [s-\tau,s]}\norm{x(\theta)-y(\theta)}_2 ds\\
\end{align*}
Now recall that $d(x,y)_{V_\alpha}=\sup_{t \ge
-\tau}\norm{x(t)-y(t)}e^{-\alpha t}$. Therefore, $\norm{x(t)-y(t)}_2
\le d(x,y)_{V_\alpha} e^{\alpha t}$ for all $t \ge -\tau$. Thus
$\sup_{\theta \in [t-\tau,t]}\norm{x(\theta)-y(\theta)}_2 \le
d(x,y)e^{\alpha t}$ for all $t \ge 0$. Thus we have the following
for all $t\ge 0$.
\begin{align*}
\norm{\Lambda x(t)-\Lambda y(t)}_2 &\le K_1 \int_{0}^t \sup_{\theta \in [s-\tau,s]}\norm{x(\theta)-y(\theta)}_2 ds\\
&\le K_1 \int_{0}^t d(x,y)_{V_\alpha} e^{\alpha s} ds\\
&= \frac{K_1}{\alpha}d(x,y)_{V_\alpha} (e^{\alpha t}-1)
\end{align*}
Therefore
\begin{align*}
d(\Lambda x,\Lambda y)_{V_\alpha}&=\sup_{t\ge -\tau}\norm{\Lambda
x(t)-\Lambda y (t)}_2 e^{-\alpha t}\\
&=\sup_{t\ge 0}\norm{\Lambda x(t)-\Lambda y (t)}_2 e^{-\alpha
t}\\
&\le \sup_{t\ge 0}\frac{K_1 }{\alpha}d(x,y)_{V_\alpha} (1-e^{-\alpha t}) \\
&\le \frac{K_1}{\alpha} d(x,y)_{V_\alpha} < d(x,y)_{V_\alpha}\\
\end{align*}
Therefore $\Lambda$ is a contraction on $V_\alpha$.
\end{proof}

\begin{lem}[Existence]\label{lem:existence}
Suppose $f$ satisfies Assumptions 1, $x \in V_\alpha$, and
$x=\Lambda x$. Then $x$ is a solution of Problem 1. i.e. $x$ is
differentiable , $x(t)=\phi(t)$ for $t \in [-\tau,0]$ and
\[\dot{x}(t)=f(x_t,t) \qquad \text{for all }t \ge 0.
\]
\end{lem}

\begin{proof}
Since $x \in V_\alpha$, $x$ is continuous and $x(t)=\phi(t)$ for $t
\in [-\tau,0]$. Furthermore, since $x=\Lambda x$, we have that
\[x(t)=\phi(0)+\int_0^t f(x_s,s)ds.
\]
Since $x$ is continuous and $f$ satisfies Assumptions 1, by
Lemma~\ref{lem:continuity} we have that $f(x_s,s)$ is a continuous
and integrable function on $[0,t]$. Now, since $x(t)-x(0)=\int_0^t
f(x_s,s)ds$, we have by the fundamental theorem of Calculus(v1) that
$x(t)-x(0)$ is differentiable and
\[\dot{x}(t)=f(x_t,t).
\]
\end{proof}

\begin{lem}[Uniqueness]\label{lem:uniqueness}
Suppose $f$ satisfies Assumptions 1, then any solution to Problem 1
is unique.
\end{lem}
\begin{proof}
Suppose $x$ and $y$ are solutions to Problem 1. Define $p$ as
follows for $t\ge 0$.
\begin{align*}
p(t)&=\sup_{\theta \in [t-\tau,t]}\norm{x(\theta)-y(\theta)}_2\\
&=\sup_{\theta \in [t-\tau,t], \theta \ge 0}\norm{x(\theta)-y(\theta)}_2\\
&=\sup_{\theta \in
[t-\tau,t], \theta \ge 0}\norm{\int_0^\theta f(x_s,s)-f(y_s,s)ds}_2\\
& \le \sup_{\theta \in
[t-\tau,t], \theta \ge 0}\int_0^\theta \norm{f(x_s,s)-f(y_s,s)}_2 ds\\
& \le \sup_{\theta \in
[t-\tau,t], \theta \ge 0} K_1 \int_0^\theta \norm{x_s-y_s}_{\mathcal{C}_\tau} ds\\
& = \sup_{\theta \in [t-\tau,t], \theta \ge 0} K_1 \int_0^\theta
\sup_{\omega \in [s-\tau,s]}\norm{x(\omega)-y(\omega)}_2 ds\\
& = \sup_{\theta \in [t-\tau,t], \theta \ge 0} K_1 \int_0^\theta
p(s) ds\\
\end{align*}
Now let $q(t)=\int_0^t p(s) ds$. Then $q(0)=0$ and
$\dot{q}(t)=p(t)\ge 0$ for all $t\ge$. Therefore, $q$ is
non-negative and monotonically increasing. Thus we have the
following.
\begin{align*}\dot{q}(t)&\le \sup_{\theta \in [t-\tau,t],
\theta \ge 0} K_1 \int_0^\theta p(s) ds\\
&\le \sup_{\theta \in [t-\tau,t], \theta \ge 0} K_1 q(\theta)=K_1
q(t)
\end{align*}
Therefore, we have that $\dot{q}(t)-K_1 q(t)\le 0$ for all $t \ge
0$. Which implies the following for $g(t)=q(t)e^{-K_1 t}$, $t\ge0$.
\[\dot{g}(t)=\frac{d}{dt}\left(q(t)e^{-K_1 t} \right)=e^{-K_1 t}\left(\dot{q}(t)-K_1
q(t)\right)\le 0
\]
Since $g(0)=q(0)=0$ and $\dot{g}(t)\le 0$ for all $t\ge0$, we have
that $g(t)\le 0$ for all $t\ge 0$. However, $g(t)\le 0$ implies
$q(t)\le 0$, and so we have that $q(t)=0$ for all $t\ge 0$.
Therefore, since $q(t)=\dot{p}(t)=0$ for all $t\ge0$ and $p(0)=0$,
we have that $p(t)=0$ for all $t\ge 0$. This implies that
$x(t)=y(t)$ for all $t\ge0$, which implies that $x=y$. Thus any
solution to problem 1 is unique.
\end{proof}
\begin{thm}
Suppose a functional $f:\mathcal{C}_\tau \times \R^+ \rightarrow
\R^n$ satisfies the following.
\begin{itemize}\item There exists a $K_1>0$ such that
\[\norm{f(x,t)-f(y,t)}_2 \le K_1 \norm{x-y}_{\mathcal{C}_\tau} \quad
\text{for all }x,y \in \mathcal{C}_\tau, t \ge 0.
\]
\item There exists a $K_2>0$ such that
\[\norm{f(0,t)}\le K_2 \quad \text{for all }t \ge 0.
\]
\item $f(x,t)$ is jointly continuous in $x$ and $t$. i.e. for every
$(x,t)$ and $\epsilon>0$, there exists a $\eta>0$ such that
\[\norm{x-y}_{\mathcal{C}_\tau}+\norm{t-s}\le \eta \; \Rightarrow \;
\norm{f(x,t)-f(y,s)}\le \epsilon.
\]
\end{itemize}
Then for any $\phi \in \mathcal{C}_\tau$, there exists a unique $x
\in \mathcal{C}[-\tau,\infty)$ such that $x$ is differentiable for
$t\ge0$, $x(t)=\phi(t)$ for $t\in [-\tau,0]$ and
\[\dot{x}(t)=f(x_t,t)\qquad \text{for all }t \ge 0.
\]
\end{thm}
\begin{proof}
We have already shown by Lemmas~\ref{lem:map}
and~\ref{lem:contraction} that $\Lambda$, as defined by $f$, is a
contraction on $V_\alpha$ for any $\alpha>K$. Therefore, there
exists some $\alpha>K$, $x \in V_\alpha$ such that $x=\Lambda x$.
Therefore, by Lemma~\ref{lem:existence}, we have that $x$ is a
solution. By Lemma~\ref{lem:uniqueness}, we have that this solution
is unique.
\end{proof}

\begin{lem}
Let
\[f(x,t):=\sum_{j=1}^N A_j x(t-\tau_j)+\int_{-\tau}^0 A(\theta)x(t-\theta)d\theta,\]
where $A(\theta)$ is bounded on $[-\tau,0]$. Then there exists some
$K>0$ such that
\[\norm{f(x,t)-f(y,t)}_2 \le K \norm{x-y}_{\mathcal{C}_\tau}.
\]
\end{lem}
\begin{proof}
For $A \in \R^{n \times m}$, let $\bar{\sigma}(A)$ denote the
maximum singular value of $A$. Then we have the following.
\begin{align*}
&\norm{f(x,t)-f(y,t)}_2\\
&=\norm{\sum_{j=1}^N A_j
(x(t-\tau_j)-y(t-\tau_j))+\int_{-\tau}^0
A(\theta)(x(t-\theta)-y(t-\theta))d\theta}_2\\
&\le \sum_{j=1}^N \norm{A_j (x(t-\tau_j)-y(t-\tau_j))}_2+
\int_{-\tau}^0 \norm{
A(\theta)(x(t-\theta)-y(t-\theta))}_2d\theta\\
&\le \sum_{j=1}^N \bar{\sigma}(A_j)\norm{ x(t-\tau_j)-y(t-\tau_j)}_2
+ \int_{-\tau}^0 \bar{\sigma}(A(\theta))\norm{
x(t-\theta)-y(t-\theta)}_2d\theta\\
&\le \sum_{j=1}^N \bar{\sigma}(A_j)\norm{ x-y}_{\mathcal{C}_\tau} +
\int_{-\tau}^0 \bar{\sigma}(A(\theta))d\theta \norm{
x-y}_{\mathcal{C}_\tau}\\
&= \left(\sum_{j=1}^N \bar{\sigma}(A_j) + \int_{-\tau}^0
\bar{\sigma}(A(\theta))d\theta\right) \norm{ x-y}_{\mathcal{C}_\tau}
\end{align*}
\end{proof}

\include{appendix_ICC}
\chapter{Appendix to
Chapter~\ref{chp:SOStheory}}\label{app:SOStheory}

\begin{lem}$M \in \bar{\Sigma}_s \bigcap \S_{2d}^n[x]$ if and only if
there exists some matrix $Q \in \R^{n_z n \times n_z n}$, where
$n_z= {n+d \choose d}$ such that $Q \ge 0$ and the following holds.
\[M(x)=\left(\bar{Z}^n_d[x]\right)^T Q \bar{Z}^n_d[x]
\]
\end{lem}
\begin{proof}
%We first show that if $M \in \bar{\Sigma}_s \bigcap \S_{2d}^n[x]$,
%then there exists $G_i\in \S_d^n[x]$ for $i=1 \ldots m$ such that
%\[M(x)=\sum_{i=1}^m G_i(x)^2
%\]
Since $M \in \bar{\Sigma}_s$, we know that there exist $G_i\in \R^{n
\times n}[x]$ for $i=1 \ldots m$ such that
\[M(x)=\sum_{i=1}^m G_i(x)^T G_i(x)
\]
Now suppose that there exists some $i'$ and $i_0$, $j_0$ such that
element $[G_{i'}]_{i_0,j_0}$ of $G_{i'}$ is of degree $\hat{d}>d$.
Then, since $M_{i_0,i_0}(x)=\sum_{i=1}^m \sum_k
[G_{i}]_{k,i_0}(x)^2$, we have that $M_{i_0,i_0}$ is of at least
degree $2 \hat{d}$. This is because a finite sum of squares of
$\R[x]$ is of degree $2\hat{d}$, where $\hat{d}$ is the maximum
degree of the squared elements. However, since $M \in \S^n_d$ by
assumption, we have a contradiction. Therefore $G_i \in \R_d^{n
\times n}[x]$. Now, for any element $G_i \in \R_d^{n \times n}[x]$,
there exists a $B \in \R^{n \times n(d_z)}$ where $d_z$ is the
length of $Z_d[x]$ such that $G_i(x)=B \bar{Z}^n_d[x]$. Therefore we
have that
\[M(x)=\sum_{i=1}^m \left(\bar{Z}^n_d[x]\right)^T B^T B \bar{Z}^n_d[x]
= \left(\bar{Z}^n_d[x]\right)^T Q \bar{Z}^n_d[x],
\]
where $Q=\sum_{i=1}^m B^T B \ge 0$.
\end{proof}

\begin{lem} $M \in \bar{\Sigma}_s$ if and only
if $y^T M(x) y \in \Sigma_s$.
\end{lem}
\begin{proof}
$(\Leftarrow)$ Suppose $f(x)=y^T M(x) y \in \Sigma_s$. Then
$f(x)=\sum_{i=1}^m g_i(x,y)^2=\sum_{i=1}^m (b_i(x)^T y)^2$ where the
$b_i$ are some vectors of functions. This follows since the $g_i$
must be homogeneous of degree $1$ in $y$. Therefore, we have the
following.
\[y^T M(x) y=\sum_{i=1}^m y^T b_i(x) b_i(x)^T y
\]
Thus $M(x)=\sum_{i=1}^m b_i(x) b_i(x)^T$. Therefore $M(x)\in
\bar{\Sigma}_s$.

$(\Rightarrow)$ Now Suppose $M \in \bar{\Sigma}_s$. Then for some
$G_i(x) \in \R^{n \times n}[x]$ for $i=1\ldots m$, we have that
$f(x)=y^T M(x) y=\sum_{i=1}^m y^T G_i^T(x)G_i(x)y=\sum_{i=1}^m
g_i(x,y)^T g_i(x,y)$, where $g_i(x,y)=G_i(x)y$. Therefore, $f \in
\Sigma_s$.
\end{proof}

\include{appendix_IQCtheory}
\chapter{Appendix to
Chapter~\ref{chp:IQCresult}}\label{app:IQCresult}

\begin{lem}$0 \le f_1(x)x \le \beta x^2$ where
$\beta=\frac{e^{\frac{\alpha}{\tau} p_0}-1}{p_0}$.
\end{lem}
\begin{proof}
Recall $f_1(y) = \min \bl\{ e^{\frac{\alpha}{\tau} y}-1,
e^{\frac{\alpha}{\tau} p_0}-1\br\}$.Let
\begin{align*}
c_1(y)=\frac{e^{\frac{\alpha}{\tau} y}-1}{y}.
\end{align*}
We show that $\dot{c}_1(y)\ge 0$ for all $y$.
\begin{align*}
\dot{c}_1(y)=\frac{(\frac{\alpha}{\tau}y-1)e^{\frac{\alpha}{\tau}
y}+1}{y^2}.
\end{align*}
$\dot{c}_1(y) \ge 0$ for all $y$ is and only if
$c_2:=\dot{c}_1(y)y^2 \ge 0$ for all $y$.
\begin{align*}
\dot{c}_2(y)=(\frac{\alpha}{\tau}ye^{\frac{\alpha}{\tau} y}.
\end{align*}
Thus $\dot{c}_2(y)<0$ for all $y<0$ and $\dot{c}_2(y)>0$ for all
$y>0$, therefore $c_2$ has a global minimum at $y=0$. Since
$c_2(0)=0$, we conclude that $c_2(y)\ge 0$ for all $y$ and thus
$\dot{c}_1(y)\ge 0$ for all $y$. Since $c_1$ is monotone increasing,
we have the following for all $y \le p_0$.
\begin{align*}
\frac{e^{\frac{\alpha}{\tau} y}-1}{y}\le\frac{e^{\frac{\alpha}{\tau}
p_0}-1}{p_0}.
\end{align*}
Thus for $y \le p_0$.
\begin{align*}
(e^{\frac{\alpha}{\tau} y}-1)y \le \beta y^2.
\end{align*}
Furthermore, $\lim_{y\rightarrow -\infty}c_1(y)=0$ therefore,
$c_1(y)\ge0$ for all $y$, thus
\begin{align*}
(e^{\frac{\alpha}{\tau} y}-1)y \ge 0.
\end{align*}
Therefore, by the definition of $f_1$, we have that $0 \le f_1(x)x
\le \beta x^2$.
\end{proof}

%%%%%%%%%%%%%%%%%%%%%%%%%%%%%%%%%%%%%%%%%%%%%%%%%%%%%%%%%%%%%%%%%%%
\begin{lem}
Let $z=\Delta_z y$ with $y \in W_{2}$, then $\lim_{t \rightarrow
\infty } z(t)=0$.
\end{lem}
\begin{proof}
Let $v=\dot{z}=\Delta y$. Suppose that $T_2>T_1>0$ and let
$H=P_{T_2}-P_{T_1}$. Then
\begin{align*}
\norm{ Hv}_2^2
&= \int_{T_1}^{T_2}\dot{z}(t)^2 dt  \\
& \leq \beta \int_{ T_1}^{T_2} \dot{z}(t)y(t) \,dt -\beta \int_{
T_1}^{T_2} \dot{z}(t)z(t) \,dt
\\
&= \beta \ip{H v}{H y} -\frac{\beta}{2} (z(T_2)^2- z(T_1)^2)
\\
&\le \beta \norm{H v}_2 \norm{H y}_2 -\frac{\beta}{2} \bl(z(T_2)^2-
z(T_1)^2\br)
\end{align*}
Hence
\begin{align*}
z(T_2)^2- z(T_1)^2
& \leq 2 \norm{H v}_2 \norm{H y}_2  -\frac{2}{\beta}\norm{H v}_2^2 \\
&\le 2\norm{H v}_2 \norm{H y}_2
\end{align*}

By Lemma~\ref{lem:z_in_L2}, $v \in L_2$. Since $\norm{v}$ and
$\norm{y}$ exist, we can use the Cauchy criterion and the above
inequality to establish that for any $\delta
> 0$, there exists a $T_{\delta}$ such that $T_2>T_1>T_{\delta}$
implies $(z(T_2)^2- z(T_1)^2)< \delta$. It is shown in
Lemma~\ref{lem:AppendixB} that this implies that for any infinite
increasing sequence $\{T_i\}$,  $\{z(T_i)^2\}$ is a Cauchy sequence
and therefore $z(t)^2$ converges to a limit. Since $z$ is
continuous, this implies that $z(t)$ also converges to a limit,
$z_{\infty}$. Since $y \in W_2$, we have $\lim_{t\rightarrow \infty}
y(t)=y_\infty=0$. Recall $f_c(a,b)=f_1(b)$ at points such that
$a>-p_0$ or $b > 0$. Suppose $z_\infty \ne 0$. If $z_\infty < 0$,
then $ y_\infty - z_\infty > 0$ and $f_c(a,b)=f_1(b)$ in some
neighborhood of $(z_\infty,y_\infty-z_\infty)$. If $z_\infty > 0$,
then $z_\infty>0\ge-p_0$ and $f_c(a,b)=f_1(b)$ in some neighborhood
of $(z_\infty,y_\infty-z_\infty)$. Since $f_1$ is continuous, we
have $\lim_{t\rightarrow \infty} \dot{z}(t)=\lim_{t\rightarrow
\infty} f_c(z(t),y(t)-z(t))=f_1(y_\infty-z_{\infty})=
f_1(-z_{\infty}) $. By inspection of the function $f_1$, we see that
$z_\infty \ne 0$ implies that $\dot{z}$ has a nonzero limit.
However, since $\dot{z}\in L_2$, it cannot have a nonzero limit.
Thus we conclude by condradiction that $z_\infty=0$.
%However, since $ \dot{z} \in L_2$, it cannot have a positive limit,
%therefore $z_\infty \ge 0$. Now suppose $z_\infty
%> 0$, then $z_\infty>-p_0$ since $p_0 \ge 0$. $f_c$ is continuous
%about this point and the same argument applies. If $p_0=0$, then we
%must consider the case $z_infty = 0$ since $f_c$ is discontinuous at
%this point. However, since $\norm{f_c(x,-x)} \le \beta \norm{x}$, we
%conclude that $\lim_{t\rightarrow \infty}
%\dot{z}(t)=\lim_{t\rightarrow \infty}
%f_c(z(t),y(t)-z(t))=f_c(z_{\infty},-z_{\infty})$ in this case as
%well. Finally, since $\dot{z} \in L_2$, if $\dot{z}$ has a limit, it
%must be $0$ and since $f_c(z_{\infty},-z_{\infty})=0$ implies
%$z_{\infty}=0$, we have $z_{\infty}=0$.
\end{proof}
%%%%%%%%%%%%%%%%%%%%%%%%%%%%%%%%%%%%%%%%%%%%%%%%%%%%%%%%%%%%%%%%%%%
\begin{lem}
If $v=\Delta y$ with $y \in W_{2}$, then $\ip{v}{\dot{y}-v} \ge
-\beta|y(0)|^2$.
\end{lem}

\begin{proof} Let $z=\Delta_z v$ and define the variable $r(t)=y(t)-z(t)$ and the set
$M=\{t: z(t)>-p_0 \; \text{ or } \; r(t) \ge 0\}$, then
\begin{align*}&\ip{v}{\dot{y}-v}=\ip{\dot{z}}{\dot{y}-\dot{z}}=\int_{0}^{\infty} \dot{z}(t) \dot{r}(t)dt \\
&=\int_M f_1(r(t))\dot{r}(t)dt \le \beta \norm{y}
\norm{\dot{y}}+\beta^2\norm{v}^2
\end{align*}
Since $y\in W_2$, we have that $y$ is absolutely continuous and thus
$r$ is absolutely continuous. Since $r,z$ are absolutely continuous
functions and since by Lemma~\ref{lem:z->0}, we have $z(t)
\rightarrow 0$, we can partition the set $M$ into the countable
union of sequential disjoint intervals $\bigcup_i I_i \bigcup I_f$
where $I_i=[T_{a,i},T_{b,i})$ with $\{T_{a,i}\}, \{T_{b,i}\}\subset
\R^+$ and $I_f=[T_{a,f},\infty)$. To see that the intervals are
closed on the left, suppose $I_i$ were open on the left. Then, since
$T_{a,i} \not\in M$, $z(T_{a,i})=-p_0$ and $r(T_{a,i})<0$. However,
since $r$ is continuous, $r(T_{a,i}+\eta)<0$ for $\eta$ sufficiently
small. Since $r(t)<0$ implies $\dot{z}(t) \le 0$, we have that
$z(T_{a,i}+\eta)\le -p_0$ and thus $T_{a,i}+\eta \not\in M$ for
$\eta$ sufficiently small, which is a contradiction. Thus all the
intervals are closed on the left. Similarly, one can show that all
the intervals are open on the right.

\begin{figure}[htb]
\centerline{\includegraphics[width=0.45\textwidth]{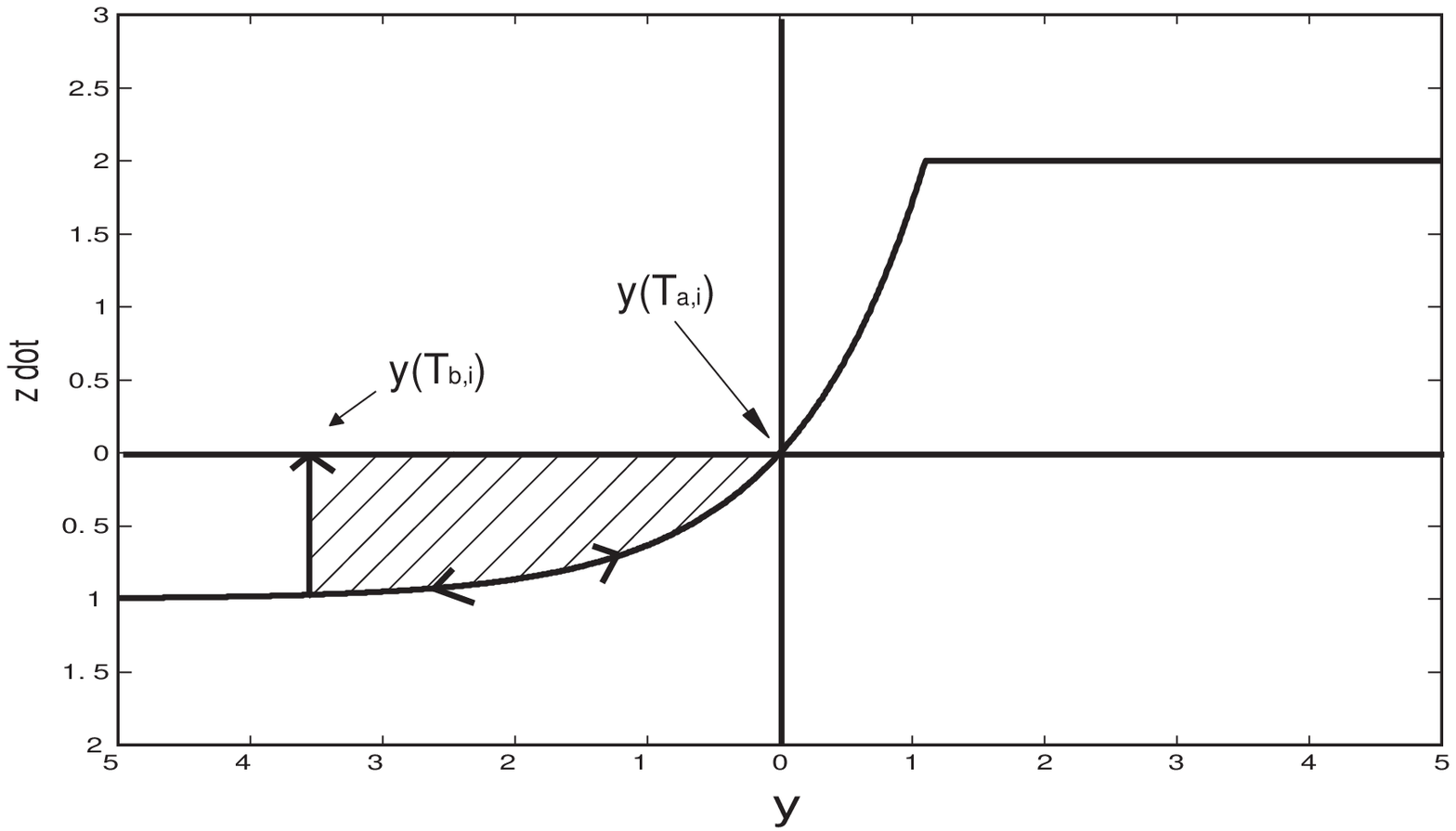}}
\caption{Value of $y,\dot{z}$ at times $T_{a,i}$ and $T_{b,i}$}
\label{fig:zdot}
\end{figure}

Now, consider time $T_{a} > 0$, where $T_{a}\in M$ defines the start
of one of the intervals described above. If $z(T_{a})>-p_0$, then
since $z$ is continuous, $z(T_{a}-\eta)>-p_0$ for all $\eta$
sufficiently small. Therefore $T_{a}-\eta \in M$ for all $\eta$
sufficiently small. This contradicts the statement that the
intervals are disjoint. We thus conclude $z(T_{a})=-p_0$ and
consequently $r(T_{a})\ge 0$ by definition of $M$. Now suppose
$r(T_{a})
> 0$. Since $r$ is continuous, $r(T_{a}-\epsilon)>0$ and consequently
$T_{a}-\epsilon \in M$ for all $\epsilon$ sufficiently small, which
contradicts the statement that the intervals are disjoint. Therefore
we conclude $r(T_{a})=0$ if $T_{a} \ne 0$. Then
\begin{align*} &\ip{v}{\dot{y}-v}=\sum_{i} \int_{I_i}
f_1(r(t))\dot{r}(t)dt+\int_{T_{a,f}}^{\infty}
f_1(r(t))\dot{r}(t)dt\\
&=\sum_{i} \int_{T_{a,i}}^{T_{b,i}} f_1(r(t))\dot{r}(t)dt
+\int_{T_{a,f}}^{\infty} f_1(r(t))\dot{r}(t)dt
\end{align*}
We will assume that $T_{a,1}=0$. If $T_{a,1}\ne 0$, we have
$r(T_{a,1})=0$ and the proof becomes simpler. Since $f_1(r)$ is
continuous in $r$ and $r(t)$ is absolutely continuous in time, by
the substitution rule we have
\begin{align*}&\ip{v}{\dot{y}-v}=\sum_{i} \int_{r(T_{a,i})}^{r(T_{b,i})}
f_1(r)dr\\
&=\int_{r(0)}^{r(T_{b,1})} f_1(r)dr +\sum_{i\ne1}
\int_0^{r(T_{b,i})} f_1(r) dr\\
&=\int_{r(0)}^0 f_1(r)dr +\sum_{i} \int_0^{r(T_{b,i})} f_1(r) dr
\end{align*}
Since $f_1 \in \text{sector}[0,\beta]$, $\int_0^{r(T_{b,i})}
f_1(r)dr \ge 0$ for any $r(T_{b,i})\in \R$. The summation converges
since it is bounded, increasing. Furthermore, since
$r(0)=y(0)-z(0)=y(0)$ and $| \int_0^y f_1(r)dr | \le f_1(y)y\le
\beta y^2$ for any $y$, we have
\[\ip{v}{\dot{y}-v}=\int_{y(0)}^0 f_1(r)dr+\sum_{i} \int_0^{r(T_{b,i})} f_1(r)dr \ge
-\beta |y(0)|^2
\]
\end{proof}
%%%%%%%%%%%%%%%%%%%%%%%%%%%%%%%%%%%%%%%%%%%%%%%%%%%%%%%%%%%%%%%%%%%%

\begin{thm} \label{thm:AppendixA} For any initial condition $x_0 \in W_{2}$ with
$x_0(-\tau) > - p_0$, $x_0(\theta) \ge -p_0$, there exists a $f \in
W_2$ and $T>0$ such that $A(x_0,t)=B_z(f,t+T)$ for $t \ge 0$.
\end{thm}
%
%\begin{lem}For any $x_0 \in W_{2}$ with $x_0(-\tau) > -p_0$, there
%exists a $f \in W_2$, $T' \ge 0$, such that
%$f(T')=x_0(-\tau)-z(T'-\tau)$
%$A_f(x_0(0),x_0,t)=B_f(x_0(0),x_0,f,t+T_1)$
%\end{lem}
\begin{proof}
Recall the map $B_z$ is defined by the following dynamics for $z=B_z
y$. $z(t)=0$ for $t \le 0$ and for $t \ge 0$
\begin{align*}
\dot{z}(t)&=f_c \bl(z(t),y(t)-z(t-\tau) \br)
\end{align*}
Now suppose we let $y(t)=f(t)+z(t-\tau)$ on a finite interval
$[0,T']$ for some $f \in W_2$. The system is still well-posed and $y
\in W_2$ if $f \in W_2$ since the derivative of $z$ is bounded. The
dynamics are now given by $z(t)=0$ for $t \le 0$ and the following
for $t \ge 0$
\begin{align*}
\dot{z}(t)&=f_c \bl( z(t),f(t) \br)
\end{align*}

\textbf{Part 1:} The first part of the proof is to construct a $f
\in W_2$ that drives the state $z(t)$ to $z(T')=x_0(0)$ for some $T'
\ge 0$. Furthermore, for continuity with Part 2, we require that
$f(T')=-x_0(-\tau)$.  If $x_0(-\tau) = 0$ and $x_0(0) = 0$, then we
are done. Otherwise, there are 4 cases to consider.\\

\paragraph{$x_0(-\tau) < 0, x_0(0)>0$}
First, let $f(t)=\epsilon t$ until
$T_1=-\frac{x_0(-\tau)}{\epsilon}$ for $\epsilon > 0$. Let
$\epsilon$ be sufficiently large so that $z(T_1)<x_0(0)$. Such an
$\epsilon$ exists since $\dot{z}$ is bounded. Let
$f(t)=-x_0(-\tau)>0$ until time $T'$ such that $z(T')=x_0(0)$. Such
$T'$ exists since $z(t)$ is now linearly
increasing. \\

\paragraph{$x_0(-\tau) > 0, x_0(0) < 0$}
This case is handled similarly to the previous one.\\

\paragraph{$x_0(-\tau) < 0, x_0(0) < 0$} For this case, let
$f(t)=-\epsilon t$ until time $T_1=\frac{1}{\epsilon}$. Then let
$f(t)=-1$ until time $T_2$ such that $z(T_2)=x_0(0)-\gamma$ for some
$\gamma>0$. Such time exists for any $\gamma<p_0-x_0(0)$ since
$z(t)$ is linearly decreasing for $z(t) \ge -p_0$. Then let
$f(t)=x_0(0)-\gamma + \lambda t$ for time $\Delta
t=\frac{1}{\lambda}(x_0(-\tau)+1)$. Make $\lambda$ sufficiently
large so that $z(T_2+\Delta t)<x_0(0)$. This is possible since
$\dot{z}$ is bounded. Finally let $f(t)=-x_0(-\tau)>0$ until time
$T'$ when $z(T')=x_0(0)$. Such $T'$
exists since $z(t)$ is now linearly increasing.\\

\paragraph{ $x_0(-\tau) >0, x_0(0)>0$ or $x_0(0)=0$ or $x_0(-\tau)=0$}
These cases are handled similarly to the previous one.\\

\textbf{Part 2:} For time $t \in [T',T'+\tau]$, let
$f(t)=-x_0(t-T')$. We then have
$y(T'+\tau)=x_0(t-T')-z(T')=x_0(-\tau)-x_0(-\tau)=0$. Let $y(t)=0$
for $t \ge T'+\tau$. Therefore $y\in W_2$ and we conclude that the
dynamics of the interconnection for time $t \in [T', T'+\tau]$ are
given by $z(T')=x_0(0)$
\begin{eqnarray}
\dot{z}(t)&=f_c(z(t),y(t)+z(t-\tau))\\
&=f_c(z(t),x_0(t-\tau))
\end{eqnarray}
And for $t \ge T'+\tau$, we have
\begin{eqnarray}
\dot{z}(t)&=f_c(z(t),y(t)+z(t-\tau))\\
&=f_c(z(t),-z(t-\tau))
\end{eqnarray}

Therefore, we have that $A(y,t+T')=B(x_0,t)$.\\
\end{proof}

\begin{lem} \label{lem:AppendixB} Suppose that for any $\delta>0$, there exists a $T_{\delta}$ such
that $T_2>T_1>T_{\delta}$ implies $z(T_2)^2-z(T_1)^2 < \delta$. Then
for any infinite, increasing sequence, $\{T_i\}$, $\{z(T_i)^2\}$ is
a Cauchy sequence. \label{lem:z2_Cauchy}
\end{lem}
\begin{proof} Proof by contradiction. Suppose that $\{z(T_i)^2\}$ is not a
Cauchy sequence. Then there exists an $\epsilon>0$ such that for any
$N > 0$, there exists $i,j > N$ such that
$|z(T_i)^2-z(T_j)^2|>\epsilon$.
%We show that this implies the
%existence of a $T_f$ such that $z(T_f)^2<0$, which is a
%contradiction.

Let $\delta=\frac{\epsilon}{2}$. By assumption there exists a
$T_{\delta}$ such that for all $T_2>T_1>T_{\delta}$,

\begin{equation}z(T_2)^2-z(T_1)^2<\frac{\epsilon}{2}
\label{eqn:app_eq1}
\end{equation}

Since $\{T_i\}$ is strictly increasing, infinite, there exists a
$N>0$ such that $T_{N} > T_{\delta}$. Since $\{z(T_i)^2\}$ is not
Cauchy, there exists $i_1,j_1 > N$ such that
$|z(T_{i_1})^2-z(T_{j_1})^2| > \epsilon$. Without loss of
generality, we may assume $i_1>j_1$. We now have one of the two
possibilities.

\begin{itemize}\item $z(T_{i_1})^2 - z(T_{j_1})^2 > \epsilon $ \item
$z(T_{i_1})^2 - z(T_{j_1})^2 < -\epsilon $
\end{itemize}
However, by Equation~\eqref{eqn:app_eq1}, the first option is not
possible since $T_i>T_j > T_N > T_{\delta}$. Therefore $z(T_{i_1})^2
< z(T_{j_1})^2 - \epsilon$. Thus for all $k > i_1$,
by~\eqref{eqn:app_eq1} and since $T_k > T_{i_1}> T_{j_1} > T_N \ge
T_\delta$, we have
\[z(T_k)^2
<z(T_{i_1})^2 + \frac{\epsilon}{2} < z(T_{j_1})^2 - \epsilon
+\frac{\epsilon}{2} = z(T_{j_1})^2 - \frac{\epsilon}{2}
\]

Now, let $N=i_1$. Again, since $\{z(T_i)^2\}$ is not Cauchy, there
exists $i_2,j_2>i_1$ such that $|z(T_{i_2})^2-z(T_{j_2})^2|
> \epsilon$. Using the method above, for all $k > i_2$ since $T_k > T_{i_2} > T_{j_2} > T_{i_1}> T_{\delta}$,
by~\eqref{eqn:app_eq1} we have
\[ z(T_k)^2
< z(T_{i_2})^2+\frac{\epsilon}{2}<z(T_{j_2})^2 - \frac{\epsilon}{2}
< z(T_{j_1})^2 - 2\frac{\epsilon}{2}
\]

Repeating $n$ times, we find a $i_n$ such that for all $k>i_n$,
 $z(T_k)^2<z(T_{j_1})^2-n \frac{\epsilon}{2}$. Since we can
assume $z(T_{j_1})^2$ is finite, let $n$ be sufficiently large and
we have the existence of an $i$ such that that $z(T_i)^2<0$, which
is a contradiction. Thus we have that $\{z(T_i)^2\}$ is a Cauchy
sequence.\end{proof}

\chapter{Appendix to
Chapter~\ref{chp:LinearCase}}\label{app:LinearCase}

\begin{thm} Suppose $M:\R \rightarrow \S^{2n}$ is a continuous matrix-valued function. Then the following are equivalent
\begin{enumerate}
\item There exists some $\epsilon>0$ such that the following holds for all
$x \in \mathcal{C}_\tau$.
\[\int_{-\tau}^0 \bmat{x(0)\\x(\theta)}^T M(\theta)
\bmat{x(0)\\x(\theta)} d \theta \ge \epsilon \norm{x}^2_2
\]
\item There exists some $\epsilon'>0$ and some continuous matrix-valued function
$T:\R \rightarrow \S^n$ such that the following holds.
\begin{align*}
&\int_{-\tau}^0 T(\theta) d \theta = 0\\
&M(\theta)+\bmat{T(\theta) & 0 \\ 0 & -\epsilon' I} \ge 0 \qquad
\text{for all } \theta \in [-\tau,0]
\end{align*}
\end{enumerate}
\end{thm}

\begin{proof}
$(2 \Rightarrow 1)$ Suppose statement 2 is true, then
%there exists
%some function $T \in \Omega$ such that
%\[M(\theta)+\bmat{T(\theta) & 0 \\ 0 & -\epsilon' I} \ge 0 \qquad \text{for all
%} \theta \in [-\tau,0]
%\]
%Then
\begin{align*}&\int_{-\tau}^0 \bmat{x(0)\\x(\theta)}^T
M(\theta) \bmat{x(0)\\x(\theta)} d \theta -\epsilon' \norm{x}_2^2\\
&= \int_{-\tau}^0 \bmat{x(0)\\x(\theta)}^T M(\theta)+\bmat{T(\theta)
& 0\\0 & -\epsilon' I} \bmat{x(0)\\x(\theta)} d \theta \ge 0.
\end{align*}
$(1 \Rightarrow 2)$ Suppose that statement 1 holds for some $M$.
Write $M$ as
\[M(\theta)=\bmat{M_{11}(\theta) & M_{12}(\theta)\\M_{12}(\theta)^T &
M_{22}(\theta)}.
\]
We first prove that $M_{22}(\theta)\ge \epsilon I$ for all $\theta
\in [-\tau,0]$. By statement 1, we have that
\[\int_{-\tau}^0 \bmat{x(0)\\x(\theta)}^T \bmat{M_{11}(\theta)
& M_{12}(\theta)\\M_{12}(\theta)^T & M_{22}(\theta)-\epsilon I}
\bmat{x(0)\\x(\theta)} d \theta \ge 0.
\]
Now suppose that $M_{22}(\theta)- \epsilon I$ is not positive
semidefinite for all $\theta \in [-\tau,0]$. Then there exists some
$x_0 \in \R^n$ and $\theta_1 \in [-\tau,0]$ such that
$x_0^T(M_{22}(\theta_1)-\epsilon I)x_0 < 0 $. By continuity of
$M_{22}$, if $\theta_1=0$ or $\theta_1=-\tau$, then there exists
some $\theta_1'\in (-\tau,0)$ such that
$x_0^T(M_{22}(\theta_1')-\epsilon I)x_0 <0$. Thus assume $\theta_1
\in (-\tau,0)$. Now, since $M_{22}$ is continuous, there exists some
$x_1$ and $\delta>0$ where $\theta_1+\delta<0$, $\theta_1-\delta
> -\tau$ and such that $x_1^T(M_{22}(\theta)-\epsilon I)x_1 \le -1$
for all $\theta \in [\theta_1-\delta,\theta_1+\delta]$. Then for
$\beta > \max \{1/(-\theta_1-\delta),1/(\tau+\theta_1-\delta)\}$,
let
\[x(\theta)=\begin{cases}
(1+\beta(\theta-(\theta_1-\delta)))x_1 & \theta \in [\theta_1-\delta-1/\beta,\theta_1-\delta]\\
x_1 & \theta\in[\theta_1-\delta,\theta_1+\delta]\\
(1-\beta(\theta-(\theta_1+\delta)))x_1 & \theta \in [\theta_1+\delta,\theta_1+\delta+1/\beta]\\
%x_0(1+\alpha \theta)& \theta \in [-1/\alpha,0]\\
0 & \text{otherwise.}
\end{cases}
\]
Then $x \in \mathcal{C}_\tau$, $x(0)=0$ and $\norm{x(\theta)}\le
\norm{x_1}$ for all $\theta \in [-\tau,0]$. Now, since $M_{22}$ is
continuous, it is bounded on $[-\tau,0]$. Therefore, there exists
some $\epsilon_2> 0$ such that $M_{22}(\theta)-\epsilon I\le
\epsilon_2 I$. Then let $\beta \ge 2\epsilon_2 \norm{x_1}^2/\delta$
and we have the following.
\begin{align*}
&\int_{-\tau}^0 \bmat{x(0)\\x(\theta)}^T M(\theta)
\bmat{x(0)\\x(\theta)} d \theta-\epsilon \norm{x}^2\\
&= \int_{-\tau}^0 x(\theta)^T (M_{22}(\theta)-\epsilon I)
x(\theta) d \theta\\
&= \int_{\theta_1-\delta}^{\theta_1+\delta} x_1^T
(M_{22}(\theta)-\epsilon I) x_1 d
\theta\\
&+\int_{\theta_1-\delta-1/\beta}^{\theta_1-\delta} x(\theta)^T
(M_{22}(\theta)-\epsilon I) x(\theta) d
\theta\\
&+\int_{\theta_1+\delta}^{\theta_1+\delta+1/\beta} x(\theta)^T
(M_{22}(\theta)-\epsilon I)
x(\theta) d \theta\\
&\le -2 \delta +\epsilon_2
\int_{\theta_1-\delta-1/\beta}^{\theta_1-\delta} \norm{x(\theta)}^2
d \theta+\epsilon_2 \int_{\theta_1+\delta}^{\theta_1+\delta+1/\beta}
\norm{x(\theta)}^2 d \theta\\
&\le -2 \delta +2 \epsilon_2 \norm{x_1}^2/\beta \le -\delta\\
\end{align*}
Therefore, by contradiction, we have that $M_{22}(\theta) \ge
\epsilon I$ for all $\theta \in [-\tau,0]$. Now we define
$\epsilon'=\epsilon/2$ and
$\tilde{M}_{22}(\theta)=M_{22}(\theta)-\epsilon'I \ge \epsilon' I $.
We now show that $\tilde{M}_{22}(\theta)^{-1}$ is continuous. We
first note that $\tilde{M}_{22}(\theta)^{-1}$ is bounded, since
\begin{align*}
\tilde{M}_{22}(\theta) &\ge \epsilon' I \quad \Rightarrow I \ge
\epsilon' \tilde{M}_{22}(\theta)^{-1} \quad \Rightarrow
\tilde{M}_{22}(\theta)^{-1} \le \frac{1}{\epsilon'} I
\end{align*}
Now since $\tilde{M}_{22}$ is continuous, for any $\beta > 0$, there
exists a $\delta>0$ such that
\[|\theta_1-\theta_2|\le \delta \Rightarrow
\norm{\tilde{M}_{22}(\theta_1)-\tilde{M}_{22}(\theta_2)}\le \beta
\]
Therefore, for any $\beta'>0$, let $\delta$ be defined as above for
$\beta=\beta'\epsilon'^2$. Then $|\theta_1-\theta_2| \le \delta$
implies
\begin{align*}
&\norm{\tilde{M}_{22}(\theta_1)^{-1}-\tilde{M}_{22}(\theta_2)^{-1}}\\
&=\norm{\tilde{M}_{22}(\theta_1)^{-1} \tilde{M}_{22}(\theta_2)
\tilde{M}_{22}(\theta_2)^{-1}-\tilde{M}_{22}(\theta_1)^{-1}
\tilde{M}_{22}(\theta_1) \tilde{M}_{22}(\theta_2)^{-1}}\\
&=\norm{\tilde{M}_{22}(\theta_1)^{-1} \left(\tilde{M}_{22}(\theta_2)
-\tilde{M}_{22}(\theta_1)\right) \tilde{M}_{22}(\theta_2)^{-1}}
\\
& \le \norm{\tilde{M}_{22}(\theta_1)^{-1}}
\norm{\tilde{M}_{22}(\theta_2) -\tilde{M}_{22}(\theta_1)}\norm{
\tilde{M}_{22}(\theta_2)^{-1}} \le \frac{1}{\epsilon'^2}
\norm{\tilde{M}_{22}(\theta_2) -\tilde{M}_{22}(\theta_1)} \le
\beta'.
\end{align*}
Therefore $\tilde{M}_{22}(\theta)^{-1}$ is continuous.
We now prove statement 2 by construction. Suppose $M$ satisfies
statement 1. Let
\[T(\theta)=T_0-(M_{11}(\theta)-
M_{12}(\theta)\tilde{M}_{22}^{-1}(\theta)M_{12}(\theta)^T).
\]
Here
\[T_0=\frac{1}{\tau}\int_{-\tau}^{0}(M_{11}(\theta)-
M_{12}(\theta)\tilde{M}_{22}^{-1}(\theta)M_{12}(\theta)^T) d \theta.
\] Then we have that $T$ is continuous and
\[\int_{-\tau}^0 T(\theta) d \theta=\tau T_0-\tau T_0=0.\]
This implies that $T \in \Omega$. We now prove that $T_0 \ge 0$. For
any vector $z_0\in \R^n$, suppose $z$ is a continuous function such
that $z(0)=(I+\tilde{M}_{22}(0)^{-1} M_{12}(0)^T)z_0$. Then let
$x(\theta)=z(\theta)-\tilde{M}_{22}(\theta)^{-1} M_{12}(\theta)^T
z_0$. Then $x$ is continuous, $x(0)=z_0$ and by statement 1 we have
the following.
\begin{align*}
 &\int_{-\tau}^0 \bmat{x(0)\\x(\theta)}^T
\bmat{M_{11}(\theta) & M_{12}(\theta)\\
M_{12}(\theta)^T & \tilde{M}_{22}(\theta)} \bmat{x(0)\\x(\theta)} d \theta\\
&=\int_{-\tau}^0 \bmat{z_0\\z(\theta)}^T \bmat{I &
0\\-\tilde{M}_{22}(\theta)^{-1} M_{12}(\theta)^T & I}^T
\bmat{M_{11}(\theta) & M_{12}(\theta)\\ M_{12}(\theta)^T &
\tilde{M}_{22}(\theta)}\\
&\qquad \qquad \qquad \qquad \qquad \qquad \qquad \quad \bmat{I &
0\\-\tilde{M}_{22}(\theta)^{-1} M_{12}(\theta)^T & I}
\bmat{z_0\\z(\theta)} d
\theta \\
&=\int_{-\tau}^0 \bmat{z_0\\z(\theta)}^T \bmat{M_{11}(\theta)-
M_{12}(\theta)\tilde{M}_{22}^{-1}(\theta)M_{12}(\theta)^T & 0\\
0 & \tilde{M}_{22}(\theta)} \bmat{z_0\\z(\theta)} d \theta \\
&=z_0^T(\int_{-\tau}^0 M_{11}(\theta)-
M_{12}(\theta)\tilde{M}_{22}^{-1}(\theta)M_{12}(\theta)^T d \theta)
z_0\\
&+ \int_{-\tau}^0 z(\theta)^T \tilde{M}_{22}(\theta) z(\theta) d \theta\\
&=z_0^T T_0 z_0 + \int_{-\tau}^0 z(\theta)^T \tilde{M}_{22}(\theta)
z(\theta) d \theta \ge \epsilon' \norm{x}_2^2
\end{align*}
We now show that this implies that $T_0 \ge 0$. Suppose there exists
some $y$ such that $y^T T_0 y <0$. Then there exists some $z_0$ such
that $z_0 T_0 z_0=-1$. Now let $\alpha>1/\tau$ and
\[z(\theta)=\begin{cases}
(I+\tilde{M}_{22}(0)^{-1} M_{12}(0)^T)z_0(1+\alpha \theta)& \theta \in [-1/\alpha,0]\\
0 & \text{otherwise.}
\end{cases}
\]
Then $z$ is continuous,
$z(0)=(I+\tilde{M}_{22}(0)^{-1}M_{12}(0)^T)z_0$ and
$\norm{z(\theta)}^2\le\norm{z(0)}^2$ for all $\theta\in[-\tau,0]$.
Recall that $\tilde{M}_{22}(\theta) \le
(\epsilon_2+\epsilon')I=\epsilon_3 I$. Let
$\alpha>2\epsilon_3\norm{z(0)^2}$ and then we have the following.
\begin{align*}
&z_0^T T_0 z_0 + \int_{-\tau}^0 z(\theta)^T M_{22}(\theta) z(\theta)
d \theta \le -1+\epsilon_3 \int_{-1/\alpha}^0 \norm{z(\theta)}^2 \\
&\le -1+\epsilon_3 \norm{z(0)}^2/\alpha < -\frac{1}{2}
\end{align*}
However, this contradicts the previous relationship. Therefore, we
have by contradiction that $T_0\ge0$. Now by using the invertibility
of $\tilde{M}_{22}(\theta)=M_{22}(\theta)-\epsilon'I\ge \epsilon' I$
and the Schur complement transformation, we have that statement 2 is
equivalent to the following for all $\theta \in [-\tau,0]$.
\begin{align*}
&\tilde{M}_{22}(\theta)\ge 0\\
&M_{11}(\theta)+T(\theta)-M_{12}(\theta)\tilde{M}_{22}(\theta)^{-1}
M_{12}(\theta)^T \ge 0
\end{align*}
We have already proven the first condition. Finally, we have the
following.
\begin{align*}
&M_{11}(\theta)+T(\theta)-M_{12}(\theta)\tilde{M}_{22}(\theta)^{-1}
M_{12}(\theta)^T=T_0 \ge 0
\end{align*}
Thus we have shown that statement 2 is true.
\end{proof}

%%%%%%%%%%%%%%%%%%%%%%%%%%%%%%%%%%%%%%%%%%%%%%%%%%%%%%%%%%%%%%%%%%%%%%
% Extension Proof 2
%%%%%%%%%%%%%%%%%%%%%%%%%%%%%%%%%%%%%%%%%%%%%%%%%%%%%%%%%%%%%%%%%%%%%%

\begin{lem} Let $M:\R \mapsto \S^{3n}$ be a continuous
matrix-valued function. Then the following are equivalent.
\begin{enumerate}
\item There exists an $\epsilon>0$ such that the following holds for
all $x \in \mathcal{C}_\tau$.
\[\int_{-\tau}^0\bmat{x(0)\\x(-\tau)\\x(\theta)}^T M(\theta)\bmat{x(0)\\x(-\tau)\\x(\theta)}\ge \epsilon \norm{x}_2^2
\]
\item There exists some $\epsilon'>0$ and a continuous matrix valued
function $T:\R \mapsto \S^{2n}$ such that the following holds.
\begin{align*}
&\int_{-\tau}^0 T(\theta) d \theta = 0\\
&M(\theta)+\bmat{T(\theta) & 0 \\ 0 & -\epsilon'I} \ge 0 \qquad
\text{for all } \theta \in [-\tau,0]
\end{align*}
\end{enumerate}
\end{lem}
\begin{proof}
$(2 \Rightarrow 1)$ Suppose there exists some function $T \in S_0$
such that
\[M(\theta)+\bmat{T(\theta) & 0 \\ 0 & -\epsilon' I} \ge 0 \qquad \text{for all
} \theta \in [-\tau,0].
\]

Then we have the following.
\begin{align*}&\int_{-\tau}^0 \bmat{x(0)\\x(-\tau)\\x(\theta)}^T
M(\theta) \bmat{x(0)\\x(-\tau)\\x(\theta)} d \theta -\epsilon' \norm{x}_2^2\\
&= \int_{-\tau}^0 \bmat{x(0)\\x(-\tau)\\x(\theta)}^T
M(\theta)+\bmat{T(\theta) & 0\\0 & -\epsilon' I}
\bmat{x(0)\\x(-\tau)\\x(\theta)} d \theta \ge 0
\end{align*}

$(1 \Rightarrow 2)$ Suppose that statement 1 holds for some $M$.
Write $M$ as
\[M(\theta)=\bmat{M_{11}(\theta) & M_{12}(\theta)\\M_{12}(\theta)^T &
M_{22}(\theta)},
\]

where $M_{22}:\R\mapsto \S^n$. We first prove that
$M_{22}(\theta)\ge \epsilon I$ for all $\theta \in [-\tau,0]$. By
statement 1, we have that

\[\int_{-\tau}^0 \bmat{x(0)\\x(-\tau)\\x(\theta)}^T \bmat{M_{11}(\theta)
& M_{12}(\theta)\\M_{12}(\theta)^T & M_{22}(\theta)-\epsilon I}
\bmat{x(0)\\x(-\tau)\\x(\theta)} d \theta \ge 0.
\]
Now suppose that $M_{22}(\theta)- \epsilon I$ is not positive
semidefinite for all $\theta \in [-\tau,0]$. Then there exists some
$x_0$ and $\theta_1 \in [-\tau,0]$ such that
$x_0^T(M_{22}(\theta_1)-\epsilon I)x_0 < 0 $. By continuity of
$M_{22}$, if $\theta_1=0$ or $\theta_1=-\tau$, then there exists
some $\theta_1'\in (-\tau,0)$ such that
$x_0^T(M_{22}(\theta_1')-\epsilon I)x_0 <0$. Thus assume $\theta_1
\in (-\tau,0)$. Now, since $M_{22}$ is continuous, there exists some
$x_1$ and $\delta>0$ where $\theta_1+\delta<0$, $\theta_1-\delta
> -\tau$ and such that $x_1^T(M_{22}(\theta)-\epsilon I)x_1 \le -1$ for
$\theta \in [\theta_1-\delta,\theta_1+\delta]$. Then for $\beta >
\max \{1/(-\theta_1-\delta),1/(\tau+\theta_1-\delta)\}$, let
\[x(\theta)=\begin{cases}
\beta(\theta-(\theta_1-\delta-1/\beta))x_1 & \theta \in [\theta_1-\delta-1/\beta,\theta_1-\delta]\\
x_1 & \theta\in[\theta_1-\delta,\theta_1+\delta]\\
(1-\beta(\theta-(\theta_1+\delta)))x_1 & \theta \in [\theta_1+\delta,\theta_1+\delta+1/\beta]\\
%x_0(1+\alpha \theta)& \theta \in [-1/\alpha,0]\\
0 & \text{otherwise.}
\end{cases}
\]

Then $x \in \mathcal{C}_\tau$, $x(-\tau)=x(0)=0$ and
$\norm{x(\theta)}^2 \le \norm{x_1}^2$ for all $\theta \in
[-\tau,0]$. Now, since $M_{22}$ is continuous, it is bounded on
$[-\tau,0]$. Therefore, there exists some $\epsilon_2> 0$ such that
$M_{22}(\theta)-\epsilon I\le \epsilon_2 I$ for $\theta \in
[-\tau,0]$.
%and since $x_1$ is fixed, there exists $b,c>0$ such that
%\begin{align*}
%c=\sup_{\theta \in [-\tau,0], \norm{x_2} \le \norm{x_1}}\left|x_2^T
%M_{22}(\theta)x_2 \right|
%\end{align*}
Let $\beta \ge 2\epsilon_2 \norm{x_1}^2/\delta$ and we have the
following.
\begin{align*}
&\int_{-\tau}^0 \bmat{x(0)\\x(-\tau)\\x(\theta)}^T M(\theta)
\bmat{x(0)\\x(-\tau)\\x(\theta)} d \theta-\epsilon \norm{x}^2 \\
&= \int_{-\tau}^0 x(\theta)^T (M_{22}(\theta)-\epsilon I)
x(\theta) d \theta\\
&= \int_{\theta_1-\delta}^{\theta_1+\delta} x_1^T
(M_{22}(\theta)-\epsilon I) x_1 d
\theta+\int_{\theta_1-\delta-1/\beta}^{\theta_1-\delta} x(\theta)^T
(M_{22}(\theta)-\epsilon I) x(\theta) d
\theta\\
&+\int_{\theta_1+\delta}^{\theta_1+\delta+1/\beta} x(\theta)^T
(M_{22}(\theta)-\epsilon I)
x(\theta) d \theta\\
&\le -2 \delta +\epsilon_2
\int_{\theta_1-\delta-1/\beta}^{\theta_1-\delta} \norm{x(\theta)}^2
d \theta+\epsilon_2 \int_{\theta_1+\delta}^{\theta_1+\delta+1/\beta}
\norm{x(\theta)}^2 d \theta\\
&\le -2 \delta +2 \epsilon_2 \norm{x_1}^2/\beta \le -\delta\\
\end{align*}
Therefore, by contradiction, we have that $M_{22}(\theta) \ge
\epsilon I$ for all $\theta \in [-\tau,0]$. Now we define
$\epsilon'=\epsilon/2$ and
$\tilde{M}_{22}(\theta)=M_{22}(\theta)-\epsilon'I \ge \epsilon' I $.
We now show that $\tilde{M}_{22}(\theta)^{-1}$ is continuous. We
first note that $\tilde{M}_{22}(\theta)^{-1}$ is bounded, since
\begin{align*}
\tilde{M}_{22}(\theta) &\ge \epsilon' I \\
&\Rightarrow I \ge \epsilon' \tilde{M}_{22}(\theta)^{-1}\\
&\Rightarrow \tilde{M}_{22}(\theta)^{-1} \le \frac{1}{\epsilon'} I
\end{align*}
Now since $\tilde{M}_{22}$ is continuous, for any $\beta > 0$, there
exists a $\delta>0$ such that
\[|\theta_1-\theta_2|\le \delta \Rightarrow
\norm{\tilde{M}_{22}(\theta_1)-\tilde{M}_{22}(\theta_2)}\le \beta.
\]
Therefore, for any $\beta'>0$, let $\delta$ be defined as above for
$\beta=\beta'\epsilon'^2$. Then $|\theta_1-\theta_2| \le \delta$
implies
\begin{align*}
&\norm{\tilde{M}_{22}(\theta_1)^{-1}-\tilde{M}_{22}(\theta_2)^{-1}}\\
&=\norm{\tilde{M}_{22}(\theta_1)^{-1} \tilde{M}_{22}(\theta_2)
\tilde{M}_{22}(\theta_2)^{-1}-\tilde{M}_{22}(\theta_1)^{-1}
\tilde{M}_{22}(\theta_1) \tilde{M}_{22}(\theta_2)^{-1}}\\
&=\norm{\tilde{M}_{22}(\theta_1)^{-1} \left(\tilde{M}_{22}(\theta_2)
-\tilde{M}_{22}(\theta_1)\right) \tilde{M}_{22}(\theta_2)^{-1}}\\
&\le \norm{\tilde{M}_{22}(\theta_1)^{-1}}
\norm{\tilde{M}_{22}(\theta_2)
-\tilde{M}_{22}(\theta_1)}\norm{ \tilde{M}_{22}(\theta_2)^{-1}}\\
&\le \frac{1}{\epsilon'^2} \norm{\tilde{M}_{22}(\theta_2)
-\tilde{M}_{22}(\theta_1)} \le \beta'.
\end{align*}
Therefore $\tilde{M}_{22}(\theta)^{-1}$ is continuous.

We now prove statement 2 by construction. Suppose $M$ satisfies
statement 1. Let

\[T(\theta)=T_0-(M_{11}(\theta)-
M_{12}(\theta)\tilde{M}_{22}^{-1}(\theta)M_{12}(\theta)^T)
\]
Where
\[T_0=\frac{1}{\tau}\int_{-\tau}^{0}(M_{11}(\theta)-
M_{12}(\theta)\tilde{M}_{22}^{-1}(\theta)M_{12}(\theta)^T) d \theta
\] Then we have that $T$ is continuous and
\[\int_{-\tau}^0 T(\theta) d \theta=\tau T_0- \tau T_0=0\]
Which implies that $T \in \Omega$. We now show that $T_0\ge 0$. For
any constant vector $z_0=\bmat{z_1 &z_2}^T$, suppose that $z$ is a
continuous function such that
\begin{align*}
z(0)&=z_1+\tilde{M}_{22}(0)^{-1} M_{12}(0)^T z_0\\
z(-\tau)&=z_2+\tilde{M}_{22}(-\tau)^{-1} M_{12}(-\tau)^T z_0.
\end{align*}
Then, let $x(\theta)=z(\theta)-\tilde{M}_{22}(\theta)^{-1}
M_{12}(\theta)^T z_0$. Then $x$ is continuous,
$\bmat{x(0)&x(-\tau)}^T=z_0$ and by statement 1 we have the
following.
\begin{align*}
 &\int_{-\tau}^0 \bmat{x(0)\\x(-\tau)\\x(\theta)}^T
\bmat{M_{11}(\theta) & M_{12}(\theta)\\
M_{12}(\theta)^T & \tilde{M}_{22}(\theta)} \bmat{x(0)\\x(-\tau)\\x(\theta)} d \theta\\
&=\int_{-\tau}^0 \bmat{z_0\\z(\theta)}^T \bmat{I &
0\\-\tilde{M}_{22}(\theta)^{-1} M_{12}(\theta)^T & I}^T
\bmat{M_{11}(\theta) & M_{12}(\theta)\\ M_{12}(\theta)^T &
\tilde{M}_{22}(\theta)} \\
& \qquad \qquad \qquad \qquad \qquad \qquad \qquad \quad \bmat{I &
0\\-\tilde{M}_{22}(\theta)^{-1} M_{12}(\theta)^T & I}
\bmat{z_0\\z(\theta)} d
\theta \\
&=\int_{-\tau}^0 \bmat{z_0\\z(\theta)}^T \bmat{M_{11}(\theta)-
M_{12}(\theta)\tilde{M}_{22}^{-1}(\theta)M_{12}(\theta)^T & 0\\
0 & \tilde{M}_{22}(\theta)} \bmat{z_0\\z(\theta)} d \theta \\
&=z_0^T(\int_{-\tau}^0 M_{11}(\theta)-
M_{12}(\theta)\tilde{M}_{22}^{-1}(\theta)M_{12}(\theta)^T d \theta)
z_0
+ \int_{-\tau}^0 z(\theta)^T \tilde{M}_{22}(\theta) z(\theta) d \theta\\
&=z_0^T T_0 z_0
+ \int_{-\tau}^0 z(\theta)^T \tilde{M}_{22}(\theta) z(\theta) d \theta \ge \epsilon' \norm{x}_2^2\\
\end{align*}

We now show that this implies that $T_0 \ge 0$. Suppose there exists
some $y$ such that $y^T T_0 y <0$. Then there exists some $z_0$ such
that $z_0 T_0 z_0=-1$. Now let $z_0=\bmat{z_1 & z_2}^T$,
$c_1=z_1+\tilde{M}_{22}(0)^{-1} M_{12}(0)^T z_0$,
$c_2=z_2+\tilde{M}_{22}(-\tau)^{-1} M_{12}(-\tau)^T z_0$,
$\alpha>2/\tau$ and
\[z(\theta)=\begin{cases}
c_1(1+\alpha \theta)& \theta \in [-1/\alpha,0]\\
c_2(1-(\theta+\tau) \alpha)& \theta \in [-\tau,-\tau+1/\alpha]\\
0 & \text{otherwise}
\end{cases}.
\]
Then $z$ is continuous, $z(0)=c_1=z_1+\tilde{M}_{22}(0)^{-1}
M_{12}(0)^T z_0$, $z(-\tau)=c_2=z_2+\tilde{M}_{22}(-\tau)^{-1}
M_{12}(-\tau)^T z_0$, and $\norm{z(\theta)}^2\le
c=\max\{\norm{c_1}^2,\norm{c_2}^2\}$ for all $\theta\in[-\tau,0]$.
Recall that $\tilde{M}_{22}(\theta) \le
(\epsilon_2+\epsilon')I=\epsilon_3 I$. Let $\alpha>2\epsilon_3 c$
and then we have
\begin{align*}
z_0^T T_0 z_0 + \int_{-\tau}^0 z(\theta)^T M_{22}(\theta) z(\theta)
d \theta &\le -1+\epsilon_3 \int_{-1/\alpha}^0 \norm{z(\theta)}^2 \\
&\le -1+\epsilon_3 c/\alpha < -\frac{1}{2}.
\end{align*}

This contradicts the previous statement. Thus we have by
contradiction that $T_0\ge0$. Now by using the invertibility of
$\tilde{M}_{22}=M_{22}(\theta)-\epsilon'I\ge \epsilon' I$ and the
Schur complement transformation, we have that statement 2 is
equivalent to the following.
\begin{align*}
&\tilde{M}_{22}(\theta)\ge 0\\
&M_{11}(\theta)+T(\theta)-M_{12}(\theta)\tilde{M}_{22}(\theta)^{-1}
M_{12}(\theta)^T \ge 0
\end{align*}
The first condition has already been proven. Finally, we have the
following.
\begin{align*}
&M_{11}(\theta)+T(\theta)-M_{12}(\theta)\tilde{M}_{22}(\theta)^{-1}
M_{12}(\theta)^T=T_0 \ge 0
\end{align*}
Thus we have shown that statement 2 is true.
\end{proof}

%%%%%%%%%%%%%%%%%%%%%%%%%%%%%%%%%%%%%%%%%%%%%%%%%%%%%%%%%%%%%%%%%%
% Extension Proof 3
%%%%%%%%%%%%%%%%%%%%%%%%%%%%%%%%%%%%%%%%%%%%%%%%%%%%%%%%%%%%%%%%%%%%%%

\begin{lem}\label{lem:extproof3lem}
Suppose $M$ is a matrix valued function which is continuous except
possibly at points $\{\tau_i\}_{i=1}^{K-1}$ where $\tau_i \in
[-\tau_K,0]$.
%Define
%\begin{align*}
%Y_1&:=\{\bmat{x_1&x_2&x_3}^T: x_3\in \mathcal{C}_{\tau}, \;
%x_1=x_3(0),\;
%x_2=x_3(-\tau)\}\\
%Y_1&:=\{\bmat{x_1&x_2&x_3}^T: \text{ $x_3$ is continuous except at
%points $\tau_i$}\}
%\end{align*}
Then the following are equivalent.

\begin{enumerate}
\item The following holds for all $x\in \mathcal{C}_{\tau_K}$.
\[\int_{-\tau_K}^0 \bmat{x(0)\\x(\theta)}^T M(\theta) \bmat{x(0)\\x(\theta)} \ge 0\]
\item The following holds for all $x$ where $x$ is continuous except
possibly at points $\{-\tau_i\}_{i=1}^{K-1}$.
\[\int_{-\tau_K}^0 \bmat{x(0)\\x(\theta)}^T M(\theta) \bmat{x(0)\\x(\theta)} \ge 0\]
\end{enumerate}

\end{lem}

\begin{proof}
Clearly statement 2 implies statement 1. Now suppose statement 2 is
false. Then there exists some piecewise continuous $x$ , $\epsilon
>0$ such that
\[V_1=\int_{-\tau_K}^0 \bmat{x(0)\\x(\theta)}M(\theta)\bmat{x(0)\\x(\theta)} d\theta\le
-\epsilon.
\]
Now let
\begin{align*}
&\hat{x}(\theta)\\
&=\begin{cases} x(-\tau_i-\hspace{-1mm}\frac{1}{\beta})\hspace{-1mm}
+\hspace{-1mm}(x(-\tau_i+\frac{1}{\beta}) -\hspace{-1mm}
x(-\tau_i-\hspace{-1mm}\frac{1}{\beta}))(\beta(\theta+\tau_i+\hspace{-1mm}\frac{1}{\beta})/2)
& \begin{matrix}\theta\hspace{-1mm} \in
\hspace{-1mm}[-\tau_i-\frac{1}{\beta},-\tau_i+\frac{1}{\beta}],\\
i=1,\ldots, K-1\end{matrix}\\x(\theta) &
\text{otherwise.}\end{cases}
\end{align*}
Then $\hat{x}$ is continuous. Since $x$ is piecewise continuous on
$[-\tau_K,0]$, it is bounded on this interval. Now let
$c=\max_{\theta \in [-\tau_{K},0]}{\norm{x(\theta)}^2}$. Then
$\norm{\hat{x}(\theta)}^2 \le c$ for all $\theta \in [-\tau_k,0]$.
%and thus
%\begin{align*}\norm{\hat{x}}_2^2&=\int_{-\tau_K}^0\norm{\hat{x}(\theta)}^2d\theta
%=\norm{x}_2^2+\sum_{i=1}^K
%\int_{-\tau_i-1/\beta}^{-\tau_i+1/\beta}\left(\norm{\hat{x}(\theta)}^2
%-\norm{x(\theta)}^2\right)\\
%&\le\norm{x}_2^2 + 4 K c/\beta
%\end{align*}
%Then $\norm{\hat{x}}_2^2\le \norm{x}_2^2+4Kc/\beta \le 2
%\norm{x}_2^2$ for $\beta>4Kc/\norm{x}_2^2$.
Now since $M$ is piecewise continuous on $[-\tau_k,0]$, there exists
some $\epsilon'>0$ such that $M(\theta)<\epsilon' I$ for $\theta \in
[-\tau_k,0]$. Therefore we have the following.
\begin{align*}
&\int_{-\tau}^0 \bmat{\hat{x}(0)\\\hat{x}(\theta)}^T
M(\theta)\bmat{\hat{x}(0)\\\hat{x}(\theta)} d\theta\\
&=V_1+\sum_{i=1}^{K-1} \int_{-\tau_i-1/\beta}^{-\tau_i+1/\beta}
\left(\bmat{\hat{x}(0)\\\hat{x}(\theta)}M(\theta)\bmat{\hat{x}(0)\\\hat{x}(\theta)}-\bmat{x(0)\\x(\theta)}
M(\theta)\bmat{x(0)\\x(\theta)}\right) d \theta\\
&\le V_1 + \sum_{i=1}^{K-1} \int_{-\tau_i-1/\beta}^{-\tau_i+1/\beta}
\epsilon'\left(\left(\norm{\hat{x}(\theta)}^2+\norm{x(0)}^2\right)
+\left(\norm{x(\theta)}^2+\norm{x(0)}^2\right)\right)d\theta\\
&\le V_1 + 2\epsilon'(K-1)(4c)/\beta\\
&\le -\epsilon +
8\epsilon'(K-1)c/\beta \\
&\le-\epsilon/2
\end{align*}
Which holds for $\beta \ge \frac{16\epsilon'(K-1)c}{\epsilon}$. Thus
statement 2 is false implies statement 1 is false. Therefore
statement 1 implies statement 2.
\end{proof}

\begin{lem}
Suppose $S_i:\R \mapsto \S^{2n}$, $i=1 \ldots K$ are continuous
symmetric matrix valued functions with domains
$[-\tau_i,-\tau_{i-1}]$ where $\tau_K > \tau_i > \tau_{i-1}
>\tau_0=0$ for $i=2\ldots K-1$. Then the following are
equivalent.
\begin{enumerate}
\item There exists an $\epsilon>0$ such that the following holds
for all $x\in \mathcal{C}_{\tau_K}$.
\[
\sum_{i=1}^K \int_{-\tau_i}^{-\tau_{i-1}}
\bmat{x(0)\\x(\theta)}S_i(\theta)
\bmat{x(0)\\x(\theta)}d \theta \ge \epsilon \norm{x}_2^2\\
\]
\item There exists an $\epsilon'>0$ and continuous symmetric matrix
valued functions, $T_i:\R \mapsto \S^n$, such that
\begin{align*}
&S_i(\theta)+\bmat{T_i(\theta) & 0\\0&-\epsilon'I} \ge 0 \qquad \text{for }\theta \in [-\tau_i,-\tau_{i-1}]\qquad i=1 \ldots K\\
&\sum_{i=1}^K \int_{-\tau_i}^{-\tau_{i-1}}T_i(\theta)=0.
\end{align*}
\end{enumerate}
\end{lem}

\begin{proof}
$(2 \Rightarrow 1)$  Suppose there exist continuous symmetric matrix
valued functions, $T_i:\R \mapsto \S^n$, such that
\begin{align*}
&S_i(\theta)+\bmat{T_i(\theta) & 0\\0&-\epsilon'I} \ge 0 \qquad i=1 \ldots K\\
&\sum_{i=1}^K \int_{-\tau_i}^{-\tau_{i-1}}T_i(\theta)=0.
\end{align*}

Then

\begin{align*}&\sum_{i=1}^K \int_{-\tau_i}^{-\tau_{i-1}}
\bmat{x(0)\\x(\theta)}S_i(\theta)
\bmat{x(0)\\x(\theta)}d \theta - \epsilon' \norm{x}_2^2\\
&=\sum_{i=1}^K \int_{-\tau_i}^{-\tau_{i-1}}
\bmat{x(0)\\x(\theta)}S_i(\theta)+\bmat{T_i(\theta) &0 \\0&
-\epsilon' I} \bmat{x(0)\\x(\theta)}d \theta \ge 0.
\end{align*}

$(1 \Rightarrow 2)$ Suppose that statement 1 holds for some $S_i$.
Write $S_i$ as
\[S_i(\theta)=\bmat{M_{11i}(\theta) & M_{12i}(\theta)\\M_{12i}(\theta)^T &
M_{22i}(\theta)}.
\]

We first prove that $M_{22i}(\theta)\ge \epsilon I$ for all $\theta
\in [-\tau_i,-\tau_{i-1}]$, $i=1,\ldots, K$. By statement 1, we have
that

\[\sum_{i=1}^K \int_{-\tau_i}^{-\tau_{i-1}}
\bmat{x(0)\\x(\theta)}\bmat{M_{11i}(\theta) &
M_{12i}(\theta)\\M_{12i}(\theta)^T & M_{22i}(\theta)-\epsilon I}
\bmat{x(0)\\x(\theta)}d \theta \ge 0.\\
\]
Now suppose that $M_{22i}(\theta)- \epsilon I$ is not positive
semidefinite for all $\theta \in [-\tau_i,-\tau_{i-1}]$. Then there
exists some $x_0\in \R^n$ and $\theta_1 \in [-\tau_i,-\tau_{i-1}]$
such that $x_0^T(M_{22i}(\theta_1)-\epsilon I)x_0 < 0 $. By
continuity of $M_{22i}$, if $\theta_1=-\tau_i$ or
$\theta_1=-\tau_{i-1}$, then there exists some $\theta_1'\in
(-\tau_i,-\tau_{i-1})$ such that $x_0^T(M_{22i}(\theta_1')-\epsilon
I)x_0 <0$. Thus assume $\theta_1 \in (-\tau_i,-\tau_{i-1})$. Now,
since $M_{22i}$ is continuous, there exists some $x_1$ and
$\delta>0$ where $\theta_1+\delta<-\tau_{i-1}$, $\theta_1-\delta
> -\tau_i$ and such that $x_1^T (M_{22i}(\theta)-\epsilon I)x_1 \le -1$ for
$\theta \in [\theta_1-\delta,\theta_1+\delta]$. Then for $\beta >
\max \{1/(-\tau_i-\theta_1-\delta),1/(\tau_i+\theta_1-\delta)\}$,
let
\[x(\theta)=\begin{cases}
\beta(\theta-(\theta_1-\delta-1/\beta))x_1 & \theta \in [\theta_1-\delta-1/\beta,\theta_1-\delta]\\
x_1 & \theta\in[\theta_1-\delta,\theta_1+\delta]\\
(1-\beta(\theta-(\theta_1+\delta)))x_1 & \theta \in [\theta_1+\delta,\theta_1+\delta+1/\beta]\\
%x_0(1+\alpha \theta)& \theta \in [-1/\alpha,0]\\
0 & \text{otherwise.}
\end{cases}
\]

Then $x \in \mathcal{C}_{\tau_K}$, $x(0)=0$ and $\norm{x(\theta)}^2
\le \norm{x_1}^2$ for $\theta \in [-\tau_K,0]$. Now, since every
$M_{22i}$ is continuous, they are bounded on
$[-\tau_i,-\tau_{i-1}]$. Therefore, there exists some $\epsilon_2>
0$ such that $M_{22i}(\theta)-\epsilon I\le \epsilon_2 I$ for
$\theta \in [-\tau_i,-\tau_{i-1}]$, $i=1,\ldots, K$. Then let $\beta
\ge 2\epsilon_2 \norm{x_1}^2/\delta$ and we have the following.
\begin{align*}
&\sum_{j=1}^K \int_{-\tau_j}^{-\tau_{j-1}} \bmat{x(0)\\x(\theta)}^T
S_j(\theta) \bmat{x(0)\\x(\theta)} d \theta-\epsilon \norm{x}^2 =
\sum_{j=1}^K \int_{-\tau_j}^{-\tau_{j-1}} x(\theta)^T
(M_{22j}(\theta)-\epsilon I)
x(\theta) d \theta\\
&= \int_{\theta_1-\delta}^{\theta_1+\delta} x_1^T
(M_{22i}(\theta)-\epsilon I) x_1 d
\theta+\int_{\theta_1-\delta-1/\beta}^{\theta_1-\delta} x(\theta)^T
(M_{22i}(\theta)-\epsilon I) x(\theta) d
\theta\\
&+\int_{\theta_1+\delta}^{\theta_1+\delta+1/\beta} x(\theta)^T
(M_{22i}(\theta)-\epsilon I)
x(\theta) d \theta\\
&\le -2 \delta +\epsilon_2
\int_{\theta_1-\delta-1/\beta}^{\theta_1-\delta} \norm{x(\theta)}^2
d \theta+\epsilon_2 \int_{\theta_1+\delta}^{\theta_1+\delta+1/\beta}
\norm{x(\theta)}^2 d \theta\\
&\le -2 \delta +2 \epsilon_2 \norm{x_1}^2/\beta \le -\delta\\
\end{align*}
Therefore, by contradiction, we have that $M_{22i}(\theta) \ge
\epsilon I$ for all $\theta \in [-\tau_i,-\tau_{i-1}]$ and
$i=1,\ldots, K$. Now we define $\epsilon'=\epsilon/2$ and
$\tilde{M}_{22i}(\theta)=M_{22i}(\theta)-\epsilon'I \ge \epsilon' I
$. We now show that $\tilde{M}_{22i}(\theta)^{-1}$ is continuous. We
first note that $\tilde{M}_{22i}(\theta)^{-1}$ is bounded, since
\begin{align*}
\tilde{M}_{22i}(\theta) &\ge \epsilon' I \\
&\Rightarrow I \ge \epsilon' \tilde{M}_{22i}(\theta)^{-1}\\
&\Rightarrow \tilde{M}_{22i}(\theta)^{-1} \le \frac{1}{\epsilon'} I.
\end{align*}
Now since $\tilde{M}_{22i}$ is continuous, for any $\beta > 0$,
there exists a $\delta>0$ such that
\[|\theta_1-\theta_2|\le \delta \Rightarrow
\norm{\tilde{M}_{22i}(\theta_1)-\tilde{M}_{22i}(\theta_2)}\le \beta.
\]
Therefore, for any $\beta'>0$, let $\delta$ be defined as above for
$\beta=\beta'\epsilon'^2$. Then $|\theta_1-\theta_2| \le \delta$
implies the following.
\begin{align*}
&\norm{\tilde{M}_{22i}(\theta_1)^{-1}-\tilde{M}_{22i}(\theta_2)^{-1}}\\
&=\norm{\tilde{M}_{22i}(\theta_1)^{-1} \tilde{M}_{22i}(\theta_2)
\tilde{M}_{22i}(\theta_2)^{-1}-\tilde{M}_{22i}(\theta_1)^{-1}
\tilde{M}_{22i}(\theta_1) \tilde{M}_{22i}(\theta_2)^{-1}}\\
&=\norm{\tilde{M}_{22i}(\theta_1)^{-1}
\left(\tilde{M}_{22i}(\theta_2)
-\tilde{M}_{22i}(\theta_1)\right) \tilde{M}_{22i}(\theta_2)^{-1}}\\
&\le \norm{\tilde{M}_{22i}(\theta_1)^{-1}}
\norm{\tilde{M}_{22i}(\theta_2)
-\tilde{M}_{22i}(\theta_1)}\norm{ \tilde{M}_{22i}(\theta_2)^{-1}}\\
&\le \frac{1}{\epsilon'^2} \norm{\tilde{M}_{22i}(\theta_2)
-\tilde{M}_{22i}(\theta_1)} \le \beta'
\end{align*}
Therefore $\tilde{M}_{22i}(\theta)^{-1}$ is continuous.

We now prove statement 2 by construction. Let

\[T_i(\theta)=T_0-(M_{11i}(\theta)-
M_{12i}(\theta)\tilde{M}_{22i}^{-1}(\theta)M_{12i}(\theta)^T),
\]
where
\[T_0=\frac{1}{\tau_K}\sum_{i=1}^K \int_{-\tau_i}^{-\tau_{i-1}}(M_{11i}(\theta)-
M_{12i}(\theta)\tilde{M}_{22i}^{-1}(\theta)M_{12i}(\theta)^T) d
\theta.
\]

Then we have that $T_i$ is continuous and
\[\sum_{i=1}^K \int_{-\tau_i}^{-\tau_{i-1}} T_i(\theta) d \theta=\tau_K T_0 -\tau_K T_0=0.\]
We now show that $T_0\ge 0$. For any constant vector $z_0$, let $z$
be a continuous function except possibly at points
$\{-\tau_i\}_{i=1}^{K-1}$ and such that
$z(0)=(I+\tilde{M}_{221}(0)^{-1} M_{121}(0)^T)z_0$, let
$x(\theta)=z(\theta)-\tilde{M}_{22i}(\theta)^{-1} M_{12i}(\theta)^T
z_0$ for $\theta\in[-\tau_i,-\tau_{i-1}]$. Then $x$ is continuous
except possibly at points $\{-\tau_i\}_{i=1}^{K-1}$, $x(0)=z_0$ and
by statement 1, and Lemma~\ref{lem:extproof3lem} we have the
following.

\begin{align*}
 &\sum_{i=1}^K \int_{-\tau_i}^{-\tau_{i-1}} \bmat{x(0)\\x(\theta)}^T
\bmat{M_{11i}(\theta) & M_{12i}(\theta)\\
M_{12i}(\theta)^T & \tilde{M}_{22i}(\theta)} \bmat{x(0)\\x(\theta)} d \theta\\
&=\sum_{i=1}^K \int_{-\tau_i}^{-\tau_{i-1}} \bmat{z_0\\z(\theta)}^T
\bmat{I & 0\\-\tilde{M}_{22i}(\theta)^{-1} M_{12i}(\theta)^T & I}^T
\bmat{M_{11i}(\theta) & M_{12i}(\theta)\\ M_{12i}(\theta)^T &
\tilde{M}_{22i}(\theta)} \\
&
\qquad\qquad\qquad\qquad\qquad\qquad\qquad\qquad\qquad\qquad\bmat{I
& 0\\-\tilde{M}_{22i}(\theta)^{-1} M_{12i}(\theta)^T & I}
\bmat{z_0\\z(\theta)} d
\theta \\
&=\sum_{i=1}^K \int_{-\tau_i}^{-\tau_{i-1}} \bmat{z_0\\z(\theta)}^T
\bmat{M_{11i}(\theta)-
M_{12i}(\theta)\tilde{M}_{22i}^{-1}(\theta)M_{12i}(\theta)^T & 0\\
0 & \tilde{M}_{22i}(\theta)} \bmat{z_0\\z(\theta)} d \theta \\
&=z_0^T\bbbl(\sum_{i=1}^K \int_{-\tau_i}^{-\tau_{i-1}}
\hspace{-5mm}M_{11i}(\theta)-
M_{12i}(\theta)\tilde{M}_{22i}^{-1}(\theta)M_{12i}(\theta)^T d
\theta\bbbr) z_0 \hspace{-1mm} + \hspace{-1mm}\sum_{i=1}^K
\int_{-\tau_i}^{-\tau_{i-1}}
\hspace{-5mm}z(\theta)^T \tilde{M}_{22i}(\theta) z(\theta) d \theta\\
&=z_0^T T_0 z_0 + \sum_{i=1}^K \int_{-\tau_i}^{-\tau_{i-1}}
z(\theta)^T \tilde{M}_{22i}(\theta) z(\theta) d \theta \ge \epsilon'
\norm{x}_2^2
\end{align*}

We now show that this implies that $T_0 \ge 0$. Suppose there exists
some $y$ such that $y^T T_0 y <0$. Then there exists some $z_0$ such
that $z_0 T_0 z_0=-1$. Now let $\alpha>1/\tau_1$ and
\[z(\theta)=\begin{cases}
(I+\tilde{M}_{221}(0)^{-1} M_{121}(0)^T)z_0(1+\alpha \theta)& \theta \in [-1/\alpha,0]\\
0 & \text{otherwise}
\end{cases}.
\]
Then $z$ is continuous,
$z(0)=(I+\tilde{M}_{221}(0)^{-1}M_{121}(0)^T)z_0$ and
$\norm{z(\theta)}^2\le\norm{z(0)}^2$ for all $\theta\in[-\tau,0]$.
Recall that $\tilde{M}_{22i}(\theta) \le
(\epsilon_2+\epsilon')I=\epsilon_3 I$. Let
$\alpha>2\epsilon_3\norm{z(0)^2}$ and then we have
\begin{align*}
z_0^T T_0 z_0 + \sum_{i=1}^K \int_{-\tau_i}^{-\tau_{i-1}} z(\theta)^T \tilde{M}_{22i}(\theta) z(\theta) d \theta &\le -1+\epsilon_3 \int_{-1/\alpha}^0 \norm{z(\theta)}^2 \\
&\le -1+\epsilon_3 \norm{z(0)}^2/\alpha < -\frac{1}{2}.
\end{align*}

But this is in contradiction to the previous relation. Thus we have
by contradiction that $T_0\ge0$. Now by using the invertibility of
$\tilde{M}_{22i}=M_{22i}(\theta)-\epsilon'I\ge \epsilon' I$ and the
Schur complement transformation, we have that statement 2 is
equivalent to the following for $i=1,\ldots, K$.
\begin{align*}
&\tilde{M}_{22i}(\theta)\ge 0 & \theta \in [-\tau_i,-\tau_{i-1}]\\
&M_{11i}(\theta)+T_i(\theta)-M_{12i}(\theta)\tilde{M}_{22i}(\theta)^{-1}
M_{12i}(\theta)^T \ge 0 \qquad & \theta \in [-\tau_i,-\tau_{i-1}]
\end{align*}
We have already shown the first statement. Finally, we have the
following for $i=1\ldots K$.
\begin{align*}
&M_{11i}(\theta)+T_i(\theta)-M_{12i}(\theta)\tilde{M}_{22i}(\theta)^{-1}
M_{12i}(\theta)^T=T_0 \ge 0
\end{align*}
Thus we have shown that statement 2 is true.
\end{proof}

%%%%%%%%%%%%%%%%%%%%%%%%%%%%%%%%%%%%%%%%%%%%%%%%%%%%%%%%%%%%%%%%%%%%%%
% Extension Proof 4
%%%%%%%%%%%%%%%%%%%%%%%%%%%%%%%%%%%%%%%%%%%%%%%%%%%%%%%%%%%%%%%%%%%%%%
\begin{lem}\label{lem:extproof4lem}
Suppose $\{M_i\}_{i=1}^K$ are continuous symmetric matrix valued
functions with domains $[-\tau_i,-\tau_{i-1}]$ and $\tau_K > \tau_i
> \tau_{i-1} >\tau_0=0$ for $i=2\ldots K-1$.
%which
%is continuous except possibly at points $\{\tau_i\}_{i=1}^K$ where
%$\tau_i \ge \tau_{i-1}$ and $\tau_i \in [-\tau,0]$.
Then the following are equivalent.

\begin{enumerate}
\item The following holds for all $x\in \mathcal{C}_{\tau_K}$.
\[\sum_{i=1}^{K} \int_{-\tau_i}^{-\tau_{i-1}}
\bmat{x(-\tau_0)\\\vdots\\x(-\tau_K)\\x(\theta)}^T M_i(\theta)
\bmat{x(-\tau_0)\\\vdots\\x(-\tau_K)\\x(\theta)}d\theta \ge 0\]
\item The following holds for all $x$ where $x$ is continuous except
possibly at points $\{-\tau_i\}_{i=1}^{K-1}$.
\[\sum_{i=1}^{K} \int_{-\tau_i}^{-\tau_{i-1}}
\bmat{x(-\tau_0)^-\\x(-\tau_1)^+\\\vdots\\x(-\tau_K)^+\\x(\theta)}^T
M_i(\theta)
\bmat{x(-\tau_0)^-\\x(-\tau_1)^+\\\vdots\\x(-\tau_K)^+\\x(\theta)}d\theta
\ge 0\] Where $x(\tau)^-$ denotes the left-handed limit and
$x(\tau)^+$ denotes the right-handed limit.
\end{enumerate}

\end{lem}

\begin{proof}
Clearly statement 2 implies statement 1. Now suppose statement 2 is
false. Then there exists some piecewise continuous $x$ , and
$\epsilon >0$ such that
\[V_1=\sum_{i=1}^{K} \int_{-\tau_i}^{-\tau_{i-1}} \bmat{x(-\tau_0)^-\\x(-\tau_1)^+\\
\vdots\\x(-\tau_K)^+\\x(\theta)}^T M_i(\theta)
\bmat{x(-\tau_0)^-\\x(-\tau_1)^+\\\vdots\\x(-\tau_K)^+\\x(\theta)}d\theta
\le -\epsilon.\]

Now let $\beta>\max_{i=1, \ldots, K}\{2/(\tau_i-\tau_{i-1})\}$ and
\[\hat{x}(\theta)=\begin{cases}
%x(-\tau_i-\frac{1}{\beta})+(x(-\tau_i+\frac{1}{\beta})
%-x(-\tau_i-\frac{1}{\beta}))(\beta(\theta+\tau_i+\frac{1}{\beta})/2)
%& \theta \in [-\tau_i-\frac{1}{\beta},-\tau_i+\frac{1}{\beta}] \\
x(-\tau_i-\frac{1}{\beta})+(x(-\tau_i)^+ -
x(-\tau_i-\frac{1}{\beta}) )\beta(\theta+\tau_i+\frac{1}{\beta}) &
\begin{matrix}\theta \in [-\tau_i-\frac{1}{\beta},-\tau_i],\\
i=1,\ldots, K-1 \end{matrix}\\x(\theta) &
\text{otherwise.}\end{cases}
\]
Then $\hat{x}$ is continuous on $[-\tau_K,0]$, $\hat{x}(0)=x(0)^-$
and $\hat{x}(-\tau_i)=x(-\tau_i)^+$ for $i=1,\ldots, K$. Since $x$
is piecewise continuous on $[-\tau_K,0]$, it is bounded on this
interval. Now let $c=\max_{\theta \in
[-\tau_{K},0]}{\norm{x(\theta)}^2}$. Then $\norm{\hat{x}(\theta)}^2
\le c$ for all $\theta \in [-\tau_k,0]$.
%and thus
%\begin{align*}\norm{\hat{x}}_2^2&=\int_{-\tau_K}^0\norm{\hat{x}(\theta)}^2d\theta
%=\norm{x}_2^2+\sum_{i=1}^K
%\int_{-\tau_i-1/\beta}^{-\tau_i+1/\beta}\left(\norm{\hat{x}(\theta)}^2
%-\norm{x(\theta)}^2\right)\\
%&\le\norm{x}_2^2 + 4 K c/\beta
%\end{align*}
%Then $\norm{\hat{x}}_2^2\le \norm{x}_2^2+4Kc/\beta \le 2
%\norm{x}_2^2$ for $\beta>4Kc/\norm{x}_2^2$.
Now since $M$ is piecewise continuous on $[-\tau_k,0]$, there exists
some $\epsilon'>0$ such that $M(\theta)<\epsilon' I$ for $\theta \in
[-\tau_k,0]$. Therefore we have the following.
\begin{align*}
&\sum_{i=1}^{K} \int_{-\tau_i}^{-\tau_{i-1}}
\bmat{\hat{x}(-\tau_0)\\\vdots\\\hat{x}(-\tau_K)\\\hat{x}(\theta)}^T
M_i(\theta) \bmat{\hat{x}(-\tau_0)\\\vdots\\\hat{x}(-\tau_K)\\\hat{x}(\theta)}d\theta\\
&=\sum_{i=1}^{K} \int_{-\tau_i}^{-\tau_{i-1}}
\bmat{x(-\tau_0)^-\\x(-\tau_1)^+\\\vdots\\x(-\tau_K)^+\\\hat{x}(\theta)}^T
M_i(\theta) \bmat{x(-\tau_0)^-\\x(-\tau_1)^+\\\vdots\\x(-\tau_K)^+\\\hat{x}(\theta)}d\theta\\
&=V_1+\sum_{i=1}^{K-1} \int_{-\tau_i-1/\beta}^{-\tau_i}
\bbbbl(\bmat{x(-\tau_0)^-\\x(-\tau_1)^+\\\vdots\\x(-\tau_K)^+\\\hat{x}(\theta)}^T
M_{i+1}(\theta) \bmat{x(-\tau_0)^-\\x(-\tau_1)^+\\
\vdots\\x(-\tau_K)^+\\\hat{x}(\theta)}
\\
&\qquad \qquad \qquad \qquad \qquad \qquad \qquad \qquad \qquad
-\bmat{x(-\tau_0)^-\\x(-\tau_1)^+\\\vdots\\x(-\tau_K)^+\\x(\theta)}^T
M_{i+1}(\theta)\bmat{x(-\tau_0)^-\\x(-\tau_1)^+\\\vdots\\x(-\tau_K)^+\\x(\theta)}\bbbbr) d \theta\\
&\le V_1 + \sum_{i=1}^{K-1}
\int_{-\tau_i-1/\beta}^{-\tau_i}\hspace{-2mm}
\epsilon'\left(\left(\norm{\hat{x}(\theta)}^2\hspace{-1mm}
+\hspace{-1mm}\sum_{j=0}^{K}\norm{\hat{x}(-\tau_j)}^2\right)
\hspace{-1mm}+\hspace{-1mm}\left(\norm{x(\theta)}^2+\hspace{-1mm}\sum_{j=0}^{K}\norm{\hat{x}(-\tau_j)}^2\right)\right)d\theta\\
&\le V_1 + \epsilon'(K-1)(2c+2(K+1)c)/\beta\\
&\le -\epsilon +
2\epsilon'(K-1)(K+2)c/\beta \\
&\le-\epsilon/2
\end{align*}
Which holds for $\beta \ge \frac{4\epsilon'(K-1)(K+2)c}{\epsilon}$.
Thus statement 2 is false implies statement 1 is false. Therefore
statement 1 implies statement 2.
\end{proof}

%%%%%%%%%%%%%%%%%%%%%%%%%%%%%%%%%%%%%%%%%%%%%%%%%%%%%%%%%%%%%%%%%%%%%

\begin{lem}
Suppose $S_i:\R \mapsto \S^{n(K+2)}$ are continuous matrix valued
functions with domains $[-\tau_i,-\tau_{i-1}]$ for $i=1,\ldots, K$
where $\tau_K>\tau_{i}>\tau_{i-1}>\tau_0=0$ for $i=2,\ldots ,K-1$.
Then the following are equivalent.
\begin{enumerate}
\item There exists an $\epsilon>0$ such that the following holds
for all $x\in \mathcal{C}_{\tau_K}$.
\[
\sum_{i=1}^K \int_{-\tau_i}^{-\tau_{i-1}}
\bmat{x(-\tau_0)\\\vdots\\x(-\tau_K)\\x(\theta)}^T S_i(\theta)
\bmat{x(-\tau_0)\\\vdots\\x(-\tau_K)\\x(\theta)}d \theta \ge \epsilon \norm{x}_2^2\\
\]
\item There exists an $\epsilon'>0$ and continuous matrix
valued functions $T_i:\R \mapsto \S^{n(K+1)}$ such that
\begin{align*}
&S_i(\theta)+\bmat{T_i(\theta) & 0\\0&-\epsilon'I} \ge 0 \qquad \theta \in [-\tau_i,-\tau_{i-1}],\;i=1 \ldots K\\
&\sum_{i=1}^K \int_{-\tau_i}^{-\tau_{i-1}}T_i(\theta)=0.
\end{align*}
\end{enumerate}
\end{lem}

\begin{proof}
$(2 \Rightarrow 1)$  Suppose there exist continuous symmetric matrix
valued functions, $T_i:\R \mapsto \S^{n(K+1)}$, such that
\begin{align*}
&S_i(\theta)+\bmat{T_i(\theta) & 0\\0&-\epsilon'I} \ge 0 \qquad i=1 \ldots K\\
&\sum_{i=1}^K \int_{-\tau_i}^{-\tau_{i-1}}T_i(\theta)=0.
\end{align*}

Then

\begin{align*}&\sum_{i=1}^K \int_{-\tau_i}^{-\tau_{i-1}}
\bmat{x(-\tau_0)\\\vdots\\x(-\tau_K)\\x(\theta)}^T S_i(\theta)
\bmat{x(-\tau_0)\\\vdots\\x(-\tau_K)\\x(\theta)}d \theta - \epsilon \norm{x}_2^2\\
&=\sum_{i=1}^K \int_{-\tau_i}^{-\tau_{i-1}}
\bmat{x(-\tau_0)\\\vdots\\x(-\tau_K)\\x(\theta)}^T
S_i(\theta)-\bmat{T_i(\theta) &0 \\0& -\epsilon' I}
\bmat{x(-\tau_0)\\\vdots\\x(-\tau_K)\\x(\theta)} d \theta \ge 0.
\end{align*}

$(1 \Rightarrow 2)$ Suppose that statement 1 holds for some $S_i$.
Write $S_i$ as
\[S_i(\theta)=\bmat{M_{11i}(\theta) & M_{12i}(\theta)\\M_{12i}(\theta)^T &
M_{22i}(\theta)},
\]

where $M_{22i}:\R \mapsto \S^n$. We first prove that
$M_{22i}(\theta)\ge \epsilon I$ for all $\theta \in
[-\tau_i,-\tau_{i-1}]$, $i=1,\ldots ,K$. By statement 1, we have
that

\[\sum_{i=1}^K \int_{-\tau_i}^{-\tau_{i-1}}
\bmat{x(-\tau_0)\\\vdots\\x(-\tau_K)\\x(\theta)}^T
\bmat{M_{11i}(\theta) & M_{12i}(\theta)\\M_{12i}(\theta)^T &
M_{22i}(\theta)-\epsilon I}
\bmat{x(-\tau_0)\\\vdots\\x(-\tau_K)\\x(\theta)} d \theta \ge 0.\\
\]
Now suppose that $M_{22i}(\theta)- \epsilon I$ is not positive
semidefinite for all $\theta \in [-\tau_i,-\tau_{i-1}]$. Then there
exists some $x_0\in \R^n$ and $\theta_1 \in [-\tau_i,-\tau_{i-1}]$
such that $x_0^T(M_{22i}(\theta_1)-\epsilon I)x_0 < 0 $. By
continuity of $M_{22i}$, if $\theta_1=-\tau_i$ or
$\theta_1=-\tau_{i-1}$, then there exists some $\theta_1'\in
(-\tau_i,-\tau_{i-1})$ such that $x_0^T(M_{22i}(\theta_1')-\epsilon
I)x_0 <0$. Thus assume $\theta_1 \in (-\tau_i,-\tau_{i-1})$. Now,
since $M_{22i}$ is continuous, there exists some $x_1$ and
$\delta>0$ where $\theta_1+\delta<-\tau_{i-1}$, $\theta_1-\delta
> -\tau_i$ and such that $x_1^T (M_{22i}(\theta)-\epsilon I)x_1 \le -1$ for
$\theta \in [\theta_1-\delta,\theta_1+\delta]$. Then for $\beta >
\max \{1/(-\tau_i-\theta_1-\delta),1/(\tau_i+\theta_1-\delta)\}$,
let
\[x(\theta)=\begin{cases}
\beta(\theta-(\theta_1-\delta-1/\beta))x_1 & \theta \in [\theta_1-\delta-1/\beta,\theta_1-\delta]\\
x_1 & \theta\in[\theta_1-\delta,\theta_1+\delta]\\
(1-\beta(\theta-(\theta_1+\delta)))x_1 & \theta \in [\theta_1+\delta,\theta_1+\delta+1/\beta]\\
%x_0(1+\alpha \theta)& \theta \in [-1/\alpha,0]\\
0 & \text{otherwise.}
\end{cases}
\]

Then $x \in \mathcal{C}_\tau$, $x(-\tau_i)=0$ for $i=0, \ldots, K$
and $\norm{x(\theta)}^2\le \norm{x_1}^2$ for all $\theta \in
[-\tau_K,0]$. Now, since every $M_{22i}$ is continuous, each are
bounded on $[-\tau_i,-\tau_{i-1}]$. Therefore, there exists some
$\epsilon_2> 0$ such that $M_{22i}(\theta)-\epsilon I\le \epsilon_2
I$ for $\theta \in [-\tau_i,-\tau_{i-1}]$, $i=1,\ldots, K$. Then let
$\beta \ge 2\epsilon_2 \norm{x_1}^2/\delta$ and we have the
following.
\begin{align*}
&\sum_{i=1}^K \int_{-\tau_i}^{-\tau_{i-1}}
\bmat{x(-\tau_0)\\\vdots\\x(-\tau_K)\\x(\theta)}^T S_i(\theta)
\bmat{x(-\tau_0)\\\vdots\\x(-\tau_K)\\x(\theta)} d \theta-\epsilon
\norm{x}^2 \\
&= \sum_{i=1}^K \int_{-\tau_i}^{-\tau_{i-1}} x(\theta)^T
(M_{22i}(\theta)-\epsilon I)
x(\theta) d \theta\\
&= \int_{\theta_1-\delta}^{\theta_1+\delta} x_1^T
(M_{22i}(\theta)-\epsilon I) x_1 d
\theta+\int_{\theta_1-\delta-1/\beta}^{\theta_1-\delta} x(\theta)^T
(M_{22i}(\theta)-\epsilon I) x(\theta) d
\theta\\
&+\int_{\theta_1+\delta}^{\theta_1+\delta+1/\beta} x(\theta)^T
(M_{22i}(\theta)-\epsilon I)
x(\theta) d \theta\\
&\le -2 \delta +\epsilon_2
\int_{\theta_1-\delta-1/\beta}^{\theta_1-\delta} \norm{x(\theta)}^2
d \theta+\epsilon_2 \int_{\theta_1+\delta}^{\theta_1+\delta+1/\beta}
\norm{x(\theta)}^2 d \theta\\
&\le -2 \delta +2 \epsilon_2 \norm{x_1}^2/\beta \le -\delta\\
\end{align*}
Therefore, by contradiction, we have that $M_{22i}(\theta) \ge
\epsilon I$ for all $\theta \in [-\tau_i,-\tau_{i-1}]$, $i=1,\ldots,
K$. Now we define $\epsilon'=\epsilon/2$ and
$\tilde{M}_{22i}(\theta)=M_{22i}(\theta)-\epsilon'I \ge \epsilon' I
$. We now show that $\tilde{M}_{22i}(\theta)^{-1}$ is continuous. We
first note that $\tilde{M}_{22i}(\theta)^{-1}$ is bounded, since
\begin{align*}
\tilde{M}_{22i}(\theta) &\ge \epsilon' I \\
&\Rightarrow I \ge \epsilon' \tilde{M}_{22i}(\theta)^{-1}\\
&\Rightarrow \tilde{M}_{22i}(\theta)^{-1} \le \frac{1}{\epsilon'} I.
\end{align*}
Now since $\tilde{M}_{22i}$ is continuous, for any $\beta > 0$,
there exists a $\delta>0$ such that
\[|\theta_1-\theta_2|\le \delta \Rightarrow
\norm{\tilde{M}_{22i}(\theta_1)-\tilde{M}_{22i}(\theta_2)}\le \beta.
\]
Therefore, for any $\beta'>0$, let $\delta$ be defined as above for
$\beta=\beta'\epsilon'^2$. Then $|\theta_1-\theta_2| \le \delta$
implies the following.
\begin{align*}
&\norm{\tilde{M}_{22i}(\theta_1)^{-1}-\tilde{M}_{22i}(\theta_2)^{-1}}\\
&=\norm{\tilde{M}_{22i}(\theta_1)^{-1} \tilde{M}_{22i}(\theta_2)
\tilde{M}_{22i}(\theta_2)^{-1}-\tilde{M}_{22i}(\theta_1)^{-1}
\tilde{M}_{22i}(\theta_1) \tilde{M}_{22i}(\theta_2)^{-1}}\\
&=\norm{\tilde{M}_{22i}(\theta_1)^{-1}
\left(\tilde{M}_{22i}(\theta_2)
-\tilde{M}_{22i}(\theta_1)\right) \tilde{M}_{22i}(\theta_2)^{-1}}\\
&\le \norm{\tilde{M}_{22i}(\theta_1)^{-1}}
\norm{\tilde{M}_{22i}(\theta_2)
-\tilde{M}_{22i}(\theta_1)}\norm{ \tilde{M}_{22i}(\theta_2)^{-1}}\\
&\le \frac{1}{\epsilon'^2} \norm{\tilde{M}_{22i}(\theta_2)
-\tilde{M}_{22i}(\theta_1)} \le \beta'
\end{align*}
Therefore $\tilde{M}_{22i}(\theta)^{-1}$ is continuous.

We now prove statement 2 by construction. Let

\[T_i(\theta)=T_0-(M_{11i}(\theta)-
M_{12i}(\theta)\tilde{M}_{22i}^{-1}(\theta)M_{12i}(\theta)^T)
\]
Where
\[T_0=\frac{1}{\tau_K}\sum_{i=1}^K \int_{-\tau_i}^{-\tau_{i-1}}(M_{11i}(\theta)-
M_{12i}(\theta)\tilde{M}_{22i}^{-1}(\theta)M_{12i}(\theta)^T) d
\theta
\] Then we have that $T_i$ is continuous and
\[\sum_{i=1}^K \int_{-\tau_i}^{-\tau_{i-1}} T_i(\theta) d \theta=\tau_K T_0-\tau_K T_0=0.\]
Now for any constant vector $z_c=\bmat{z_0 &\cdots z_K}^T$, suppose
that $z$ is a continuous function except possibly at points
$\{-\tau_i\}_{i=1}^{K-1}$ such that
\begin{align*}
z(-\tau_i)^+&=z_i+\tilde{M}_{22i}(-\tau_i)^{-1} M_{12i}(-\tau_i)^T z_c \qquad i=1,\ldots, K,\\
z(-\tau_0)^-&=z_0+\tilde{M}_{221}(-\tau_0)^{-1} M_{121}(-\tau_0)^T
z_c.
\end{align*}
Then, let $x(\theta)=z(\theta)-\tilde{M}_{22i}(\theta)^{-1}
M_{12i}(\theta)^T z_c$ for $\theta \in [-\tau_i,-\tau_{i-1}]$. Then
$x$ is continuous except possibly at points
$\{-\tau_i\}_{i=1}^{K-1}$,
$x(-\tau_0)^-=z(-\tau_0)^--\tilde{M}_{221}(-\tau_0)^{-1}M_{121}(-\tau_0)^T
z_c=z_0$,
$x(-\tau_i)^+=z(-\tau_i)^+-\tilde{M}_{22i}(-\tau_i)^{-1}M_{12i}(-\tau_i)^T
z_c=z_i$ for $i=1,\ldots, K$ and by statement 1 and
Lemma~\ref{lem:extproof4lem}, we have the following.

\begin{align*}
 &\sum_{i=1}^K \int_{-\tau_i}^{-\tau_{i-1}} \bmat{x(-\tau_0)^-\\x(-\tau_1)^+\\\vdots\\x(-\tau_K)^+\\x(\theta)}^T
\bmat{M_{11i}(\theta) & M_{12i}(\theta)\\
M_{12i}(\theta)^T & \tilde{M}_{22i}(\theta)} \bmat{x(-\tau_0)^-\\x(-\tau_1)^+\\\vdots\\x(-\tau_K)^+\\x(\theta)} d \theta\\
&=\sum_{i=1}^K \int_{-\tau_i}^{-\tau_{i-1}} \bmat{z_c\\z(\theta)}^T
\bmat{I & 0\\-\tilde{M}_{22i}(\theta)^{-1} M_{12i}(\theta)^T & I}^T
\bmat{M_{11i}(\theta) & M_{12i}(\theta)\\ M_{12i}(\theta)^T &
\tilde{M}_{22i}(\theta)}\\
& \qquad \qquad\qquad\qquad\qquad \qquad\qquad\qquad
\qquad\quad\bmat{I & 0\\-\tilde{M}_{22i}(\theta)^{-1}
M_{12i}(\theta)^T & I} \bmat{z_c\\z(\theta)} d
\theta \\
&=\sum_{i=1}^K \int_{-\tau_i}^{-\tau_{i-1}} \bmat{z_c\\z(\theta)}^T
\bmat{M_{11i}(\theta)-
M_{12i}(\theta)\tilde{M}_{22i}^{-1}(\theta)M_{12i}(\theta)^T & 0\\
0 & \tilde{M}_{22i}(\theta)} \bmat{z_c\\z(\theta)} d \theta \\
&=z_c^T(\sum_{i=1}^K \int_{-\tau_i}^{-\tau_{i-1}}M_{11i}(\theta)-
M_{12i}(\theta)\tilde{M}_{22i}^{-1}(\theta)M_{12i}(\theta)^T d
\theta) z_c\\
& + \sum_{i=1}^K \int_{-\tau_i}^{-\tau_{i-1}}
z(\theta)^T \tilde{M}_{22i}(\theta) z(\theta) d \theta\\
&=z_c^T T_0 z_c
+ \sum_{i=1}^K \int_{-\tau_i}^{-\tau_{i-1}}
z(\theta)^T \tilde{M}_{22i}(\theta) z(\theta) d \theta \ge \epsilon' \norm{x}_2^2\\
\end{align*}

We now show that this implies that $T_0 \ge 0$. Suppose there exists
some $y$ such that $y^T T_0 y <0$. Then there exists some $z_c$ such
that $z_c T_0 z_c=-1$. Now let
$\alpha>2/\max_i\norm{\tau_i-\tau_{i-1}}$ and
\[z(\theta)=\begin{cases}
(z_0+\tilde{M}_{221}(0)^{-1} M_{121}(0)^T z_c)(1+\alpha \theta)& \theta \in [-1/\alpha,0]\\
(z_i+\tilde{M}_{22i}(-\tau_i)^{-1} M_{12i}(-\tau_i)^T z_c)(1+\alpha \theta)& \theta \in [-\tau_i,-\tau_i+1/\alpha] \quad i=1\ldots K\\
0 & \text{otherwise.}
\end{cases}
\]
Then $z$ is continuous except at points $\{-\tau_i\}_{i=1}^{K-1}$,
$z(0)^-=z_0+\tilde{M}_{221}(0)^{-1}M_{121}(0)^T z_c$,
$z(-\tau_i)^+=z_i+\tilde{M}_{22i}(-\tau_i)^{-1} M_{12i}(-\tau_i)^T
z_c$ and $\norm{z(\theta)}^2\le c$ for all $\theta\in[-\tau,0]$
where
\[c=\max\left\{\max_{i=0\ldots K}\{\norm{z_i+\tilde{M}_{22i}(-\tau_i)^{-1}
M_{12i}(-\tau_i)^T z_c}^2\},\norm{z_0+\tilde{M}_{221}(0)^{-1}
M_{121}(0)^T z_c}\right\}.\]

Recall that $\tilde{M}_{22i}(\theta) \le
(\epsilon_2+\epsilon')I=\epsilon_3 I$. Let $\alpha>2\epsilon_3
c(K+1)$ and then we have
\begin{align*}
&z_0^T T_0 z_0 + \sum_{i=1}^K \int_{-\tau_i}^{-\tau_{i-1}} z(\theta)^T \tilde{M}_{22i}(\theta) z(\theta) d \theta \\
&=z_0^T T_0 z_0 + \int_{-1/\alpha}^0z(\theta)^T
\tilde{M}_{221}(\theta) z(\theta)d\theta
+ \sum_{i=1}^K \int_{-\tau_i}^{-\tau_{i}+1/\alpha} z(\theta)^T \tilde{M}_{22i}(\theta) z(\theta) d \theta \\
&\le -1+\epsilon_3 \int_{-1/\alpha}^0 \norm{z(\theta)}^2 + \sum_{i=1}^K \epsilon_3 \int_{-\tau_i}^{-\tau_{i}+1/\alpha} \norm{z(\theta)}^2 d \theta \\
&\le -1+\epsilon_3 c(K+1)/\alpha < -\frac{1}{2}.
\end{align*}

But this contradicts the previous relation. Thus we have by
contradiction that $T_0\ge0$. Now by using the invertibility of
$\tilde{M}_{22i}=M_{22i}(\theta)-\epsilon'I\ge \epsilon' I$ and the
Schur complement transformation, we have that statement 2 is
equivalent to the following for $i=1, \ldots, K$.
\begin{align*}
&\tilde{M}_{22i}(\theta)\ge 0 & \theta \in [-\tau_i,-\tau_{i-1}]\\
&M_{11i}(\theta)+T_i(\theta)-M_{12i}(\theta)\tilde{M}_{22i}(\theta)^{-1}
M_{12i}(\theta)^T \ge 0 \qquad & \theta \in [-\tau_i,-\tau_{i-1}]
\end{align*}
We have already proven the first statement. Finally, we have the
following.
\begin{align*}
&M_{11i}(\theta)+T_i(\theta)-M_{12i}(\theta)\tilde{M}_{22i}(\theta)^{-1}
M_{12i}(\theta)^T=T_0 \ge 0
\end{align*}
Thus we have shown that statement 2 is true.
\end{proof}
%%%%%%%%%%%%%%%%%%%%%%%%%%%%%%%%%%%%%%%%%%%%%%%%%%%%%%%%%%%%%%%%%%%%%%%
% Finite-Rank Proof
%%%%%%%%%%%%%%%%%%%%%%%%%%%%%%%%%%%%%%%%%%%%%%%%%%%%%%%%%%%%%%%%%%%%%%%

\begin{lem}\label{lem:FRcontinuous}
Suppose $M:\R^2 \mapsto \R^{n \times n}$ is a matrix valued function
which is continuous except possibly at points $\{\tau_i\}_{i=1}^K$
where $\tau_i\le \tau_{i-1}$ for $i=1,\ldots,K$ and $\tau_0=0$.
%Define
%\begin{align*}
%Y_1&:=\{\bmat{x_1&x_2&x_3}^T: x_3\in \mathcal{C}_{\tau}, \;
%x_1=x_3(0),\;
%x_2=x_3(-\tau)\}\\
%Y_1&:=\{\bmat{x_1&x_2&x_3}^T: \text{ $x_3$ is continuous except at
%points $\tau_i$}\}
%\end{align*}
Then the following are equivalent.
\begin{enumerate}
\item The following holds for all $x\in \mathcal{C}_{\tau_K}$.
\[\int_{-\tau_K}^0 \int_{-\tau_K}^0 x(\theta)^T M(\theta,\omega) x(\theta)d\theta d \omega \ge 0\]
\item The following holds for all $x$ where $x$ is continuous except
possibly at points $\{-\tau_i\}_{i=1}^{K-1}$.
\[\int_{-\tau_K}^0 \int_{-\tau_K}^0 x(\theta)^T M(\theta,\omega) x(\theta)d\theta d \omega \ge 0\]
\end{enumerate}
\end{lem}

\begin{proof}
Clearly statement 2 implies statement 1. Now suppose statement 2 is
false. Then there exists some piecewise continuous $x$ , $\epsilon
>0$ such that
\[V_1=\int_{-\tau_K}^0 \int_{-\tau_K}^0 x(\theta)^T M(\theta,\omega) x(\omega) d\theta d \omega \le
-\epsilon.
\]
Now let
\begin{align*}
&\hat{x}(\theta)\\
&=\begin{cases}x(-\tau_i-\frac{1}{\beta})+\beta(x(-\tau_i+\frac{1}{\beta})
-x(-\tau_i-\frac{1}{\beta}))(\theta+\tau_i+\frac{1}{\beta})/2 &
\begin{matrix}\theta \in
[-\tau_i-\frac{1}{\beta},-\tau_i+\frac{1}{\beta}]\\
i=1,\ldots, K-1\end{matrix}\\x(\theta) &
\text{otherwise.}\end{cases}
\end{align*}
Then $\hat{x}$ is continuous. Since $x$ is piecewise continuous on
$[-\tau_K,0]$ it is bounded on this interval. Now let
$c=\max_{\theta \in [-\tau_{K},0]}{\norm{x(\theta)}^2}$. Then
$\norm{\hat{x}(\theta)}^2 \le c$ for all $\theta \in [-\tau_K,0]$.
%and thus
%\begin{align*}\norm{\hat{x}}_2^2&=\int_{-\tau_K}^0\norm{\hat{x}(\theta)}^2d\theta
%=\norm{x}_2^2+\sum_{i=1}^K
%\int_{-\tau_i-1/\beta}^{-\tau_i+1/\beta}\left(\norm{\hat{x}(\theta)}^2
%-\norm{x(\theta)}^2\right)\\
%&\le\norm{x}_2^2 + 4 K c/\beta
%\end{align*}
%Then $\norm{\hat{x}}_2^2\le \norm{x}_2^2+4Kc/\beta \le 2
%\norm{x}_2^2$ for $\beta>4Kc/\norm{x}_2^2$.
Now since $M$ is piecewise continuous on $[-\tau_K,0]$, it is
bounded on this interval, i.e. there exists some $\epsilon'>0$ such
that $y^T M(\theta,\omega) z<\epsilon'
\max\{\norm{y}^2,\norm{z}^2\}$ for $\theta,\omega \in [-\tau_K,0]$.
Therefore we have the following.
\begin{align*}
&\int_{-\tau_K}^0 \int_{-\tau_K}^0 \hat{x}(\theta)^T
M(\theta,\omega) \hat{x}(\omega) d\theta d \omega\\
&=V_1+\sum_{i,j=1}^{K-1}
\int_{-\tau_i-1/\beta}^{-\tau_i+1/\beta}\int_{-\tau_j-1/\beta}^{-\tau_j
+1/\beta} \left(\hat{x}(\theta)^T M(\theta,\omega) \hat{x}(\omega)
-x(\theta)^T M(\theta,\omega) x(\omega)\right) d \theta d \omega\\
&\le V_1 + \sum_{i,j=1}^{K-1}
\int_{-\tau_i-1/\beta}^{-\tau_i+1/\beta}\int_{-\tau_j-1/\beta}^{-\tau_j
+1/\beta} (2 \epsilon'
c) d\theta d \omega\\
&= -\epsilon + 4\epsilon'(K-1)^2 c/\beta < -\epsilon/2\\
\end{align*}
Which holds for $\beta \ge \frac{8\epsilon'(K-1)^2 c}{\epsilon}$.
Thus statement 2 is false implies statement 1 is false. Therefore
statement 1 implies statement 2.
\end{proof}

\begin{lem} Let $M$ be a matrix valued function
$M: \R^2 \rightarrow \S^{n}$ which is discontinuous only at points
$\theta,\omega=-\tau_{i}$ for $i=1, \ldots, K-1$ where the $\tau_i$
are increasing and $\tau_0=0$. Then $M \in \tilde{H}_2^+$ if and
only if there exists some continuous matrix valued function $R: \R^2
\rightarrow \R^{nK \times nK}$ such that $R \in H_2^+$ and the
following holds where $I_i=[-\tau_i,-\tau_{i-1}]$,
$\Delta_i=\tau_i-\tau_{i-1}$.
\begin{align*}
&M(\theta,\omega)=M_{ij}(\theta,\omega) \qquad \text{for all }\;
\theta \in
I_i, \quad \omega \in I_j\\
&M_{ij}(\theta,\omega)=
R_{ij}\left(\frac{\tau_K}{\Delta_i}\theta+\tau_{i-1}\frac{\tau_K}{\Delta_i},
\frac{\tau_K}{\Delta_j}\omega+\tau_{j-1}\frac{\tau_K}{\Delta_j}\right)\\
&R(\theta,\omega)=\bmat{R_{11}(\theta,\omega) & \ldots&  R_{1K}(\theta,\omega) \\
\vdots & & \vdots\\
R_{K1}(\theta,\omega) &\ldots & R_{KK}(\theta,\omega)}
\end{align*}
\end{lem}

\begin{proof}
($\Leftarrow$) Define $\theta_i(\theta) =
\frac{\tau_K}{\Delta_i}\theta+\tau_{i-1}\frac{\tau_K}{\Delta_i}$ and
$\omega_i(\omega) =
\frac{\tau_K}{\Delta_i}\omega+\tau_{i-1}\frac{\tau_K}{\Delta_i}$.
Then $\theta_i(-\tau_i)=-\tau_K$, $\theta_i(-\tau_{i-1})=0$,
$\theta(\theta_i)=\frac{\Delta_i}{\tau_K}\theta_i-\tau_{i-1}$ and
$\frac{d \theta}{d \theta_i}=\frac{\Delta_i}{\tau_K}$. The same
relation holds between $\omega$ and $\omega_i$. Because $M$ is
piecewise continuous, $M_{ij}$ as defined above are continuous. Then
for any $x \in \mathcal{C}_{\tau_K}$, let
$x_i(\theta)=\frac{\Delta_i}{\tau_K}
x\left(\frac{\Delta_i}{\tau_K}\theta-\tau_{i-1}\right)$. Then $x_i
\in \mathcal{C}_{\tau_K}$ and we have the following.
\begin{align*}
&\int_{-\tau_K}^0 \int_{-\tau_K}^0 x(\theta)^T
M(\theta,\omega)x(\omega)\\
&=\sum_{i,j=1}^{K}\int_{I_i} \int_{I_j} x(\theta)^T
M_{ij}(\theta,\omega)x(\omega)d \theta d \omega\\
&=\sum_{i,j=1}^{K}\int_{-\tau_K}^{0} \int_{-\tau_K}^{0}
x_i(\theta_i)^T R_{ij}(\theta_i,\omega_j
)x_j(\omega_j)d \theta_i d \omega_j\\
&=\int_{-\tau_K}^0 \int_{-\tau_K}^0
\bmat{x_1(\theta)\\\vdots\\x_K(\theta)}^T R(\theta,\omega) \bmat{x_1(\omega)\\
\vdots \\x_K(\omega)}d \theta d \omega \ge 0
\end{align*}
Thus $R\in H_2^+$ implies $M \in H_2^+$.

($\Rightarrow$) Now suppose $M \in H_2^+$. Then for any $nK$
dimensional vector valued function in $\mathcal{C}_{\tau_K}$,
$x(\theta)=\bmat{x_1(\theta)^T&  \cdots & x_K(\theta)^T}^T$, define
$\tilde{x}(\theta)= \frac{\tau_K}{\Delta_i}
x_i\left(\frac{\tau_K}{\Delta_i}\theta
+\tau_{i-1}\frac{\tau_K}{\Delta_i}\right)$ for $\theta \in I_i$.
Then $\tilde{x}$ is continuous except possibly at points
$\{-\tau_i\}$ and we have the following.
\begin{align*}
&\int_{-\tau_K}^0 \int_{-\tau_K}^0 x(\theta)^T
R(\theta,\omega)x(\omega)d \theta d \omega\\
&=\sum_{i,j=1}^{K}\int_{-\tau_K}^{0} \int_{-\tau_K}^{0}
x_i(\theta_i)^T R_{ij}(\theta_i,\omega_j
)x_j(\omega_j)d \theta_i d \omega_j\\
&=\sum_{i,j=1}^{K}\int_{I_i} \int_{I_j} \tilde{x}(\theta)^T
M_{ij}(\theta,\omega)\tilde{x}(\omega)d \theta d \omega\\
&=\int_{-\tau_K}^0 \int_{-\tau_K}^0 \tilde{x}(\theta)^T
M(\theta,\omega)\tilde{x}(\omega) d\theta d \omega
\end{align*}
where $M$ is as defined above. Since $M$ is bounded and positive on
$x\in \mathcal{C}_{\tau_K}$, Lemma~\ref{lem:FRcontinuous} shows that
$M$ is also positive on the space of functions which are continuous
except possibly at points $\{-\tau_i\}_{i=1}^{K-1}$. Therefore, we
have that $R \in H_2^+$.
\end{proof}

\include{appendix_NonlinearCase}

% bibliography.tex should include either
% \bibliographystyle{...}
% \bibliography{mythesis}
% or some other way of doing the bibliography
\bibliographystyle{IEEEtranS}
\bibliography{thesis_mmpeet}

\end{document}